\documentclass{amsart}
\usepackage{amssymb}
\usepackage{amsmath}
\usepackage{amsfonts}

\newtheorem{theorem}{Theorem}
\theoremstyle{plain}

\newtheorem{corollary}{Corollary}

\newtheorem{definition}{Definition}

\newtheorem{remark}{Remark}

\theoremstyle{definition}
\newtheorem{example}{Example}
\numberwithin{equation}{section}
\input{tcilatex}
\usepackage{graphicx}
\usepackage{graphicx}
\usepackage{graphicx}

\begin{document}
\title[Random IFSs]{A Fractal Valued Random Iteration Algorithm and Fractal
Hierarchy}
\author{Michael Barnsley}
\address{335 Pennbrooke Trace, Duluth,\\
GA 30097, USA}
\email{Mbarnsley@aol.com}
\author{John Hutchinson}
\address{Department of Mathematics\\
Australian National University}
\email{john.hutchinson@anu.edu.au}
\author{\"{O}rjan Stenflo}
\address{Department of Mathematics, Stockholm University, SE-10691
Stockholm, Sweden}
\email{stenflo@math.su.se}
\date{ \today }

\begin{abstract}
We describe new families of random fractals, referred to as
\textquotedblleft $V$-variable\textquotedblright , which are intermediate
between the notions of deterministic and of standard random fractals. The
parameter $V$ describes the degree of \textquotedblleft
variability\textquotedblright : at each magnification level any $V$-variable
fractals has at most $V$ key \textquotedblleft forms\textquotedblright\ or
\textquotedblleft shapes\textquotedblright . $V$-variable\ random fractals
have the surprising property that they can be computed using a forward
process. More precisely, a version of the usual Random Iteration Algorithm,
operating on sets (or measures) rather than points, can be used to sample
each family. To present this theory, we review relevant results on fractals
(and fractal measures), both deterministic and random. Then our new results
are obtained by constructing an iterated function system (a super IFS) from
a collection of standard IFSs together with a corresponding set of
probabilities. The attractor of the super IFS is called a superfractal; it
is a collection of $V$-variable random fractals (sets or measures) together
with an associated probability distribution on this collection. When the
underlying space is for example $\mathbb{R}^{2}$, and the transformations
are computationally straightforward (such as affine transformations), the
superfractal can be sampled by means of the algorithm, which is highly
efficient in terms of memory usage. The algorithm is illustrated by some
computed examples. Some variants, special cases, generalizations of the
framework, and potential applications are mentioned.
\end{abstract}

\maketitle

\section{Introduction and Notation}

\subsection{Fractals and Random Fractals.}

A theory of deterministic fractal sets and measures, using a
\textquotedblleft backward\textquotedblright\ algorithm, was developed in
Hutchinson \cite{Hu81}. A different approach using a \textquotedblleft
forward\textquotedblright\ algorithm was developed in Barnsley and Demko 
\cite{BaDe}.

Falconer \cite{Fa86}, Graf \cite{Graf} and Mauldin and Williams \cite{MaWi}
randomized each step in the backward construction algorithm to obtain random
fractal sets. Arbeiter \cite{Arbe} introduced and studied random fractal
measures; see also Olsen \cite{Olse}. Hutchinson and R\"{u}chendorff \cite%
{HuRa} and \cite{HuRF} introduced new probabilistic techniques which allowed
one to consider more general classes of random fractals. For further
material see Z\"{a}hle \cite{Zahl}, Patzschke and Z\"{a}hle \cite{PaZa}, and
the references in all of these.

This paper begins with a review of material on deterministic and random
fractals generated by IFSs, and then introduces the class of $V$-variable
fractals which in a sense provides a link between deterministic and
\textquotedblleft standard\textquotedblright\ random fractals.

Deterministic fractal sets and measures are defined as the attractors of
certain iterated function systems (IFSs), as reviewed in Section \ref{secIFS}%
. Approximations in practical situations quite easily can be computed using
the associated random iteration algorithm. Random fractals are typically
harder to compute because one has to first calculate lots of fine random
detail at low levels, then one level at a time, build up the higher levels.

In this paper we restrict the class of random fractals to ones that we call
random $V$-variable fractals. Superfractals are sets of $V$-variable
fractals. They can be defined using a new type of IFS, in fact a
\textquotedblleft super\textquotedblright\ IFS\ made of a finite number $N$
of IFSs, and there is available a novel random iteration algorithm: each
iteration produces new sets, lying increasingly close to V-variable fractals
belonging to the superfractal, and moving ergodically around the
superfractal.

Superfractals appear to be a new class of geometrical object, their elements
lying somewhere between fractals generated by IFSs with finitely many maps,
which correspond to $V=N=1$, and realizations of the most generic class of
random fractals, where the local structure around each of two distinct
points are independent, corresponding to $V=\infty $. They seem to allow
geometric modelling of some natural objects, examples including
realistic-looking leaves, clouds, and textures; and good approximations can
be computed fast in elementary computer graphics examples. They are
fascinating to watch, one after another, on a computer screen, diverse, yet
ordered enough to suggest coherent natural phenomena and potential
applications.

Areas of potential applications include computer graphics and rapid
simulation of trajectories of stochastic processes \ The forward algorithm
also enables rapid computation of good approximations to random (including
\textquotedblleft fully\textquotedblright\ random) processes, where
previously there was no available efficient algorithm.

\subsection{An Example.\label{AnExample}}

Here we give an illustration of an application of the theory in this paper.
By means of this example we introduce informally V-variable fractals and
superfractals. We also explain why we think these objects are of special
interest and deserve attention.

We start with two pairs of contractive affine transformations, $%
\{f_{1}^{1},f_{2}^{1}\}$ and $\{f_{1}^{2},f_{2}^{2}\}$, where $%
f_{m}^{n}:\square $ $\rightarrow \square $ with $\square :=[0,1]\times
\lbrack 0,1]\subset \mathbb{R}^{2}$. We use two pairs of screens, where each
screen corresponds to a copy of $\square $ and represents for example a
computer monitor. We designate one pair of screens to be the Input Screens,
denoted by $(\square _{1},\square _{2})$. The other pair of screens is
designated to be the Output Screens, denoted by $(\square _{1^{\prime
}},\square _{2^{\prime }})$.

Initialize by placing an image on each of the Input Screens, as illustrated
in Figure \ref{twofish1}, and clearing both of the Output Screens. We
construct an image on each of the two Output Screens as follows.

(i) Pick randomly one of the pairs of functions $\{f_{1}^{1},f_{2}^{1}\}$ or 
$\{f_{1}^{2},f_{2}^{2}\}$, say$\{f_{1}^{n_{1}},f_{2}^{n_{1}}\}$. Apply $%
f_{1}^{n_{1}}$ to one of the images on $\square _{1}$ or $\square _{2}$,
selected randomly, to make an image on $\square _{1^{\prime }}$. Then apply $%
f_{2}^{n_{1}}$ to one of the images on $\square _{1}$ or $\square _{2}$,
also selected randomly, and overlay the resulting image $I$ on the image now
already on $\square _{1^{\prime }}$. (For example, if black-and-white images
are used, simply take the union of the black region of $I$ with the black
region on $\square _{1^{\prime }}$, and put the result back onto $\square
_{1^{\prime }}$.)

(ii) Again pick randomly one of the pairs of functions $%
\{f_{1}^{1},f_{2}^{1}\}$ or $\{f_{1}^{2},f_{2}^{2}\}$, say $%
\{f_{1}^{n_{2}},f_{2}^{n_{2}}\}$. Apply $f_{1}^{n_{2}}$ to one of the images
on $\square _{1}$, or $\square _{2}$, selected randomly, to make an image on 
$\square _{2^{\prime }}$. Also apply $f_{2}^{n_{2}}$ to one of the images on 
$\square _{1}$, or $\square _{2}$, also selected randomly, and overlay the
resulting image on the image now already on $\square _{2^{\prime }}$.

(iii) Switch Input and Output, clear the new Output Screens, and repeat
steps (i), and (ii).

(iv) Repeat step (iii) many times, to allow the system to settle into its
\textquotedblleft stationary state\textquotedblright .

What kinds of images do we see on the successive pairs of screens, and what
are they like in the \textquotedblleft stationary state\textquotedblright ?
What does the theory developed in this paper tell us about such situations?

As a specific example, let us choose 
\begin{equation}
f_{1}^{1}(x,y)=(\frac{1}{2}x-\frac{3}{8}y+\frac{5}{16},\frac{1}{2}x+\frac{3}{%
8}y+\frac{3}{16})\text{,}  \label{trans_one}
\end{equation}%
\begin{equation}
f_{2}^{1}(x,y)=(\frac{1}{2}x+\frac{3}{8}y+\frac{3}{16},-\frac{1}{2}x+\frac{3%
}{8}y+\frac{11}{16})\text{,}  \label{trans_two}
\end{equation}%
\begin{equation*}
f_{1}^{2}(x,y)=(\frac{1}{2}x-\frac{3}{8}y+\frac{5}{16},-\frac{1}{2}x-\frac{3%
}{8}y+\frac{13}{16})\text{,}
\end{equation*}%
\begin{equation}
f_{2}^{2}(x,y)=(\frac{1}{2}x+\frac{3}{8}y+\frac{3}{16},\frac{1}{2}x-\frac{3}{%
8}y+\frac{5}{16})\text{.}  \label{trans_four}
\end{equation}

We describe how these transformations act on the triangle $ABC$ in the
diamond $ABCD$, where $A=(\frac{1}{4},\frac{1}{2})$, $B=(\frac{1}{2},\frac{3%
}{4})$, $C=(\frac{3}{4},\frac{1}{2})$, and $D=(\frac{1}{2},\frac{1}{4})$.
Let $B_{1}=$ $(\frac{9}{32},\frac{23}{32})$, $B_{2}=(\frac{23}{32},\frac{23}{%
32})$, $B_{3}=(\frac{9}{32},\frac{9}{32})$, and $B_{4}=(\frac{23}{32},\frac{%
23}{32})$. See Figure \ref{FigureA}. \FRAME{ftbpFU}{4.0465in}{3.0372in}{0pt}{%
\Qcb{Triangles used to define the four transformations $%
f_{1}^{1},f_{2}^{1},f_{1}^{2},$ and $f_{2}^{2}.$}}{\Qlb{FigureA}}{figurea.ps%
}{\raisebox{-3.0372in}{\includegraphics[height=3.0372in]{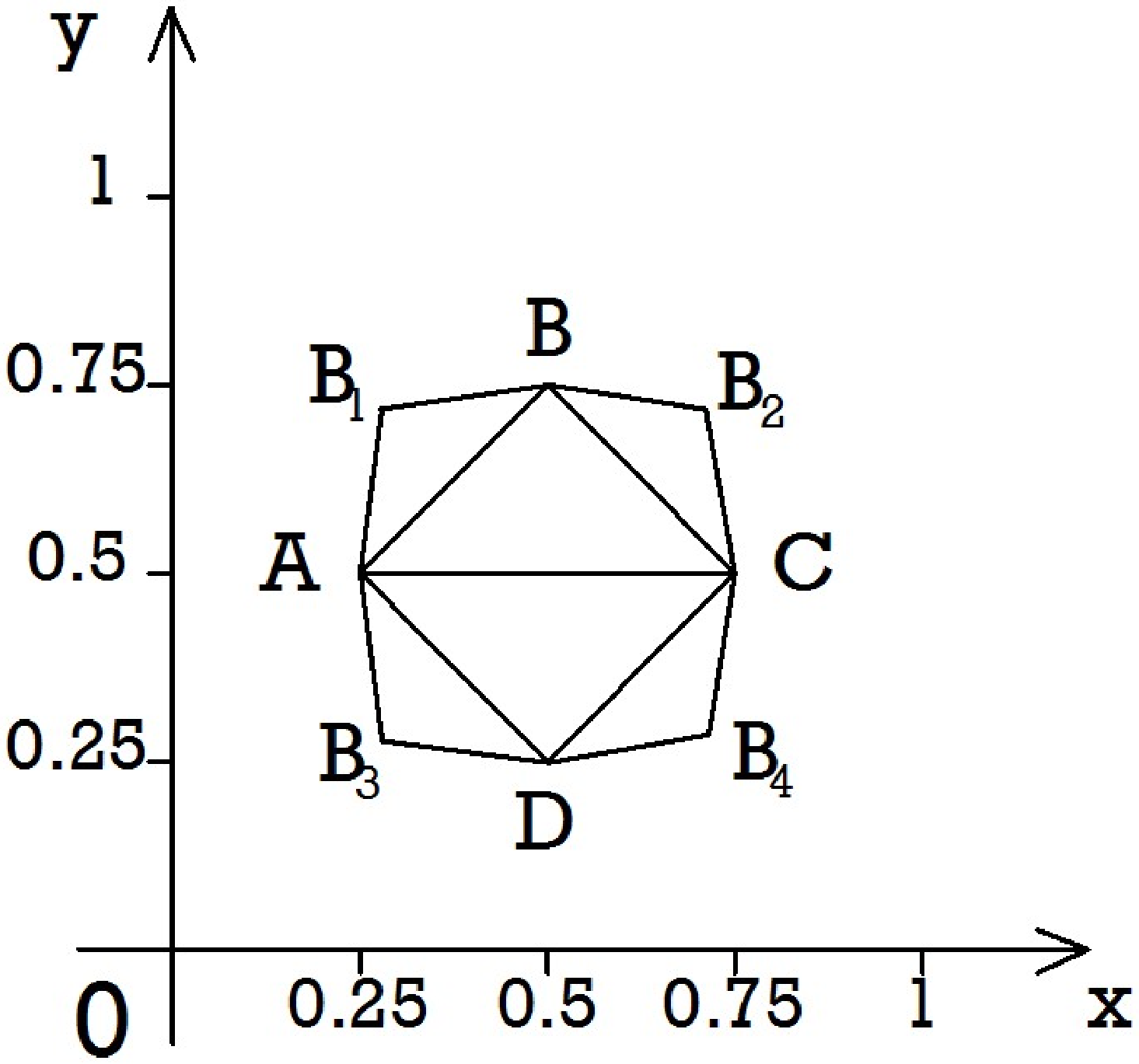}}} Then we have 
\begin{equation*}
f_{1}^{1}(A)=A\text{, }f_{1}^{1}(B)=B_{1}\text{, }f_{1}^{1}(C)=B\text{;}
\end{equation*}%
\begin{equation*}
f_{2}^{1}(A)=B\text{, }f_{2}^{1}(B)=B_{2}\text{, }f_{2}^{1}(C)=C\text{;}
\end{equation*}%
\begin{equation*}
f_{1}^{2}(A)=A\text{, }f_{1}^{2}(B)=B_{3}\text{, }f_{1}^{2}(C)=D\text{;}
\end{equation*}%
\begin{equation*}
f_{2}^{2}(A)=D\text{, }f_{2}^{2}(B)=B_{4}\text{, }f_{2}^{2}(C)=C\text{.}
\end{equation*}

In Figure \ref{twofish1} we show an initial pair of images, two jumping
fish, one on each of the two screens $\square _{1}$ and $\square _{2}$. In
Figures \ref{twofish2}, \ref{twofish3}, \ref{twofish4}, \ref{twofish5}, \ref%
{twofish6}, \ref{twofish7}, and \ref{twofish8}, we show the start of the
sequence of pairs of images obtained in a particular trial, for the first
seven iterations. Then in Figures \ref{final1}, \ref{final2}, and \ref%
{final3}, we show three successive pairs of computed screens, obtained after
more than twenty iterations. These latter images are typical of those
obtained after twenty or more iterations, very diverse, but always
representing continuous \textquotedblleft random\textquotedblright\ paths in 
$\mathbb{R}^{2}$; they correspond to the $``$stationary
state\textquotedblright , at the resolution of the images. More precisely,
with probability one the empirically obtained distribution on such images
over a long experimental run corresponds to the stationary state
distribution.\FRAME{ftbpFU}{4.2947in}{1.9752in}{0pt}{\Qcb{An initial image
of a jumping fish on each of the two screens $\square _{1}$ and $\square _{2}
$. }}{\Qlb{twofish1}}{twofish1.ps}{\raisebox{-1.9752in}{\includegraphics[height=1.9752in]{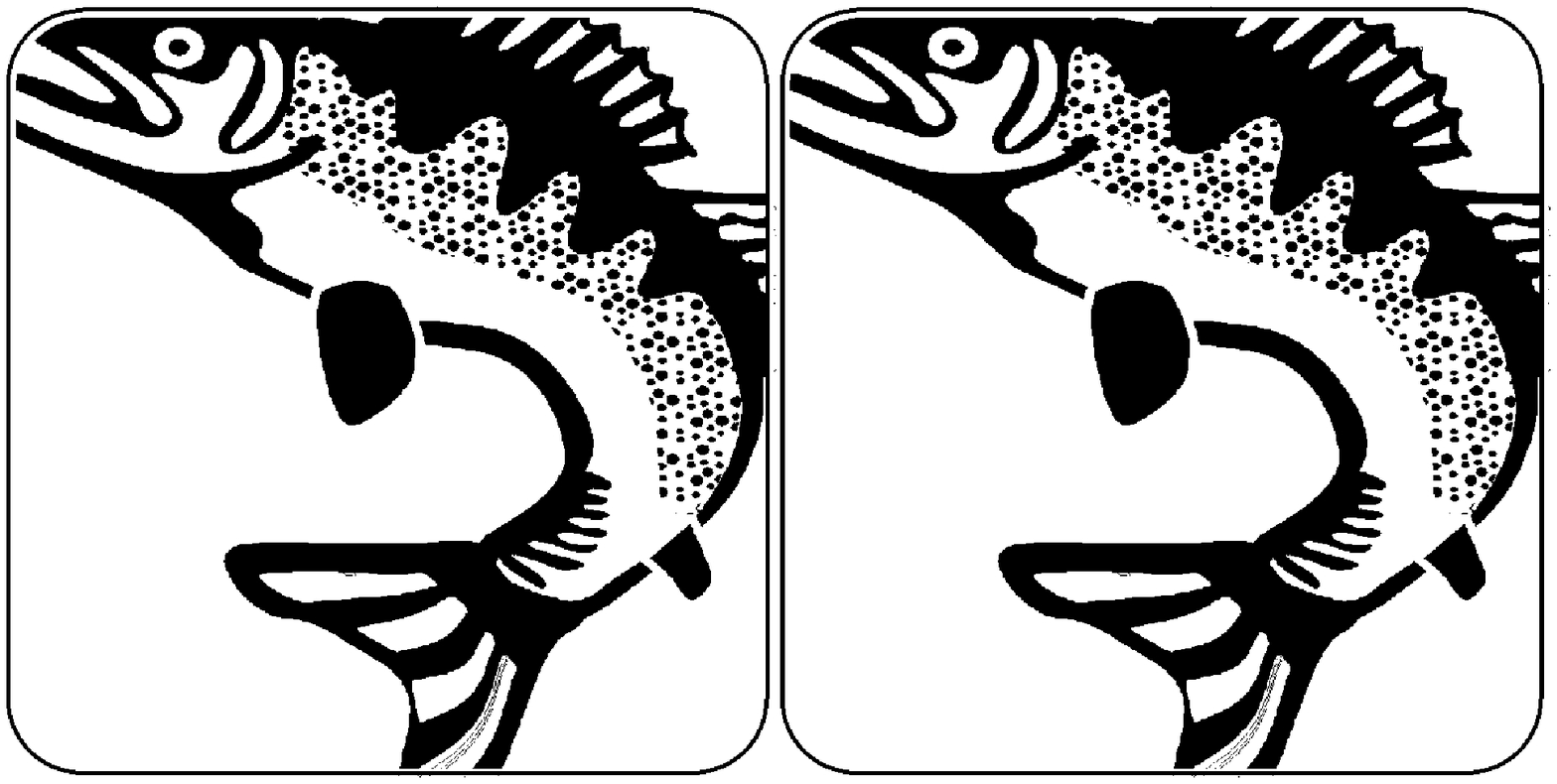}}}\FRAME{ftbpFU}{4.2947in}{1.9752in}{%
0pt}{\Qcb{The pair of images after one iteration.}}{\Qlb{twofish2}}{%
twofish2.ps}{\raisebox{-1.9752in}{\includegraphics[height=1.9752in]{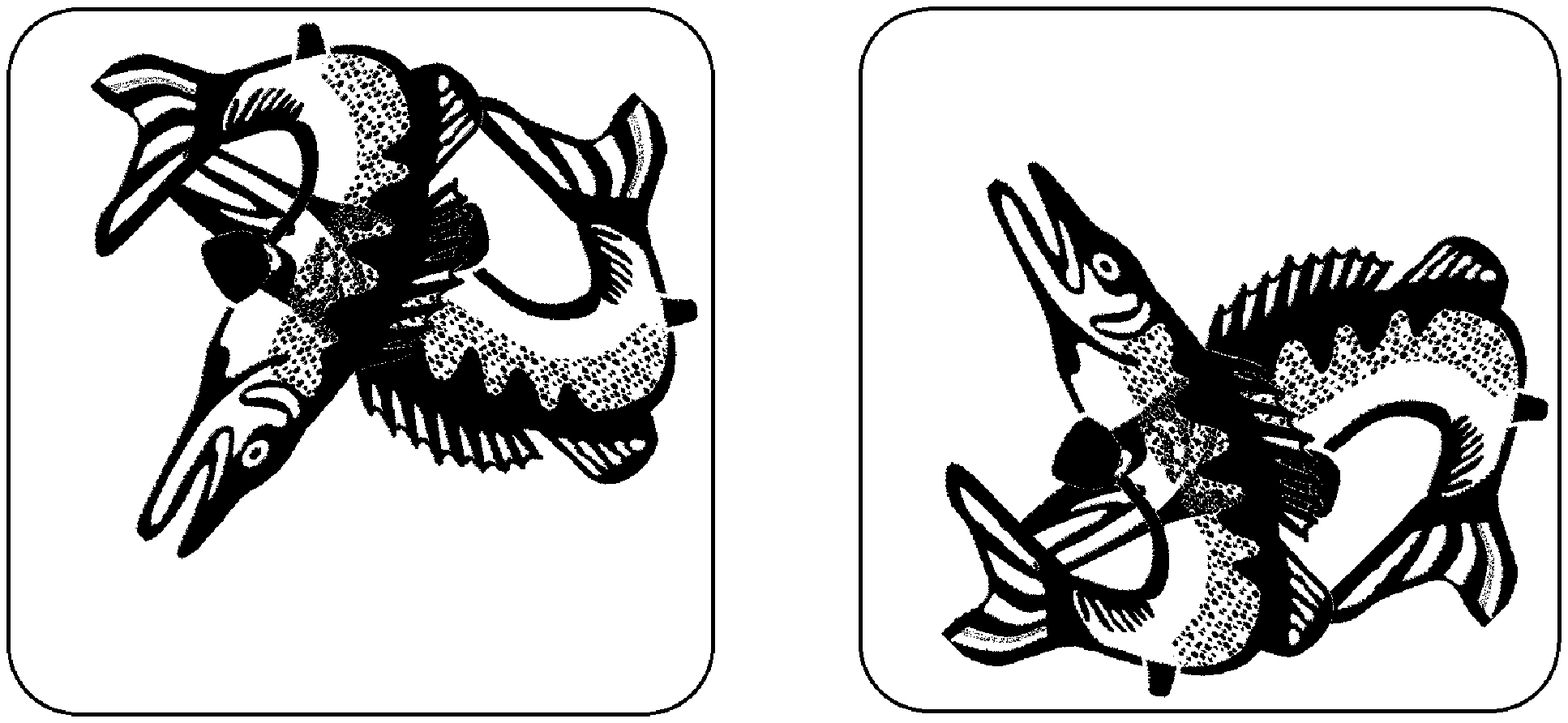}}}\FRAME{ftbpFU}{4.2947in}{1.9752in}{%
0pt}{\Qcb{The two images after two iterations.}}{\Qlb{twofish3}}{twofish3.ps%
}{\raisebox{-1.9752in}{\includegraphics[height=1.9752in]{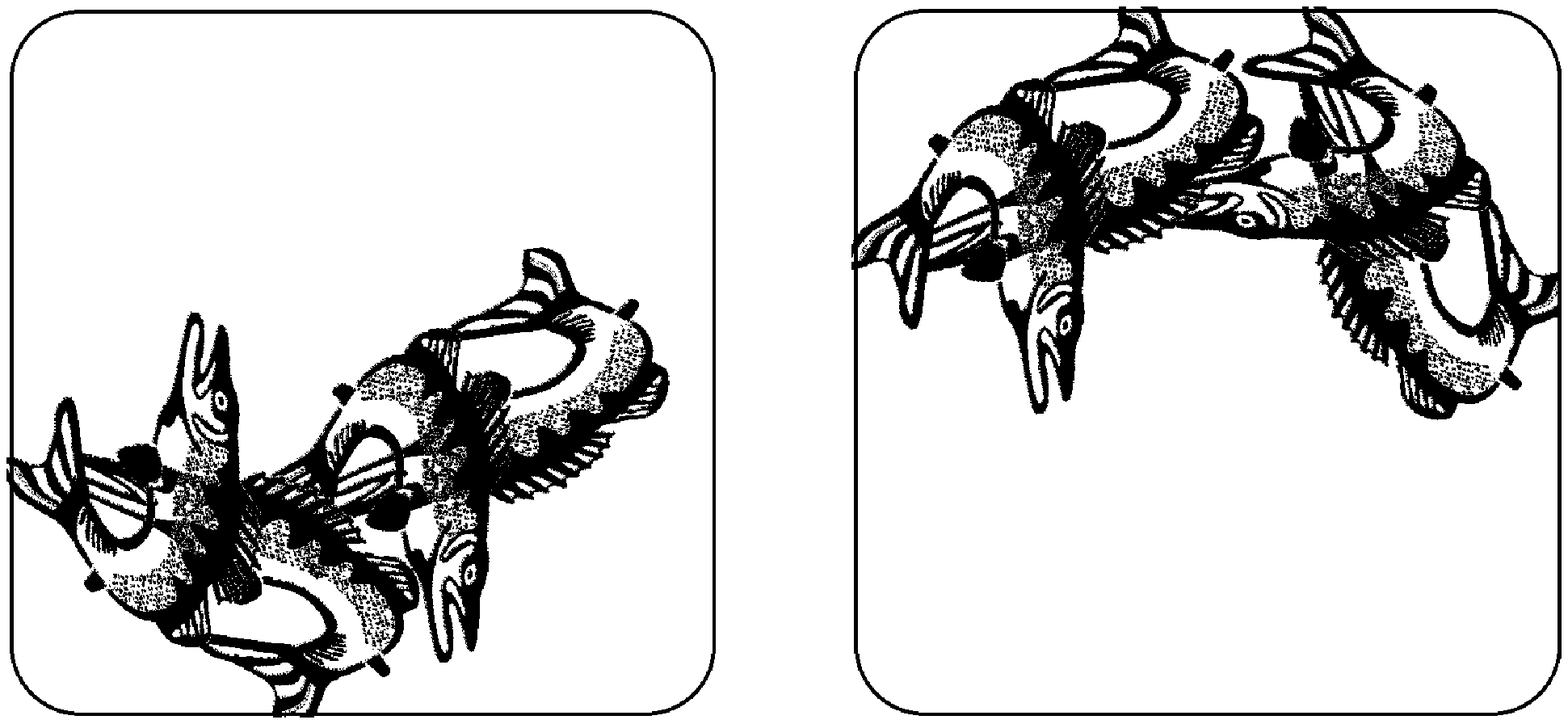}}}\FRAME{ftbpFU}{4.2947in}{1.9752in}{0pt}{\Qcb{The two images after
three iterations. Both images are the same.}}{\Qlb{twofish4}}{twofish4.ps}{%
\raisebox{-1.9752in}{\includegraphics[height=1.9752in]{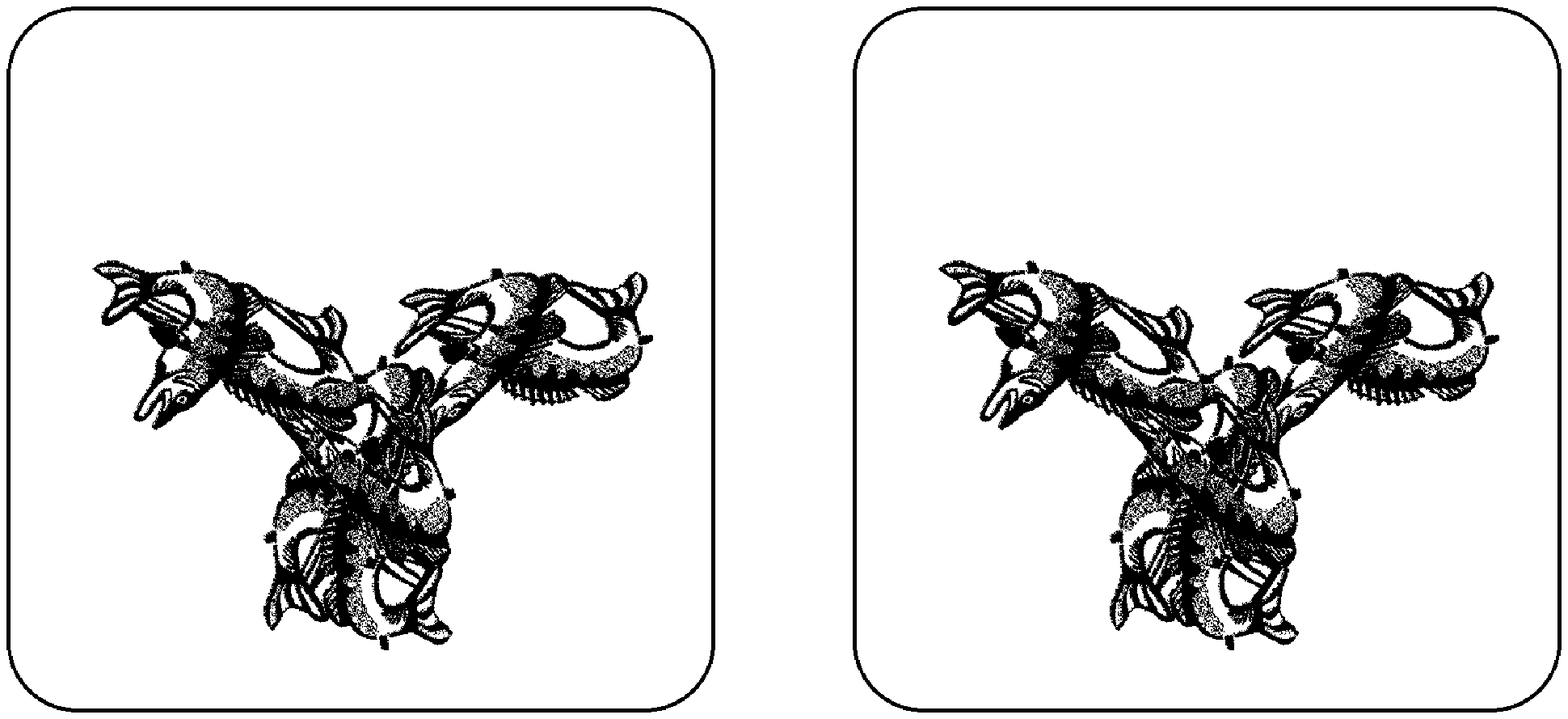}}}\FRAME{ftbpFU}{4.2947in}{1.9752in}{0pt}{\Qcb{The two images after
four iterations. Both images are again the same, a braid of fish.}}{\Qlb{%
twofish5}}{twofish5.ps}{\raisebox{-1.9752in}{\includegraphics[height=1.9752in]{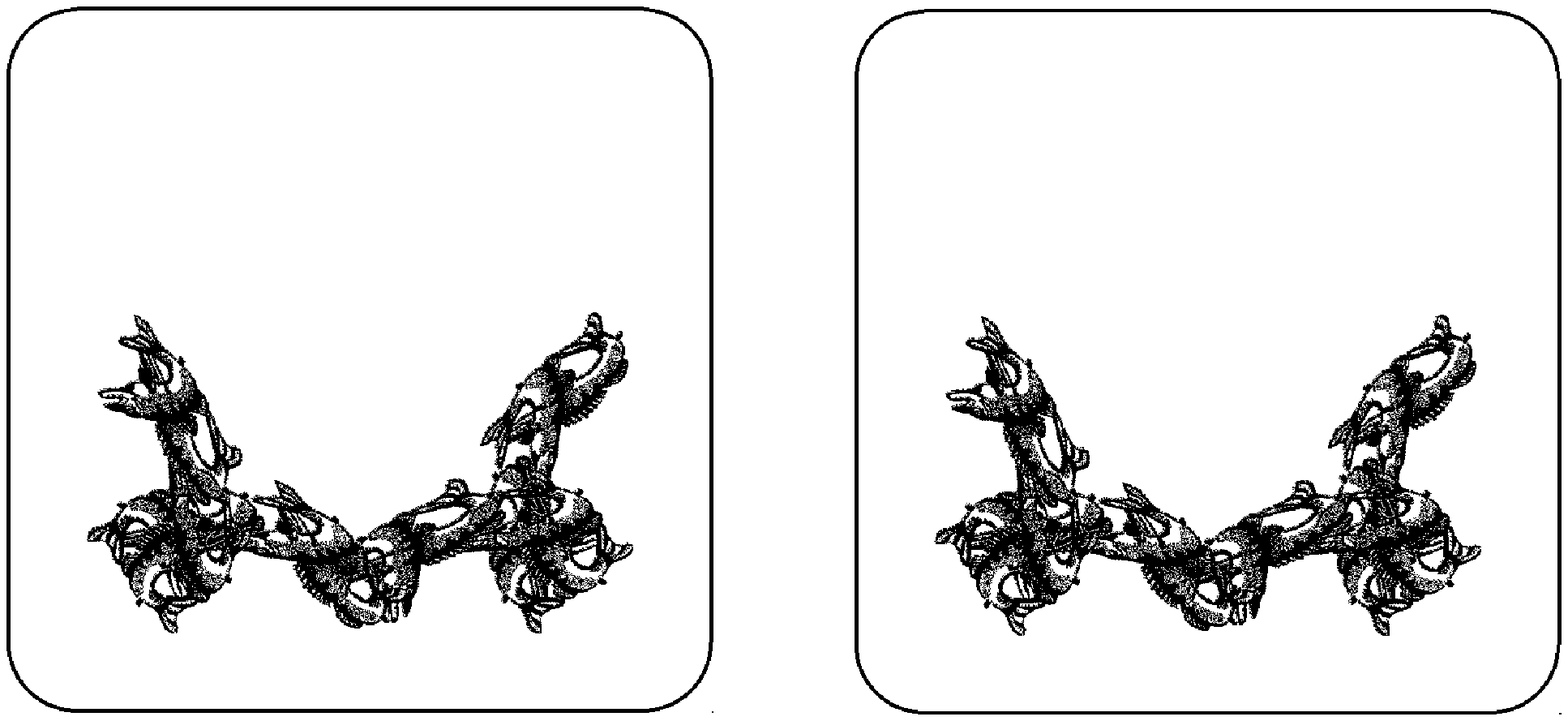}}}\FRAME{ftbpFU}{4.2947in}{1.9752in}{%
0pt}{\Qcb{The two images after five iterations. The two images are the same.}%
}{\Qlb{twofish6}}{twofish6.ps}{\raisebox{-1.9752in}{\includegraphics[height=1.9752in]{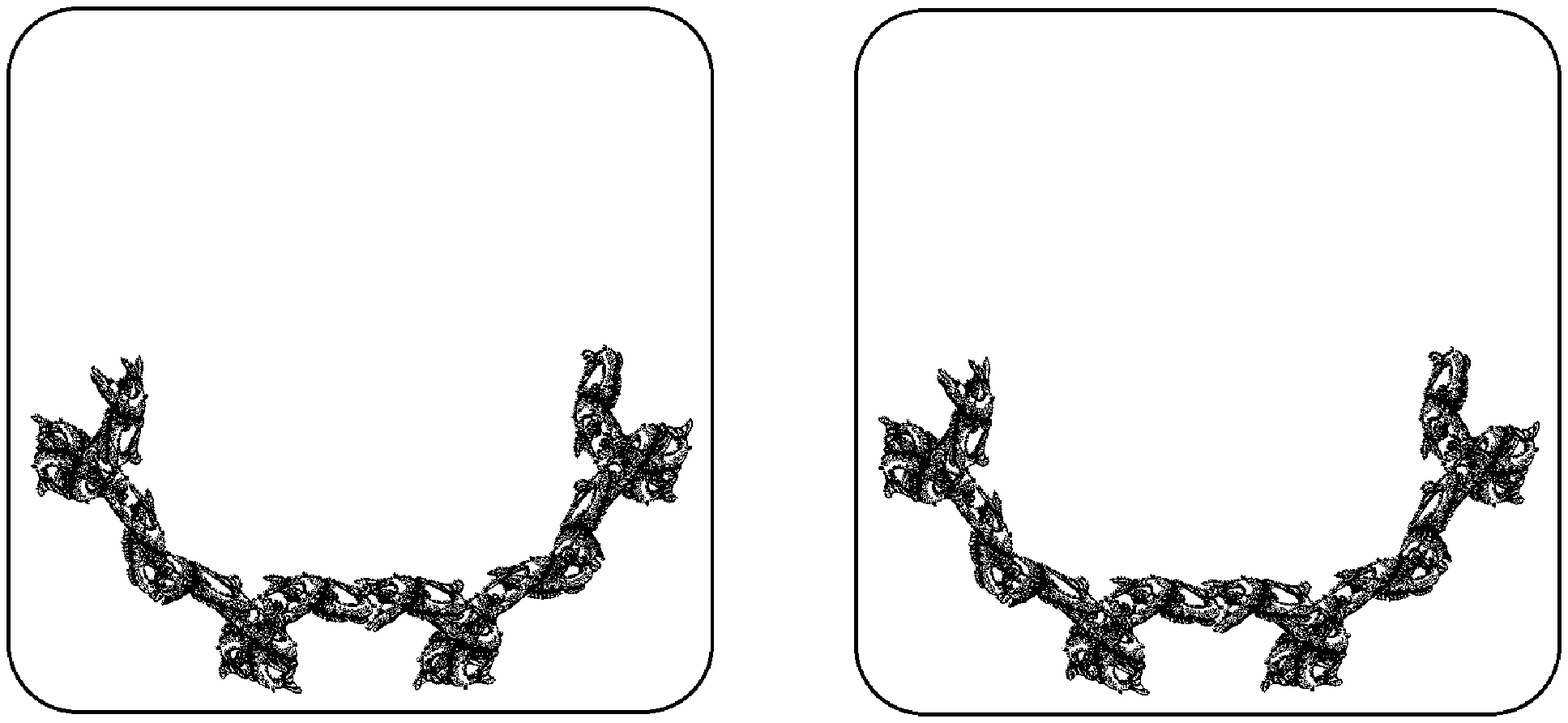}}}\FRAME{ftbpFU}{4.2947in}{1.9752in}{%
0pt}{\Qcb{The two images after six iterations.}}{\Qlb{twofish7}}{twofish7.ps%
}{\raisebox{-1.9752in}{\includegraphics[height=1.9752in]{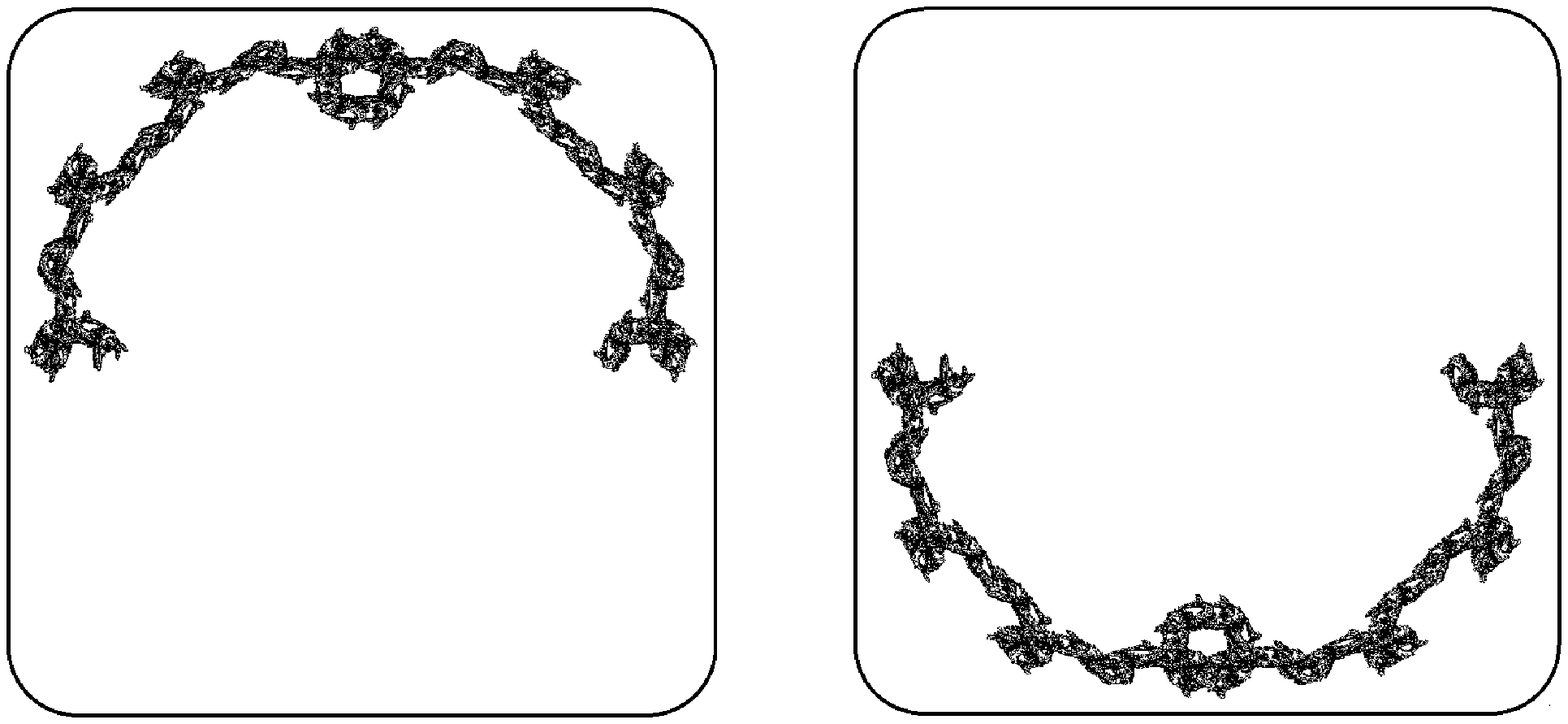}}}\FRAME{ftbpFU}{4.2947in}{1.9752in}{0pt}{\Qcb{The two images after
seven iterations.}}{\Qlb{twofish8}}{twofish8.ps}{\raisebox{-1.9752in}{\includegraphics[height=1.9752in]{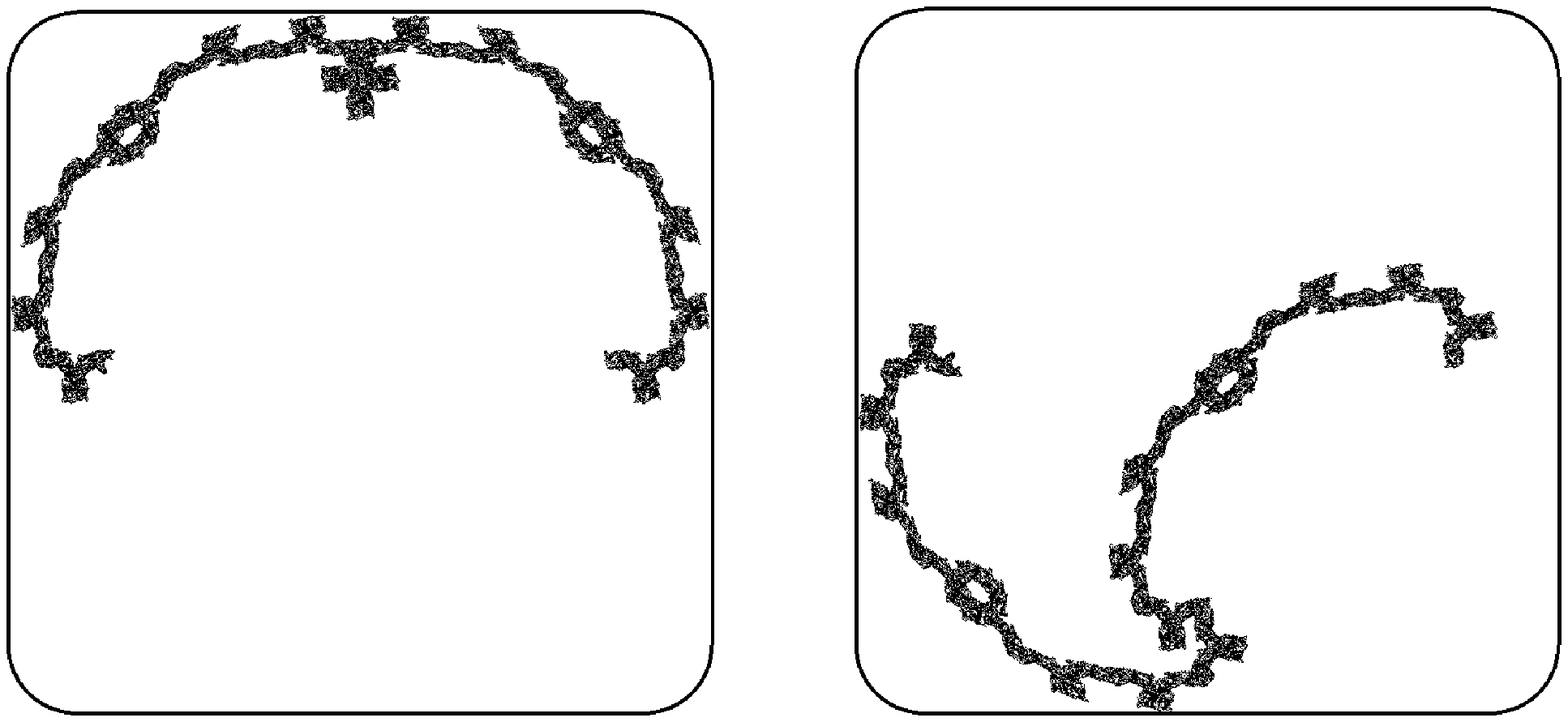}}}\FRAME{ftbpFU}{4.2947in}{1.9752in}{0pt}{\Qcb{Images on the two
screens $\square _{1}$ and $\square _{2}$ after a certain number $L>20$ of
iterations. Such pictures are typical of the "stationary state" at the
printed resolution.}}{\Qlb{final1}}{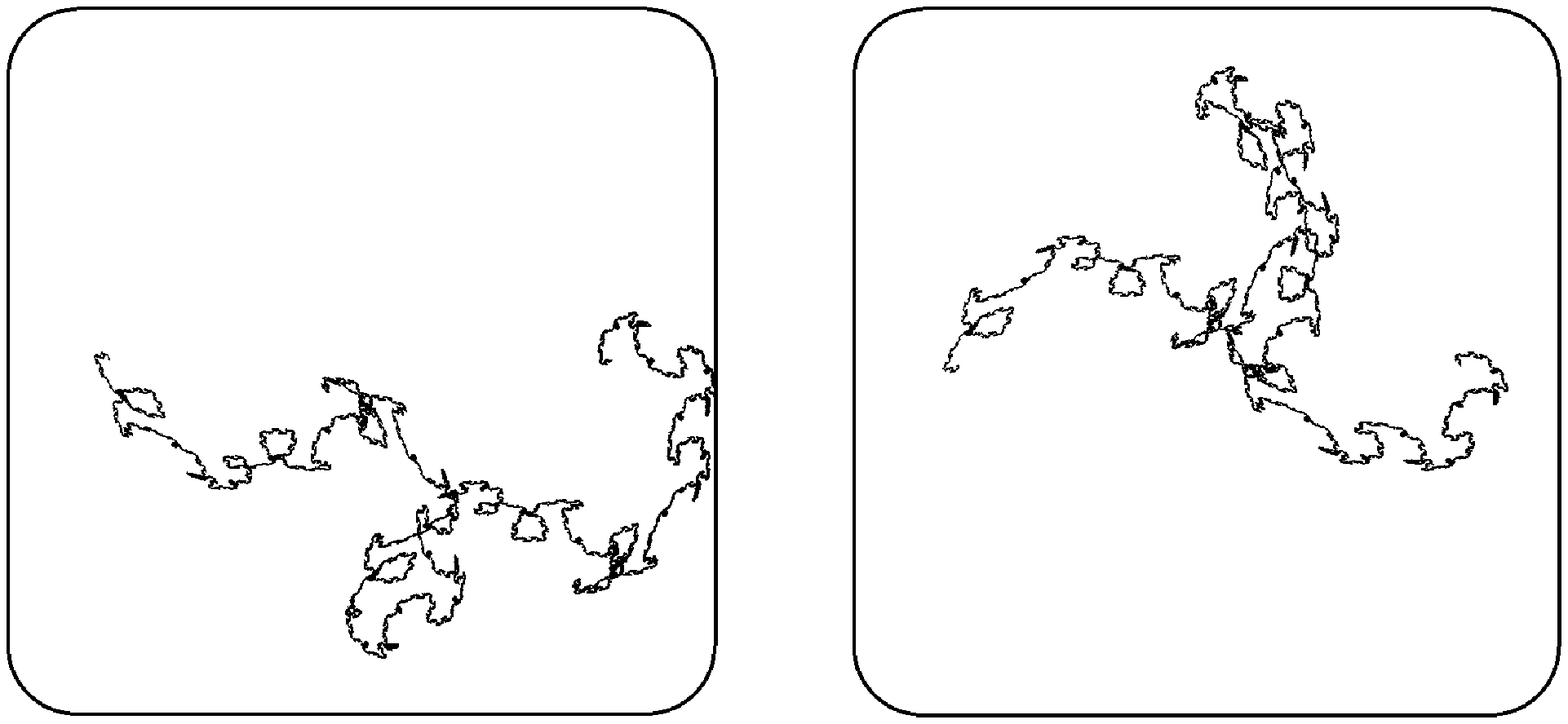}{\raisebox{-1.9752in}{\includegraphics[height=1.9752in]{final1.ps}}}\FRAME{ftbpFU}{4.2947in}{1.9752in}{0pt}{\Qcb{Images on the two
screens after $L+1$ iterations.}}{\Qlb{final2}}{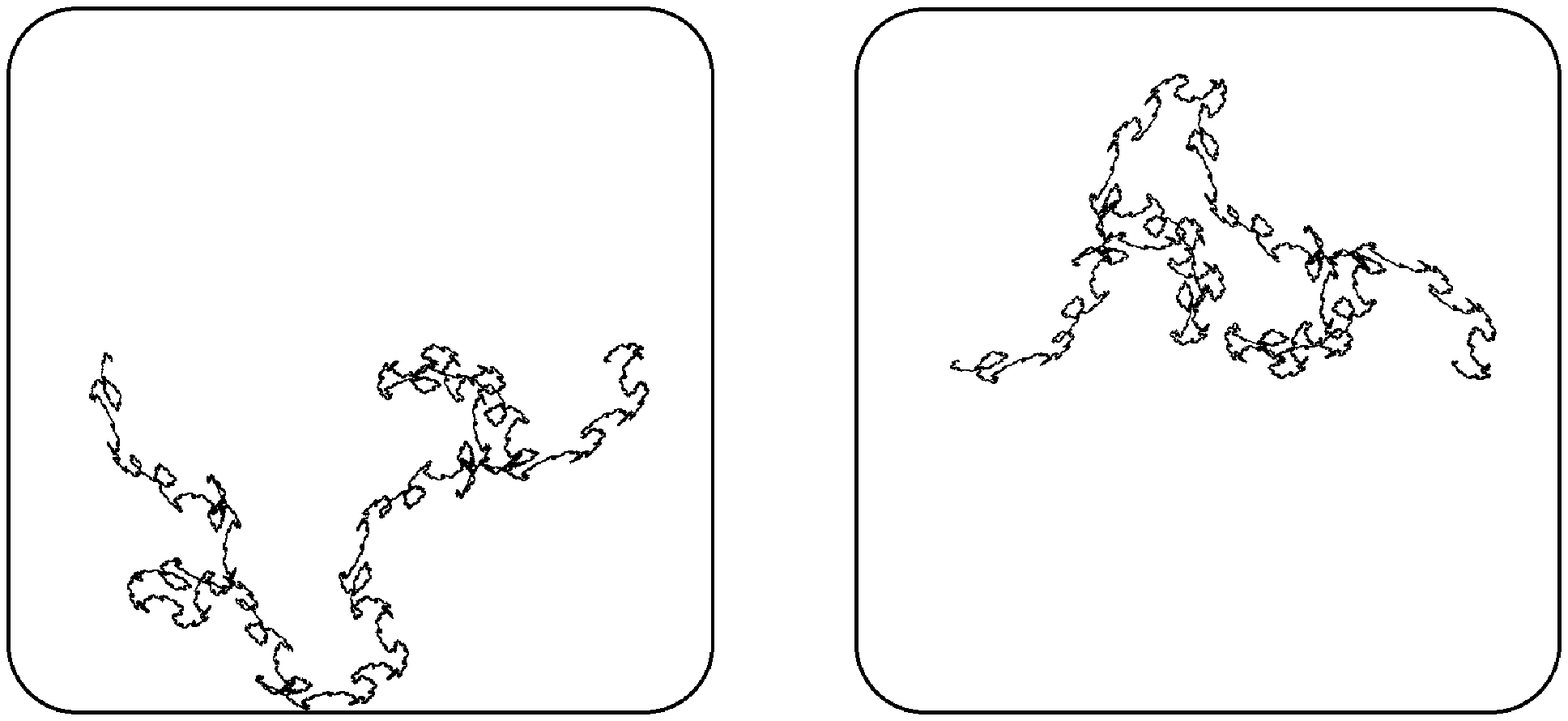}{\raisebox{-1.9752in}{\includegraphics[height=1.9752in]{final2.ps}}}\FRAME{ftbpFU}{4.2947in}{1.9752in}{0pt}{\Qcb{Images on the two
screens after $L+2$ iterations.}}{\Qlb{final3}}{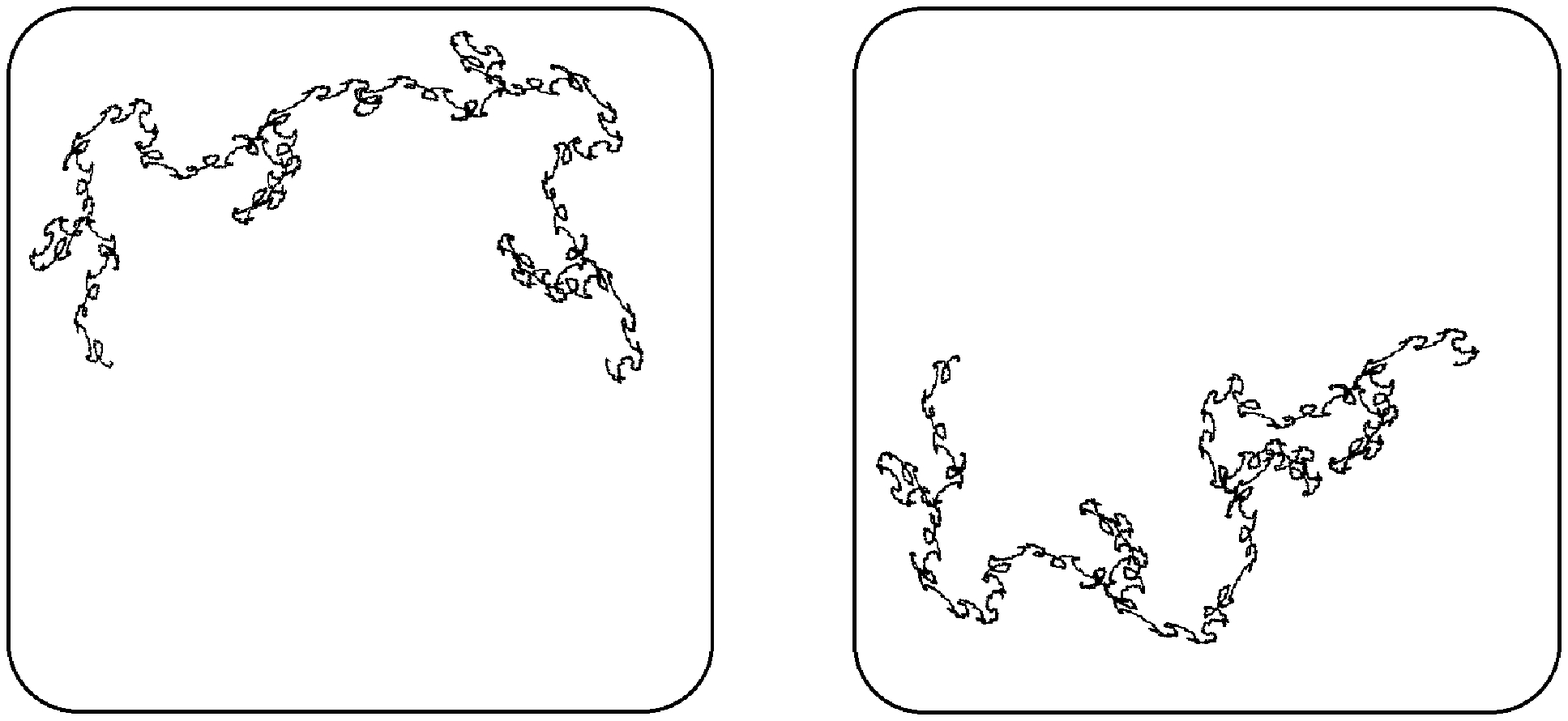}{\raisebox{-1.9752in}{\includegraphics[height=1.9752in]{final3.ps}}}

Notice how the two images in Figure \ref{final2} consist of the union of
shrunken copies of the images in Figure \ref{final1}, while the curves in
Figure \ref{final3} are made from two shrunken copies of the curves in
Figure \ref{final2}.

This example illustrates some typical features of the theory in this paper.
(i) New images are generated,\textit{\ one per iteration per screen}. (ii)
After sufficient iterations for the system to have settled into its
\textquotedblleft stationary state\textquotedblright , each image looks like
a finite resolution rendering of a fractal set that typically changes from
one iteration to the next; each fractal belongs to the same family, in the
present case a family of continuous curves. (iii) In fact, it follows from
the theory that the pictures in this example correspond to curves with this
property: for any $\epsilon >0$ the curve is the union of \textquotedblleft
little\textquotedblright\ curves, ones such that the distance apart of any
two points is no more than $\epsilon $, each of which is an affine
transformation of one of at most two continuous closed paths in $\mathbb{R}%
^{2}$. (iv) We will show that the successive images, or rather the abstract
objects they represent, eventually all lie arbitrarily close to an object
called a \textit{superfractal.} The superfractal is the attractor of a
superIFS which induces a natural invariant probability measure on the
superfractal.\textit{\ The images produced by the algorithm are distributed
according to this measure}. (v) The images produced in the \textquotedblleft
stationary state\textquotedblright\ are independent of the starting images.
For example, if the initial images in the example had been of a dot or a
line instead of a fish, and the same sequence of random choices had been
made, then the images produced in Figures \ref{final1}, \ref{final2}, and %
\ref{final3} would have been the same at the printed resolution.

One similarly obtains V-variable fractals and their properties using $V$,
rather than two, screens and otherwise proceeding similarly. In (iii) each
of the sets of diameter at most $\epsilon $ is an affine transformation of
at most $V$ sets in $\mathbb{R}^{2}$, where these sets again depend upon $%
\epsilon $ and the particular image.

This example and the features just mentioned suggest that superfractals are
of interest because they provide a natural mathematical bridge between
deterministic and random fractals, and because they may lead to practical
applications in digital imaging, special effects, computer graphics, as well
as in the many other areas where fractal geometric modelling is applied.

\subsection{The structure of this paper.}

The main contents of this paper, while conceptually not very difficult,
involves potentially elaborate notation because we deal with iterated
function systems (IFSs) made of IFSs, and probability measures on spaces of
probability measures. So a material part of our effort has been towards a
simplified notation. Thus, below, we set out some notation and conventions
that we use throughout.

The core machinery that we use is basic IFS theory, as described in \cite%
{Hu81} and \cite{BaDe}. So in Section \ref{secIFS} we review relevant parts
of this theory, using notation and organization that extends to and
simplifies later material. To keep the structural ideas clear, we restrict
attention to IFSs with strictly contractive transformations and constant
probabilities. Of particular relevance to this paper, we explain what is
meant by \textit{the random iteration algorithm. }We illustrate the theorems
with simple applications to two-dimensional computer graphics, both to help
with understanding and to draw attention to some issues related to
discretization that apply \textit{a fortiori }in computations\ of V-variable
fractals.

We begin Section \ref{secthree} with the definition of a \textit{superIFS},
namely an IFS made of IFSs. We then introduce associated trees, in
particular labelled trees, the space of code trees $\Omega $, and
construction trees; then we review standard random fractals using the
terminology of trees and superIFSs.

In Section \ref{four} we study a special class of code trees, called $V$%
-variable trees, where $V$ is an integer. What are these trees like? At each
level they have at most $V$ distinct subtrees! In fact these trees are
described with the aid of an IFS $\{\Omega ^{V};\eta ^{a},\mathcal{P}%
^{a},a\in \mathcal{A}\}\mathcal{\ }$where $\mathcal{A}$ is a finite index
set, $\mathcal{P}^{a}$s are probabilities, and each $\eta ^{a}$ is a
contraction mapping from $\Omega ^{V}$ to itself. The IFS enables one to put
a measure attractor on the set of $V$-variable trees, such that they can be
sampled by means of the random iteration algorithm. We describe the mappings 
$\eta ^{a}$ and compositions of them using certain finite doubly labelled
trees. This, in turn, enables us to establish the convergence, as $%
V\rightarrow \infty ,$ of the probability measure on the set of $V$-variable
trees, associated with the IFS and the random iteration algorithm, to a
corresponding natural probability distribution on the space $\Omega $.

In Section \ref{five} the discussion of trees in Section \ref{four} is
recapitulated twice over: the same basic IFS theory is applied in two
successively more elaborate settings, yielding the formal concepts of
V-variable fractals and superfractals. More specifically, in Section \ref%
{five}.1, the superIFS is used to define an IFS\ of functions that map $V$%
-tuples of compact sets into $V$-tuples of compact sets; the attractor of
this IFS is a set of $V$-tuples of compact sets; these compact sets are
named $V$-variable fractals and the set of these $V$-variable fractals is
named a superfractal. We show that these $V$-variable fractals can be
sampled by means of the random iteration algorithm, adapted to the present
setting; that they are distributed according to a certain stationary measure
on the superfractal; and that this measure converges to a corresponding
measure on the set of \textquotedblleft fully\textquotedblright\ random
fractals as $V\rightarrow \infty $, in an appropriate metric. We also
provide a continuous mapping from the set of $V$-variable trees to the set
of $V$-variable fractals, and characterize the $V$-variable fractals in
terms of a property that we name \textquotedblleft $V$-variability%
\textquotedblright . Section \ref{five}.2 follows the same lines as in
Section \ref{five}.1, except that here the superIFS is used to define an
IFS\ that maps $V$-tuples of measures to $V$-tuples of measures; this leads
to the definition and properties of $V$-variable fractal measures. In
Section \ref{five}.3 we describe how to compute the fractal dimensions of
V-variable fractals in certain cases and compare them, in a case involving
Sierpinski triangles, with the fractal dimensions of deterministic fractals,
\textquotedblleft fully\textquotedblright\ random fractals, and
\textquotedblleft homogeneous\textquotedblright\ random fractals that
correspond to $V=1$ and are a special case of a type of random fractal
investigated by Hambly and others \cite{Ham}, \cite{BaHa}, \cite{kif}, \cite%
{stenflo2}.

In Section \ref{six} we describe some potential applications of the theory
including new types of space-filling curves for digital imaging, geometric\
modelling and texture rendering in digital content creation, and random
fractal interpolation for computer aided design systems. In Section \ref%
{seven} we discuss generalizations and extensions of the theory, areas of
ongoing research, and connections to the work of others.

\subsection{Some Notation.}

We use notation and terminology consistent with \cite{BaDe}.

Throughout we reserve the symbols $M$, $N$, and $V$ for positive integers.
We will use the variables $m\in \{1,2,...,M\}$, $n\in \{1,2,...,N\}$, and $%
v\in \{1,2,...,V\}$.

Throughout we use an underlying metric space $(\mathbb{X},d_{\mathbb{X}})$
which is assumed to be compact unless otherwise stated. We write $\mathbb{X}%
^{V}$ to denote the compact metric space%
\begin{equation*}
\underset{V\text{ TIMES }}{\underbrace{\mathbb{X}\times \mathbb{X}\times
...\times \mathbb{X}}}\text{.}
\end{equation*}%
with metric%
\begin{equation*}
d(x,y)=d_{\mathbb{X}^{V}}(x,y)=\max \text{ }\{d_{\mathbb{X}}(x_{v},y_{v})%
\text{ }|\text{ }v=1,2,...,V\}\text{, }\forall x,y\in \mathbb{X}^{V}\text{,}
\end{equation*}%
where $x=(x_{1},x_{2},...,x_{V})$ and $y=(y_{1},y_{2},...,y_{V})$.

In some applications, to computer graphics for example, $(\mathbb{X},d_{%
\mathbb{X}})$ is a bounded region in $\mathbb{R}^{2}$ with the Euclidean
metric, in which case we will usually be concerned with affine or projective
maps.

Let $\mathbb{S=S(X)}$ denote the set of all subsets of $\mathbb{X}$, and let 
$\ C\in \mathbb{S}$. We extend the definition of a function $f:\mathbb{%
X\rightarrow X}$ to $f:\mathbb{S\rightarrow S}$ by 
\begin{equation*}
f(C)=\{f(x)\text{ }|\text{ }x\in C\}
\end{equation*}

Let $\mathbb{H}$=$\mathbb{H}$($\mathbb{X}$) denote the set of non-empty
compact subsets of $\mathbb{X}$. Then if $f:\mathbb{X}\rightarrow \mathbb{X}$
we have $f$ $:\mathbb{H\rightarrow H}$. We use $d_{\mathbb{H}}$ to denote
the Hausdorff metric on $\mathbb{H}$ implied by the metric $d_{\mathbb{X}}$
on $\mathbb{X}$. This is defined as follows. Let $A$\ and $B$ be two sets in 
$\mathbb{H}$, define the distance \textit{from} $A$ \textit{to} $B$ to be 
\begin{equation}
\mathcal{D}(A,B)=\max \{\min \{d_{\mathbb{X}}(x,y)\text{ }|\text{ }y\in B\}%
\text{ }|\text{ }x\in A\},  \label{set_distance}
\end{equation}%
and define the \textit{Hausdorff metric} by%
\begin{equation*}
d_{\mathbb{H}}(A,B)=\max \{\mathcal{D}(A,B),\mathcal{D}(B,A)\}.
\end{equation*}%
Then $\mathbb{(H}$,$d_{\mathbb{H}})$ is a compact metric space. We will
write $(\mathbb{H}^{V},d_{\mathbb{H}^{V}})$ to denote the $V$-dimensional
product space constructed from $\mathbb{(H}$,$d_{\mathbb{H}})$ just as $(%
\mathbb{X}^{V},d_{\mathbb{X}^{V}})$ is constructed from $(\mathbb{X},d_{%
\mathbb{X}})$. When we refer to continuous, Lipschitz, or strictly
contractive functions acting on $\mathbb{H}^{V}$ we assume that the
underlying metric is $d_{\mathbb{H}^{V}}$.

We will in a number of places start from a function acting on a space, and
extend its definition to make it act on other spaces, while leaving the
symbol unchanged as above.

Let $\mathbb{B=B(X)}$ denote the set of Borel subsets of $\mathbb{X}$. Let $%
\mathbb{P=P(X)}$. In some applications to computer imaging one sets $\mathbb{%
X}=[0,1]\times \lbrack 0,1]\subset \mathbb{R}^{2}$ and identifies a black
and white image with a member of $\mathbb{H(X)}$. Greyscale images are
identified with members of $\mathbb{P(X)}$. Probability measures on images
are identified with $\mathbb{P(H(X))}$ or $\mathbb{P(P(X))}$.

Let $d_{\mathbb{P(X)}}$ denote the \textit{Monge Kantorovitch metric} on $%
\mathbb{P(X)}$. This is defined as follows. Let $\mu $ and $\nu $ be any
pair of measures in $\mathbb{P}$. Then 
\begin{equation*}
d_{\mathbb{P}}(\mu ,\nu )=\sup \left\{ \int\limits_{\mathbb{X}}\text{ }fd\mu
-\int\limits_{\mathbb{X}}\text{ }fd\nu \Bigm|f:\mathbb{X\rightarrow R}\text{%
, }|f(x)-f(y)|\leq d_{\mathbb{X}}(x,y)\text{ }\forall \text{ }x,y\in \mathbb{%
X}\right\} .
\end{equation*}%
Then ($\mathbb{P},d_{\mathbb{P}})$ is a compact metric space. The distance
function $d_{\mathbb{P}}$ metrizes the topology of weak convergence of
probability measures on $\mathbb{X}$, \cite{dud}. We define the push-forward
map $f:\mathbb{P(X)\rightarrow P(X)}$ by 
\begin{equation*}
f(\mu )=\mu \circ f^{-1}\text{ \ }\forall \text{ }\mu \in \mathbb{P(X)}\text{%
.}
\end{equation*}%
Again here we have extended the domain of action of the function $f:\mathbb{%
X\rightarrow X}.$

We will use such spaces as $\mathbb{P(H}^{V}\mathbb{)}$ and $\mathbb{P((P(X))%
}^{V})$ (or $\mathbb{H(H}^{V})$ and $\mathbb{H(P}^{V})$ depending on the
context). These spaces may at first seem somewhat Baroque, but as we shall
see, they are very natural. In each case we assume that the metric of a
space is deduced from the space from which it is built, as above, down to
the metric on the lowest space $\mathbb{X}$, and often we drop the subscript
on the metric without ambiguity. So for example, we will write%
\begin{equation*}
d(A,B)=d_{\mathbb{H((P(X))}^{V})}(A,B)\text{ \ }\forall \text{ }A,B\in 
\mathbb{H((P(X))}^{V})\text{.}
\end{equation*}

We also use the following common symbols:

$\mathbb{N}=\{1,2,3,...\}$, $\mathbb{Z=\{}...-2,-1,0,1,2,...\},$ and $%
\mathbb{Z}^{+}\mathbb{=\{}0,1,2,...\}$.

When $S$ is a set, $|S|$ denotes the number of elements of $S.$

\section{\label{secIFS}Iterated Function Systems}

\subsection{Definitions and Basic Results.}

In this section we review relevant information about IFSs. To clarify the
essential ideas we consider the case where all mappings are contractive, but
indicate in Section \ref{five} how these ideas can be generalized. The
machinery and ideas introduced here are applied repeatedly later on in more
elaborate settings.

Let 
\begin{equation}
F=\{\mathbb{X};f_{1},f_{2},...,f_{M};p_{1},p_{2},...,p_{M}\}  \label{IFS}
\end{equation}%
denote an IFS with probabilities. The functions $f_{m}:\mathbb{X\rightarrow X%
}$ are contraction mappings with fixed Lipschitz constant $0\leq l<1$; that
is 
\begin{equation*}
d(f_{m}(x),f_{m}(y))\leq l\cdot d(x,y)\text{ }\forall x,y\in \mathbb{X}%
,\forall m\in \{1,2,...,M\}\text{.}
\end{equation*}%
The $p_{m}$'s are probabilities, with 
\begin{equation*}
\sum\limits_{m=1}^{M}p_{m}=1\text{, }p_{m}\geq 0\text{ }\forall m\text{.}
\end{equation*}

We define mappings $F:\mathbb{H(X)}\rightarrow \mathbb{H(X)}$ and $F:\mathbb{%
P(X)}\rightarrow \mathbb{P(X)}$ by 
\begin{equation*}
F(K)=\bigcup\limits_{m=1}^{M}f_{m}(K)\text{ \ }\forall K\in \mathbb{H}\text{,%
}
\end{equation*}%
and%
\begin{equation*}
F(\mu )=\sum\limits_{m=1}^{M}p_{m}f_{m}(\mu )\text{ \ }\forall \mu \in 
\mathbb{P}\text{.}
\end{equation*}%
In the latter case note that the weighted sum of probability measures is
again a probability measure.

\begin{theorem}
\label{Hutchinson1}\cite{Hu81}The mappings $F:\mathbb{H(X)}\rightarrow 
\mathbb{H(X)}$ and $F:\mathbb{P(X)}\rightarrow \mathbb{P(X)}$ are both
contractions with factor $0\leq l<1$. That is,%
\begin{equation*}
d(F(K),F(L))\leq l\cdot d(K,L)\text{ }\forall \text{ }K,L\in \mathbb{H(X)}%
\text{,}
\end{equation*}%
and%
\begin{equation*}
d(F(\mu ),F(\nu ))\leq l\cdot d(\mu ,\nu )\text{ }\forall \text{ }\mu ,\nu
\in \mathbb{P(X)}\text{.}
\end{equation*}%
As a consequence, there exists a unique nonempty compact set $A\in \mathbb{%
H(X)}$ such that%
\begin{equation*}
F(A)=A\text{,}
\end{equation*}%
and a unique measure $\mu \in \mathbb{P(X)}$ such that 
\begin{equation*}
F(\mu )=\mu \text{.}
\end{equation*}%
The support of $\mu $ is contained in, or equal to $A$, with equality when
all of the probabilities $p_{m}$ are strictly positive.
\end{theorem}

\begin{definition}
The set $A$ in Theorem \ref{Hutchinson1} is called the set attractor of the
IFS $F$, and the measure $\mu $ is called the measure attractor of $F$.
\end{definition}

We will use the term \textit{attractor of an IFS} to mean either the set
attractor or the measure attractor. We will also refer informally to the set
attractor of an IFS as a \textit{fractal} and to its measure attractor as a 
\textit{fractal measure, }and to either as a \textit{fractal}. Furthermore,
we say that the set attractor of an IFS is a \textit{deterministic} fractal.
This is in distinction to \textit{random} fractals, and in particular to $V$%
\textit{-variable} random fractals which are the main goal of this paper.

There are two main types of algorithms for the practical computation of
attractors of IFS that we term \textit{deterministic algorithms }and \textit{%
random iteration algorithms}, also known as backward and forward algorithms,
c.f.\ \cite{DiaFree}. These terms should not be confused with the type of
fractal that is computed by means of the algorithm. Both deterministic and
random iteration algorithms may be used to compute deterministic fractals,
and as we discuss later, a similar remark applies to our $V$-variable
fractals.

Deterministic algorithms are based on the following:

\begin{corollary}
\label{detalgcor}Let $A_{0}\in \mathbb{H(X)}$, or $\mu _{0}\in \mathbb{P(X)}$%
, and define recursively 
\begin{equation*}
A_{k}=F(A_{k-1})\text{, or }\mu _{k}=F(\mu _{k-1})\text{, }\forall k\in 
\mathbb{N}\text{, }
\end{equation*}%
respectively; then 
\begin{equation}
\underset{k\rightarrow \infty }{\lim }A_{k}=A\text{, or }\underset{%
k\rightarrow \infty }{\lim }\mu _{k}=\mu \text{,}  \label{alg1}
\end{equation}%
respectively. The rate of convergence is geometrical; for example, 
\begin{equation*}
d(A_{k},A)\leq l^{k}\cdot d(A_{0},A)\text{ }\forall k\in \mathbb{N}\text{.}
\end{equation*}
\end{corollary}

In practical applications to two-dimensional computer graphics, the
transformations and the spaces upon which they act must be discretized. The
precise behaviour of computed sequences of approximations to an attractor of
an IFS depends on the details of the implementation and is generally quite
complicated; for example, the discrete IFS may have multiple attractors, see 
\cite{Peruggia}, Chapter 4. The following example gives the flavour of such
applications.

\begin{example}
\label{determ_example}In Figure \ref{texfishx} we illustrate a practical
deterministic algorithm, based on the first formula in Equation (\ref{alg1})
starting from a simple IFS on the unit square $\square \subset $ $\mathbb{R}%
^{2}$. The IFS is $F=\{\square
;f_{1}^{1},f_{2}^{1},f_{2}^{2};0.36,0.28,0.36\}$ where the transformations
are defined in Equations (\ref{trans_one}), (\ref{trans_two}), and (\ref%
{trans_four}). The successive images, from left to right, from top to
bottom, represent $A_{0}$, $A_{2}$, $A_{5}$, $A_{7}$, $A_{20}$, and $A_{21}$%
. In the last two images, the sequence appears to have converged at the
printed resolution to representations of the set attractor. Note however
that the initial image is partitioned into two subsets corresponding to the
colours red and green. Each successive computed image is made of pixels
belonging to a discrete model for $\square $ and consists of red pixels and
green pixels. Each pixel corresponds to a set of points in $\mathbb{R}^{2}$.
But for the purposes of computation only one point corresponding to each
pixel is used. When both a point in a red pixel and point in a green pixel
belonging to say $A_{n}$ are mapped under $F$ to points in the same pixel in 
$A_{n+1}$ a choice has to be made about which colour, red or green, to make
the new pixel of $A_{n+1}$. Here we have chosen to make the new pixel of $%
A_{n+1}$ the same colour as that of the pixel containing the last point in $%
A_{n}$, encountered in the course of running the computer program, to be
mapped to the new pixel. The result is that, although the sequence of
pictures converge to the set attractor of the IFS, the colours themselves do
not settle down, as illustrated in Figure \ref{texture_effect}. We call this
\textquotedblleft the texture effect\textquotedblright , and comment on it
in Example \ref{texture_example}. In printed versions of the figures
representing $A_{20}$, and $A_{21}$ the red and green pixels are somewhat
blended.\FRAME{ftbpFU}{4.8058in}{3.141in}{0pt}{\Qcb{An illustration of the
deterministic algorithm.}}{\Qlb{texfishx}}{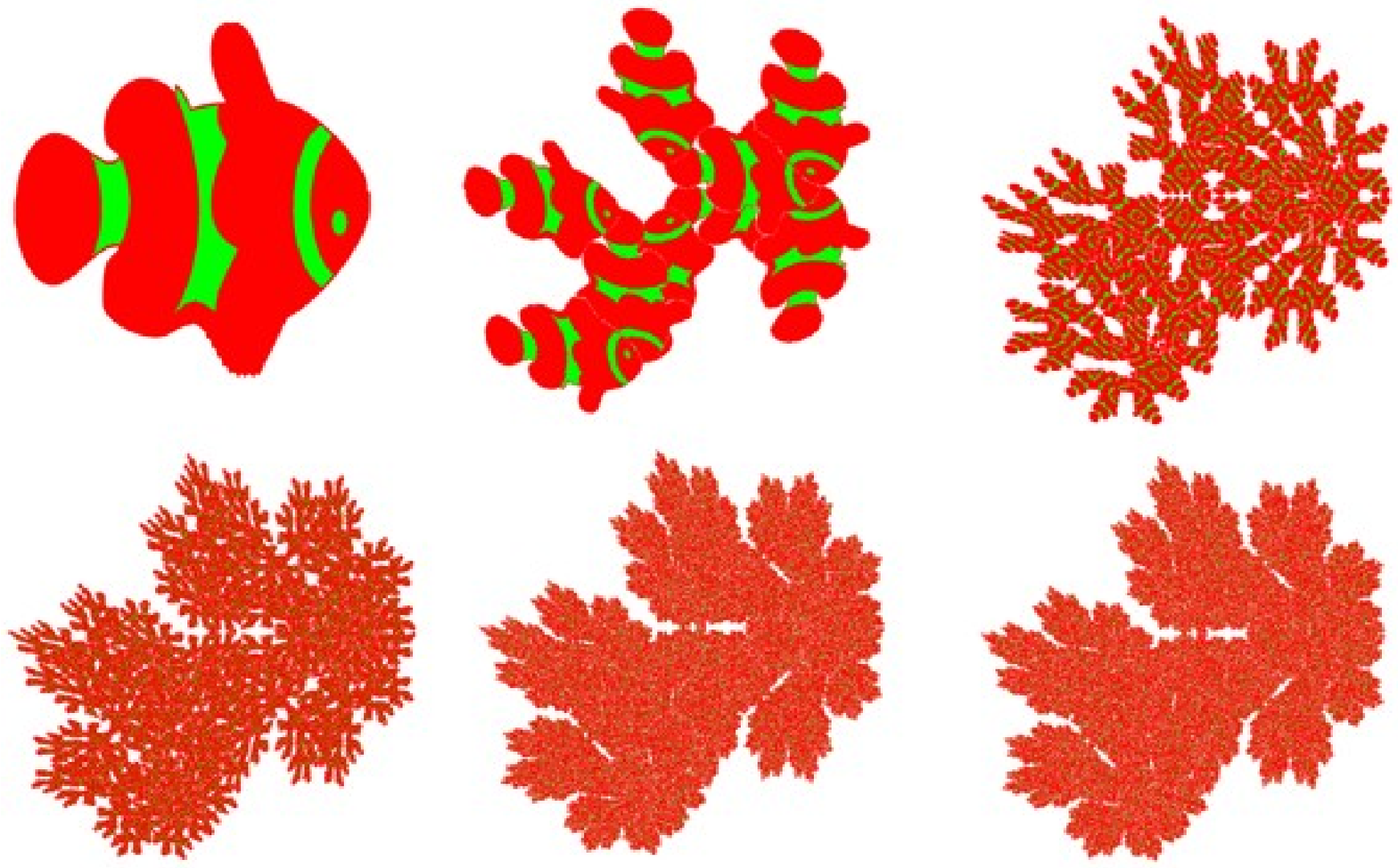}{\raisebox{-3.141in}{\includegraphics[height=3.141in]{texfishx.ps}}}
\end{example}

The following theorem is the mathematical justification and description of 
\textit{the random iteration algorithm}. It follows from Birkhoff's ergodic
theorem and our assumption of contractive maps. A more general version of it
is proved in \cite{JEergodic}.

\begin{theorem}
\label{RIA}Specify a starting point $x_{1}\in \mathbb{X}$. Define a random
orbit of the IFS to be $\{x_{l}\}_{l=1}^{\infty }$ where $%
x_{l+1}=f_{m}(x_{l})$ with probability $p_{m\text{ }}$. Then for almost all
random orbits $\{x_{l}\}_{l=1}^{\infty }$ we have:%
\begin{equation}
\mu (B)=\underset{l\rightarrow \infty }{\lim }\frac{|\{B\cap
\{x_{1},x_{2},...,x_{l}\}|}{l}\text{.}  \label{weak}
\end{equation}%
for all $B\in \mathbb{B(X)}$ such that $\mu (\partial B)=0,$ where $\partial
B$ denotes the boundary of $B$.
\end{theorem}

\begin{remark}
This is equivalent by standard arguments to the following: for any $x_{1}\in 
\mathbb{X}$ and almost all random orbits the sequence of point measures $%
\frac{1}{l}(\delta _{x_{1}}+\delta _{x_{2}}+...+\delta _{x_{l}})$ converges
in the weak sense to $\mu $, see for example \cite{Bi68}, pages 11 and 12.
(Weak convergence of probability measures is the same as convergence in the
Monge Kantorovitch metric, see \cite{dud}, pages 310 and 311.)
\end{remark}

The random iteration algorithm can be applied to the computation of
two-dimensional computer graphics. It has benefits compared to deterministic
iteration of low memory requirements, high accuracy --- as the iterated
point can be kept at much higher precision than the resolution of the
computed image --- and it allows the efficient computation of zooms into
small parts of an image. However, as in the case of deterministic
algorithms, the images produced depend on the computational details of image
resolution, the precision to which the points $\{x_{1},x_{2},...,x_{l}\}$
are computed, the contractivity of the transformations, the way in which
Equation (\ref{weak}) is implemented, choices of colours, etc. Different
implementations can produce different results.

\begin{example}
\label{greeny}Figure \ref{InvMeas3} shows a \textquotedblleft
picture\textquotedblright\ of the invariant measure of the IFS in Example %
\ref{determ_example}, computed using a discrete implementation of the random
iteration algorithm, as follows.\ Pixels corresponding to a discrete model
for $\square \subset $ $\mathbb{R}^{2}$ are assigned the colour white.
Successive floating point coordinates of points in $\square $ are computed
by random iteration and the first (say) one hundred points are discarded.
Thereafter, as each new point is calculated, the pixel to which it belongs
is set to black. This phase of the computation continues until the pixels
cease to change, and produces a black image of the support of the measure,
the set attractor of the IFS, against a white background. Then the random
iteration process is continued, and as each new point is computed the green
component of the pixel to which the latest point belongs is brightened by a
fixed amount. Once a pixel is at brightest green, its value is not changed
when later points belong to it. The computation is continued until a balance
is obtained between that part of the image which is brightest green and that
which is lightest green, and then stopped. \FRAME{ftbpFU}{2.5469in}{2.3921in%
}{0pt}{\Qcb{\textquotedblleft Picture\textquotedblright\ of the measure
attractor of an IFS with probabilities produced by the random iteration
algorithm. The measure is depicted in shades of green, from 0 (black) to 255
(bright green).}}{\Qlb{InvMeas3}}{invmeas4.ps}{\raisebox{-2.3921in}{\includegraphics[height=2.3921in]{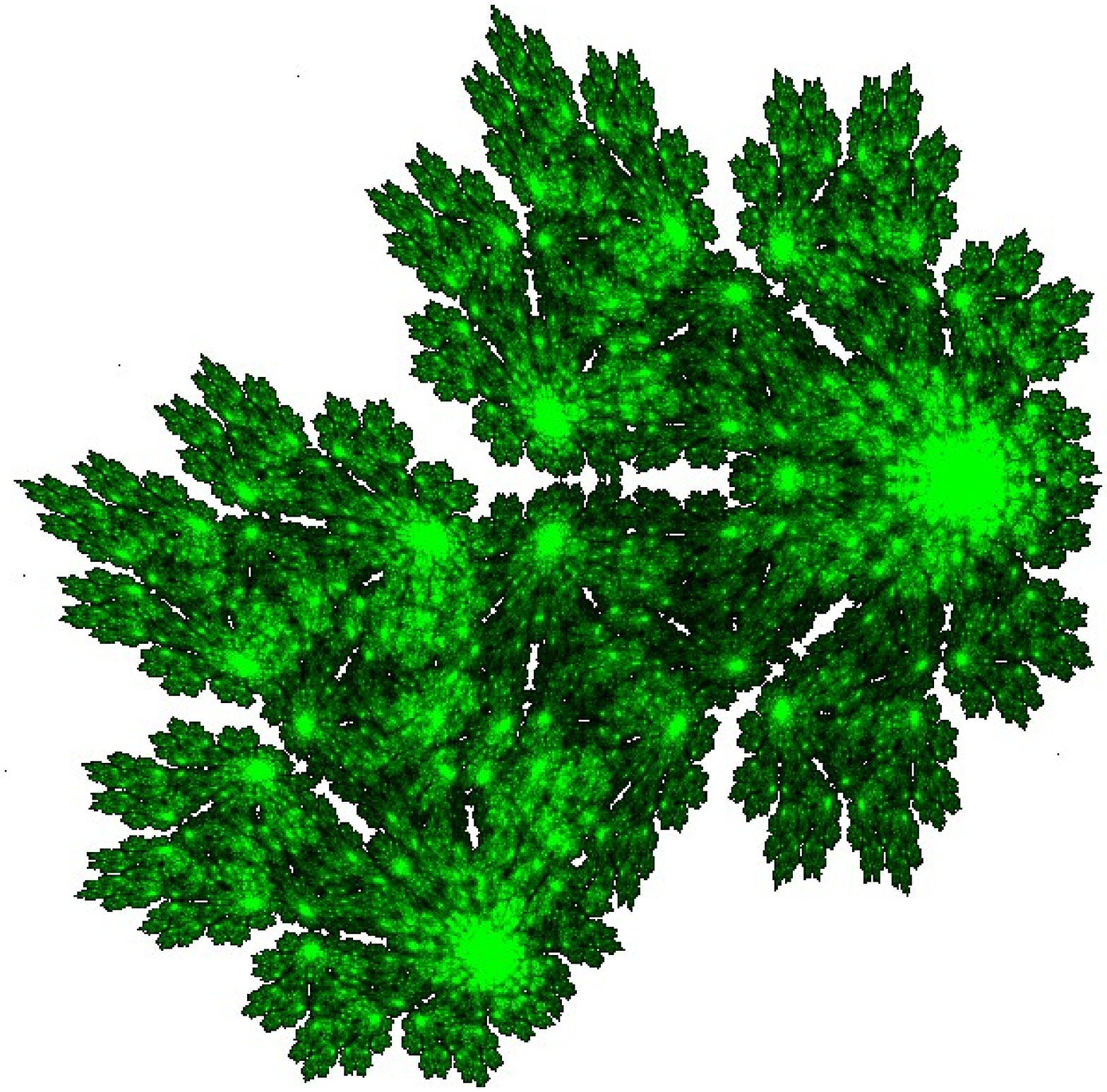}}}
\end{example}

The following theorem expresses the ergodicity of the IFS $F$. The proof
depends centrally on the uniqueness of the measure attractor. A variant of
this theorem, weaker in the constraints on the IFS but stronger in the
conditions on the set $B$, and stated in the language of stochastic
processes, is given by \cite{JEergodic}. We prefer the present version for
its simple statement and direct measure theoretic proof.

\begin{theorem}
\label{ergodic}Suppose that $\mu $ is the unique measure attractor for the
IFS $F$. Suppose $B\in \mathbb{B(X)}$ is such that $f_{m}(B)\subset B$ \ $%
\forall m\in \{1,2,...,M\}$. Then $\mu (B)=0$ or $1$.
\end{theorem}

\begin{proof}
Let us define the measure $\mu \lfloor B$ ($\mu $ restricted by $B$) by ($%
\mu \lfloor B)(E)=\mu (B\cap E)$. The main point of the proof is to show
that $\mu \lfloor B$ is invariant under the IFS $F$. (\ A similar result
applies to $\mu \lfloor B^{C}$ where $B^{C}$ denotes the complement of $B$.)

If $E\subset B^{C}$, for any $m$, since $f_{m}(B)\subset B$, 
\begin{equation*}
f_{m}(\mu \lfloor B)(E)=\mu (B\cap f_{m}^{-1}(E))=\mu (\emptyset )=0\text{.}
\end{equation*}

Moreover, 
\begin{equation}
\mu (B)=f_{m}(\mu \lfloor B)(\mathbb{X})=f_{m}(\mu \lfloor B)(B)\text{.}
\label{two}
\end{equation}

It follows that 
\begin{eqnarray*}
\mu (B) &=&\sum\limits_{m=1}^{M}p_{m}f_{m}\mu
(B)=\sum\limits_{m=1}^{M}p_{m}f_{m}(\mu \lfloor
B)(B)+\sum\limits_{m=1}^{M}p_{m}f_{m}(\mu \lfloor B^{C})(B) \\
&=&\mu (B)+\sum\limits_{m=1}^{M}p_{m}f_{m}(\mu \lfloor B^{C})(B)\text{ (from
(\ref{two})).}
\end{eqnarray*}

Hence%
\begin{equation}
\sum\limits_{m=1}^{M}p_{m}f_{m}(\mu \lfloor B^{C})(B)=0\text{.}
\label{three}
\end{equation}%
Hence for any measurable set $E\subset \mathbb{X}$ 
\begin{eqnarray*}
(\mu \lfloor B)(E) &=&\mu (B\cap E)=\sum\limits_{m=1}^{M}p_{m}f_{m}\mu
(B\cap E) \\
&=&\sum\limits_{m=1}^{M}p_{m}f_{m}(\mu \lfloor B)(B\cap
E)+\sum\limits_{m=1}^{M}p_{m}f_{m}(\mu \lfloor B^{C})(B\cap E) \\
&=&\sum\limits_{m=1}^{M}p_{m}f_{m}(\mu \lfloor B)(E)+0\text{ \ (using (\ref%
{three})).}
\end{eqnarray*}%
Thus $\mu \lfloor B$ is invariant and so is either the zero measure or for
some constant $c\geq 1$ we have $c\mu \lfloor B=\mu $ (by uniqueness)$=\mu
\lfloor B+\mu \lfloor B^{C}$. This implies $\mu \lfloor B^{C}=0$ and in
particular $\mu (B^{C})=0$ and $\mu (B)=1$.
\end{proof}

\begin{example}
\label{texture_example} Figure \ref{texture_effect} shows close-ups on the
two images at the bottom left in Figure \ref{texfishx}, see Example \ref%
{determ_example}. At each iteration it is observed that the pattern of red
and green pixels changes in a seemingly random manner. A similar\ texture
effect is often observed in other implementations and in applications of $V$%
-variable fractals to computer graphics. Theorem \ref{ergodic} provides a
simple model explanation for this effect as follows. Assume that the red
pixels and the green pixels both correspond to sets of points of positive
measure, both invariant under F. Then we have a contradiction to the
Corollary above. So neither the red nor the green set can be invariant under
F. Hence, either one of the sets disappears --- which occurs in some other
examples --- or the pixels must jump around. A similar argument applied to
powers of $F$ shows that the way the pixels jump around cannot be periodic,
and hence must be $``$random\textquotedblright . (A more careful explanation
involves numerical and statistical analysis of the specific computation.)%
\FRAME{ftbpFU}{4.5455in}{2.4915in}{0pt}{\Qcb{Close-up on the same small
region in each of the bottom left two images in Figure \protect\ref{texfishx}%
, showing the texture effect; the distribution of red and green pixels
changes with each iteration.}}{\Qlb{texture_effect}}{sample1.ps}{\raisebox{-2.4915in}{\includegraphics[height=2.4915in]{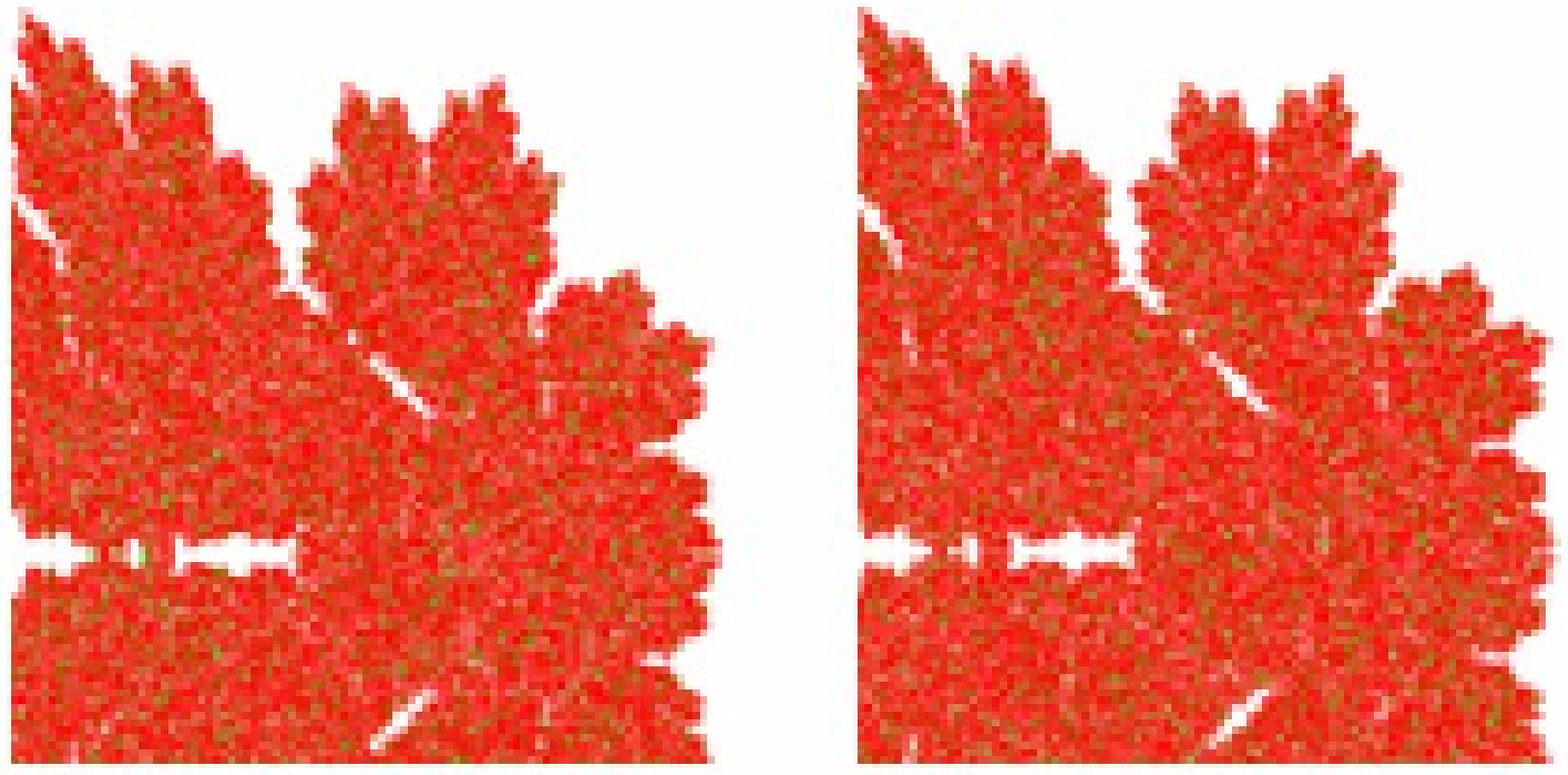}}}
\end{example}

\subsection{Fractal Dimensions.}

In the literature there are many different definitions of a theoretical
quantity called the \textquotedblleft fractal dimension\textquotedblright\
of a subset of $\mathbb{X}$. A mathematically convenient definition of the
fractal dimension of a set $S\subset \mathbb{X}$ is the Hausdorff dimension.
This is always well-defined. Its numerical value often but not always
coincides with the values provided by other definitions, when they apply.

Fractal dimensions are useful for a number of reasons. They can be used to
compare and classify sets and measures and they have some natural invariance
properties. For example the Hausdorff dimension of a set $S$ is invariant
under any bi-Lipshitz transformation; that is, if $f:\mathbb{X\rightarrow X}$
is such that there are constants $c_{1}$and $c_{2}$ in $(0,\infty )$ with $%
c_{1}\cdot d(x,y)\leq d(f(x),f(y))\leq c_{1}\cdot d(x,y)$ $\forall $ $x,y\in 
\mathbb{X}$ then the Hausdorff dimension of $S$ is the same as that of $f(S)$%
. Fractal dimensions are useful in fractal image modelling: for example,
empirical fractal dimensions of the boundaries of clouds can be used as
constraints in computer graphics programs for simulating clouds. Also, as we
will see below, the specific value of the Hausdorff dimension of the set
attractor $A$ of an IFS can yield the probabilities for most efficient
computation of $A$ using the random iteration algorithm. For these same
reasons fractal dimensions are an important concept for V-variable fractals
and superfractals.

The following two definitions are discussed in \cite{Falconer} pp. 25 et seq.

\begin{definition}
Let $S\subset \mathbb{X}$, $\delta >0$, and $0\leq s<\infty $. Let 
\begin{equation*}
H_{\delta }^{s}(S)=\inf \left\{ \sum\limits_{i=1}^{\infty }|U_{i}|^{s}\Bigm|%
\{U_{i}\}\text{ is a }\delta -\text{cover of }S\right\} ,
\end{equation*}%
where $|U_{i}|^{s}$ denotes the $s^{th}$ power of the diameter of the set $%
U_{i}$, and where a $\delta -$cover of $S$ is a covering of $S$ by subsets
of $\mathbb{X}$ of diameter less than $\delta $. Then the $s$-dimensional
Hausdorff measure of the set $S$ is defined to be 
\begin{equation*}
H^{s}(S)=\underset{\delta \rightarrow 0}{\lim }H_{\delta }^{s}(S).
\end{equation*}
\end{definition}

The s-dimensional Hausdorff measure is a Borel measure but is not normally
even $\sigma $-finite.

\begin{definition}
The Hausdorff dimension of the set $S$ $\subset \mathbb{X}$ is defined to be%
\begin{equation*}
\dim _{H}S=\inf \{s\text{ }|\text{ }H^{s}(S)=0\}.
\end{equation*}
\end{definition}

The following quantity is often called the \textit{fractal dimension} of the
set $S$. It can be approximated in practical applications, by estimating the
slope of the graph of the logarithm of the number of \textquotedblleft
boxes\textquotedblright\ of side length $\delta $ that intersect $S,$ versus
the logarithm of $\delta .$

\begin{definition}
The box-counting dimension of the set $S$ $\subset \mathbb{X}$ is defined to
be 
\begin{equation*}
\dim _{B}S=\lim_{\delta \rightarrow 0}\frac{\log N_{\delta }(S)}{\log
(1/\delta )}
\end{equation*}%
if and only if this limit exists, where $N_{\delta }(S)$ is the smallest
number of sets of diameter $\delta >0$ that can cover $S$.
\end{definition}

In order to provide a precise calculation of box-counting and Hausdorff
dimension of the attractor of an IFS we need the following condition.

\begin{definition}
The IFS $F$\ is said to obey the open set condition if there exists a
non-empty open set $O$ such that 
\begin{equation*}
F(O)=\bigcup\limits_{m=1}^{M}f_{m}(O)\subset O,
\end{equation*}%
and%
\begin{equation*}
f_{m}(O)\cap f_{l}(O)=\varnothing \text{ if }m\neq l.
\end{equation*}
\end{definition}

The following theorem provides the Hausdorff dimension of the attractor of
an IFS in some special cases.

\begin{theorem}
\label{dimension_theorem}Let the IFS $F$ consist of similitudes, that is $%
f_{m}(x)=s_{m}O_{m}x+t_{m}$ where $O_{m}$ is an orthonormal transformation
on $\mathbb{R}^{K}$, $s_{m}\in (0,1)$, and $t_{m}\in \mathbb{R}^{K}$. Also
let $F$ obey the open set condition, and let $A$ denote the set attractor of 
$F$. Then 
\begin{equation*}
\dim _{H}A=\dim _{B}A=D
\end{equation*}%
where $D$ is the unique solution of 
\begin{equation}
\sum_{m=1}^{M}s_{m}^{D}=1\text{.}  \label{equation_dimension}
\end{equation}%
Moreover,%
\begin{equation*}
0<H^{D}(A)<\infty \text{.}
\end{equation*}
\end{theorem}

\begin{proof}
This theorem, in essence, was first proved by Moran in 1946, \cite{moran}. A
full proof is given in \cite{Falconer} p.118.
\end{proof}

A\ good choice for the probabilities, which ensures that the points obtained
from the random iteration algorithm are distributed uniformly around the set
attractor in case the open set condition applies, is $p_{m}=s_{m}^{D}$. Note
that the choice of $D$ in Equation \ref{equation_dimension} is the unique
value which makes $(p_{1},p_{2},...,p_{M})$ into a probability vector.

\subsection{Code Space\label{codespace}}

A good way of looking at an IFS $F$ as in (\ref{IFS}) is in terms of the
associated \textit{code space} $\Sigma =\{1,2,...,M\}^{\infty }$. Members of 
$\Sigma $ are infinite sequences from the alphabet $\{1,2,...,M\}$ and
indexed by $\mathbb{N}$. We equip $\Sigma $ with the metric $d_{\Sigma }$
defined for $\omega \neq \varkappa $ by 
\begin{equation*}
d_{\Sigma }(\omega ,\varkappa )=\frac{1}{M^{k}}\text{,}
\end{equation*}%
where $k$ is the index of the first symbol at which $\omega $ and $\varkappa 
$ differ. Then $(\Sigma ,d_{\Sigma })$ is a compact metric space.

\begin{theorem}
\label{Add}Let $A$ denote the set attractor of the IFS $F$. Then there
exists a continuous onto mapping $F:\Sigma \rightarrow A$, defined for all $%
\sigma _{1}\sigma _{2}\sigma _{3}...\in \Sigma $ by 
\begin{equation*}
F(\sigma _{1}\sigma _{2}\sigma _{3}...)=\underset{k\rightarrow \infty }{\lim 
}f_{\sigma _{1}}\circ f_{\sigma _{2}}\circ ...\circ f_{\sigma _{k}}(x)\text{.%
}
\end{equation*}%
The limit is independent of $x\in \mathbb{X}$ and the convergence is uniform
in $x$.
\end{theorem}

\begin{proof}
This result is contained in \cite{Hu81} Theorem 3.1(3).
\end{proof}

\begin{definition}
The point $\sigma _{1}\sigma _{2}\sigma _{3}...\in \Sigma $ is called an
address of the point $F(\sigma _{1}\sigma _{2}\sigma _{3}...)\in A$.
\end{definition}

Note that $F:\Sigma \rightarrow A$ is not in general one-to-one.

The following theorem characterizes the measure attractor of the IFS $F$ as
the push-foward, under $F:\Sigma \rightarrow A$, of an elementary measure $%
\rho \in \mathbb{P}(\Sigma )$, the measure attractor of a fundamental IFS on 
$\Sigma $.

\begin{theorem}
\label{shiftmeas}For each $m\in \{1,2,...,M\}$ define the shift operator $%
s_{m}:\Sigma \rightarrow \Sigma $ by 
\begin{equation*}
s_{m}(\sigma _{1}\sigma _{2}\sigma _{3}...)=m\sigma _{1}\sigma _{2}\sigma
_{3}...
\end{equation*}%
$\forall \sigma _{1}\sigma _{2}\sigma _{3}...\in \Sigma $. Then $s_{m}$ is a
contraction mapping with contractivity factor $\frac{1}{M}$. Consequently%
\begin{equation*}
S:=\{\Sigma ;s_{1,}s_{2},...,s_{M};p_{1},p_{2},...,p_{M}\}
\end{equation*}%
is an IFS. Its set attractor is $\Sigma $. Its measure attractor is the
unique measure $\pi \in \mathbb{P(}\Sigma )$ such that 
\begin{equation*}
\pi \{\omega _{1}\omega _{2}\omega _{3}...\in \Sigma |\omega _{1}=\sigma
_{1},\omega _{2}=\sigma _{2},...,\omega _{k}=\sigma _{k}\}=p_{\sigma
_{1}}\cdot p_{\sigma _{2}}\cdot ...\cdot p_{\sigma _{k}}
\end{equation*}%
$\forall k\in \mathbb{N}$, $\forall \sigma _{1},\sigma _{1},...,\sigma
_{k}\in \{1,2,...,M\}$.

If $\mu $ is the measure attractor of the IFS $F$, with $F:\Sigma
\rightarrow A$ defined as in Theorem \ref{Add}, then 
\begin{equation*}
\mu =F(\pi )\text{.}
\end{equation*}
\end{theorem}

\begin{proof}
This result is \cite{Hu81} Theorem 4.4(3) and (4).
\end{proof}

We call $S$ the \textit{shift} IFS on code space. It has been well studied
in the context of information theory and dynamical systems, see for example 
\cite{Bixx}, and results can often be lifted to the IFS $F$. For example,
when the IFS is non-overlapping, the entropy (see \cite{BaDuXi} for the
definition) of the stationary stochastic process associated with $F$ is the
same as that associated with the corresponding shift IFS, namely: $-\sum
p_{m}\log p_{m}$.

\section{$\label{secthree}$Trees of Iterated Function Systems and Random
Fractals}

\subsection{SuperIFSs}

Let $(X,d_{X})$ be a compact metric space, and let $M$ and $N$ be positive
integers. For $n\in \{1,2,...,N\}$ let $F^{n}$ denote the IFS%
\begin{equation*}
F^{n}=\{\mathbb{X}%
;f_{1}^{n},f_{2}^{n},...f_{M}^{n};p_{1}^{n},p_{2}^{n},...p_{M}^{n}\}
\end{equation*}%
where each $f_{m}^{n}:X\rightarrow X$ is a Lipshitz function with Lipschitz
constant $0\leq l<1$ and the $p_{m}^{n}$ 's are probabilities with 
\begin{equation*}
\sum\limits_{m=1}^{M}p_{m}^{n}=1\text{, }p_{m}^{n}\geq 0\text{ }\forall 
\text{ }m,n\text{.}
\end{equation*}%
Let 
\begin{equation}
\mathcal{F}=\{\mathbb{X};F^{1},F^{2},...,F^{N};P_{1},P_{2},...,P_{N}\}\text{,%
}  \label{superIFS}
\end{equation}%
where the $P_{n}$'s are probabilities with 
\begin{equation*}
\sum\limits_{n=1}^{N}P_{n}=1\text{, }P_{n}\geq 0\text{ \ }\forall n\in
\{1,2,...,N\}\text{.}
\end{equation*}%
$P_{n}>0$ $\forall n\in \{1,2,...,N\}$ and $\sum_{n=1}^{N}P_{n}=1$.

As we will see in later sections, given any positive integer $V$ we can use
the set of IFSs $\mathcal{F}$ to construct a single IFS acting on $\mathbb{%
H(X)}^{V}$. In such cases we call $\mathcal{F}$ a \textit{superIFS}.
Optionally, we will drop the specific reference to the probabilities.

\subsection{\label{treesec}Trees}

We associate various collections of trees with $\mathcal{F}$ and the
parameters $M$ and $N$.

Let $T$ denote the ($M$-fold) set of finite sequences from $\{1,2,...,M\}$,
including the empty sequence $\emptyset $. Then $T$ is called a \textit{tree 
}and the sequences are called the \textit{nodes }of the tree\textit{. F}or $%
i=i_{1}i_{2}...i_{k}\in T$ let $|i|=k$. The number $k$ is called the \textit{%
level} of the \textit{node }$\sigma $. The bottom node $\emptyset $ is at
level zero. If $j=j_{1}j_{2}...j_{l}\in T$ then $ij$ is the concatenated
sequence $i_{1}i_{2}...i_{k}j_{1}j_{2}...j_{l}$.

We define a \textit{\ level}-$k$ ($M$-fold) tree, or a tree of height $k$, $%
T_{k}$ to be the set of nodes of $T$ of level less than or equal to $k$.

A \textit{labelled tree }is a function whose domain is a tree or a level-$k$
tree. A \textit{limb }of a tree is either an ordered pair of nodes of the
form $(i,im)$ where $i\in T$ and $m\in \{1,2,...,M\}$, or the pair of nodes $%
(\emptyset ,\emptyset )$, which is also called the \textit{trunk}. In
representations of labelled trees, as in Figures \ref{Twotree} and \ref%
{col_tree4}, limbs are represented by line segments and we attach the labels
either to the nodes where line segments meet or to the line segments
themselves, or possibly to both the nodes and limbs when a labelled tree is
multivalued. For a two-valued labelled tree $\tau $ we will write 
\begin{equation*}
\tau (i)=(\tau (\text{\textit{node }}i),\tau (\text{\textit{limb }}i))\text{
for }i\in T,
\end{equation*}%
to denote the two components.

A \textit{code} tree is a labelled tree whose range is $\{1,2,...,N\}$. We
write 
\begin{equation*}
\Omega =\{\tau \text{ }|\text{ }\tau :T\rightarrow \{1,2,...,N\}\}
\end{equation*}%
for the set of all infinite code trees.

We define the \textit{subtree} $\widetilde{\tau }:T\rightarrow \{1,2,...,N\}$
of a labelled tree $\tau :T\rightarrow \{1,2,...,N\},$ corresponding to a
node\textit{\ }$i=i_{1}i_{2}...i_{k}\in T,$ by 
\begin{equation*}
\widetilde{\tau }(j)=\tau (ij)\text{ }\forall \text{ }j\in T\text{.}
\end{equation*}%
In this case we say that $\widetilde{\tau }$ is a subtree of $\tau $ \textit{%
at level }$k$. (One can think\ of a subtree as a branch of a tree.)

\FRAME{ftbpFU}{4.7245in}{3.5362in}{0pt}{\Qcb{Pictorial representation of a
level-4 2-fold tree labelled by the sequences corresponding to its nodes.
The labels on the fourth level are shown for every second node. The line
segments between the nodes and the line segment below the bottom node are
referred to as \textit{limbs}. The bottom limb is also called the \textit{%
trunk.}}}{\Qlb{Twotree}}{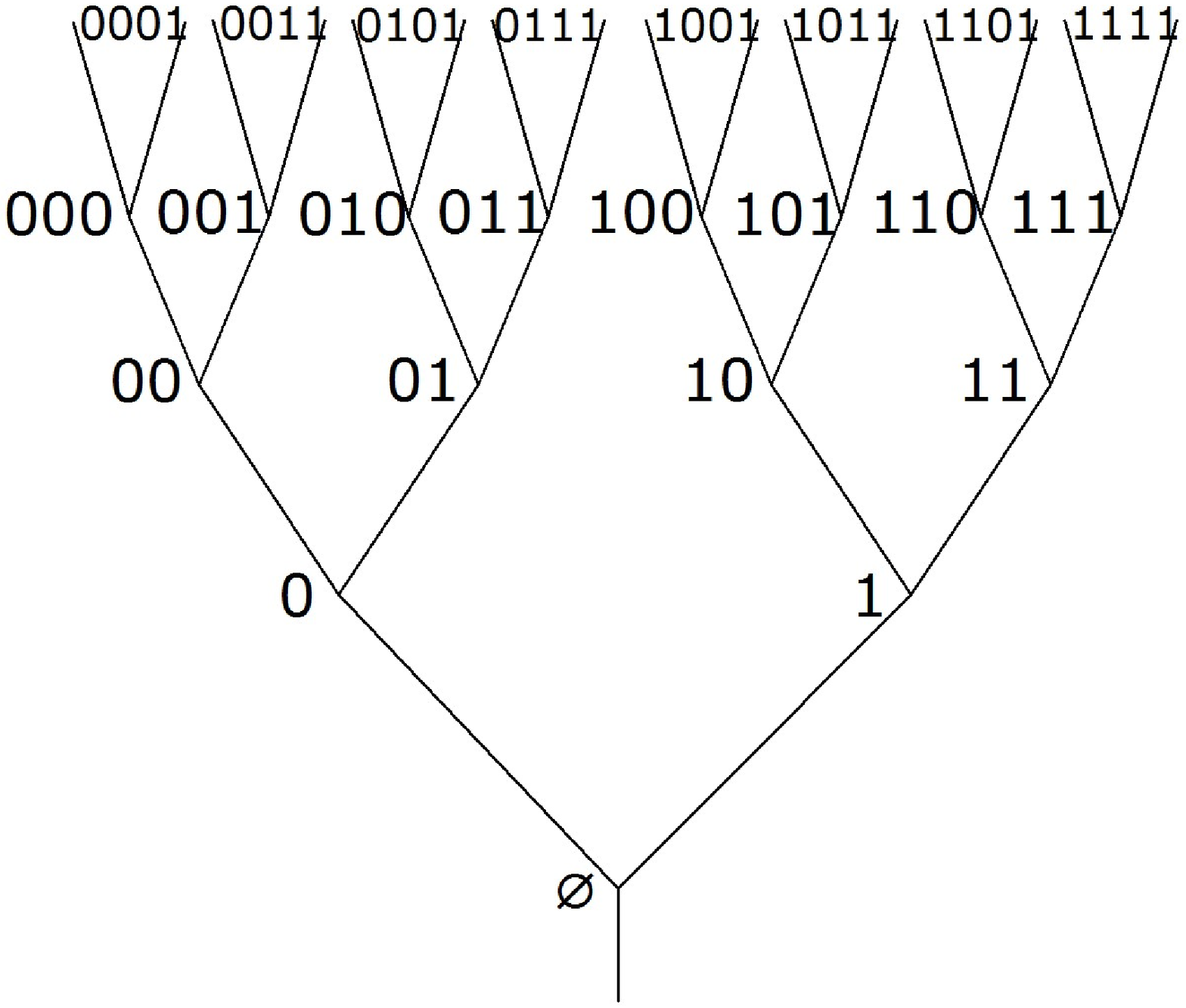}{\raisebox{-3.5362in}{\includegraphics[height=3.5362in]{twotree.drw.ps}}}

Suppose $\tau $ and $\sigma $ are labelled trees with $|\tau |\leq |\sigma |$%
, where we allow $|\sigma |=\infty $. We say that $\sigma $ extends $\tau $,
and $\tau $ is an \textit{initial segment} of $\sigma $, if $\sigma $ and $%
\tau $ agree on their common domain, namely at nodes up to and including
those at level $|\tau |$. We write 
\begin{equation*}
\tau \prec \sigma \text{.}
\end{equation*}%
If $\tau $ is a level-$k$ code tree, the corresponding \textit{cylinder set}
is defined by%
\begin{equation*}
\lbrack \tau ]=[\tau ]_{\Omega }:=\{\sigma \in \Omega :\tau \prec \sigma \}%
\text{.}
\end{equation*}

We define a metric on $\Omega $ by, for $\omega \neq \varkappa $, 
\begin{equation*}
d_{\Omega }(\omega ,\varkappa )=\frac{1}{M^{k}}
\end{equation*}%
if $k$ is the least integer such that $\omega (i)\neq $ $\varkappa (i)$ for
some $i\in T$ with $|i|=k$. Then ($\Omega ,d_{\Omega })$ is a compact metric
space. Furthermore,%
\begin{equation*}
diam(\Omega )=1,
\end{equation*}%
and%
\begin{equation}
diam([\tau ])=\frac{1}{M^{k+1}}  \label{size}
\end{equation}%
whenever $\tau $ is a level-$k$ code tree.

The probabilities $P_{n}$ associated with $\mathcal{F}$ in Equation (\ref%
{superIFS}) induce a natural probability distribution 
\begin{equation*}
\rho \in \mathbb{P}(\Omega )
\end{equation*}%
on $\Omega $. It is defined on cylinder sets $[\tau ]$ by 
\begin{equation}
\rho ([\tau ])=\prod\limits_{1\leq |i|\leq |\tau |}P_{\tau (i)}\text{.}
\label{treeprobs}
\end{equation}%
That is, the random variables $\tau (i)$, with nodal values in $%
\{1,2,...,N\} $, are chosen i.i.d. via $\Pr (\tau (i)=n)=P_{n}$. Then $\rho $
is extended in the usual way to the $\sigma $-algebra $\mathbb{B(}\Omega )$
generated by the cylinder sets. Thus we are able to speak of \textit{the set
of random code trees }$\Omega $\textit{\ with probability distribution }$%
\rho $, and of \textit{selecting trees }$\sigma \in \Omega $\textit{\
according to }$\rho $\textit{.}

A \textit{construction tree} for $\mathcal{F}$ is a code tree wherein the
symbols $1,2,...,$ and $N$ are replaced by the respective IFSs $%
F^{1},F^{2},...,$ and $F^{N}$. A construction tree consists of nodes and
limbs, where each node is labelled by one of the IFSs belonging to $\mathcal{%
F}$. We will associate the $M$ limbs that lie above and meet at a node with
the constituent functions of the IFS of that node; taken in order.

We use the notation $\mathcal{F(}\Omega )$ for the set of construction trees
for $\mathcal{F}$. For $\sigma \in \Omega $ we write $\mathcal{F}(\sigma )$
to denote the corresponding construction tree. We will use the \textit{same
notation} $\mathcal{F}(\sigma )$ to denote the random fractal set associated
with the construction tree $\mathcal{F}(\sigma )$, as described in the next
section.

\subsection{Random Fractals\label{randomsec}}

In this section we describe the canonical random fractal sets and measures
associated with $\mathcal{F}$ in (\ref{superIFS}).

Let $\mathcal{F}$ be given as in (\ref{superIFS}), let $k\in \mathbb{N}$,
and define%
\begin{equation*}
\mathcal{F}_{k}:\Omega \times \mathbb{H(X)\rightarrow H(X)}\text{,}
\end{equation*}%
by 
\begin{equation}
\mathcal{F}_{k}(\sigma )(K\mathbb{)}\mathbb{=}\bigcup\limits_{\{i\in T|\text{
}|i|=k\}}f_{i_{1}}^{\sigma (\emptyset )}\circ f_{i_{2}}^{\sigma
(i_{1})}\circ ...\circ f_{i_{k}}^{\sigma (i_{1}i_{2}...i_{k-1})}(K)
\label{frasets1}
\end{equation}%
$\forall $ $\sigma \in \Omega $ and $K\in $ $\mathbb{H(X)}$. (The set $%
\mathcal{F}_{k}(\sigma )(K\mathbb{)}$ is obtained by taking the union of the
compositions of the functions occurring on the branches of the construction
tree $\mathcal{F}(\sigma )$ starting at the bottom and working up to the $k^{%
\text{th}}$ level, acting upon $K$.)

In a similar way, with measures in place of sets, and unions of sets
replaced by sums of measures weighed by probabilities, we define%
\begin{equation*}
\mathcal{\tilde{F}}_{k}:\Omega \times \mathbb{P(X)\rightarrow P(X)}
\end{equation*}%
by%
\begin{equation}
\mathcal{\tilde{F}}_{k}(\sigma )(\varsigma \mathbb{)}\mathbb{=}%
\sum\limits_{\{i\in T|\text{ }|i|=k\}}(p_{i_{1}}^{\sigma (\emptyset )}\cdot
p_{i_{2}}^{\sigma (i_{1})}\cdot ...\cdot p_{i_{k}}^{\sigma
(i_{1}i_{2}...i_{k-1})})\tilde{f}_{i_{1}}^{\sigma (\emptyset )}\circ \tilde{f%
}_{i_{2}}^{\sigma (i_{1})}\circ ...\circ \tilde{f}_{i_{k}}^{\sigma
(i_{1}i_{2}...i_{k-1})}(\varsigma )  \label{frameas1}
\end{equation}%
$\forall $ $\sigma \in \Omega $ and $\varsigma \in $ $\mathbb{P(X)}$. Note
that the $\mathcal{\tilde{F}}_{k}(\sigma )(\varsigma \mathbb{)}$ all have
unit mass because the $p_{m}^{n}$ sum (over $m$) to unity.

\begin{theorem}
\label{codetree}\label{functions}Let sequences of functions $\{\mathcal{F}%
_{k}\}$ and $\{\mathcal{\tilde{F}}_{k}\}$ be defined as above. Then both the
limits 
\begin{equation*}
\mathcal{F}(\sigma )=\underset{k\rightarrow \infty }{\lim }\{\mathcal{F}%
_{k}(\sigma )(K)\}\text{, and }\mathcal{\tilde{F}}(\sigma )=\underset{%
k\rightarrow \infty }{\lim }\{\mathcal{\tilde{F}}_{k}(\sigma )(\varsigma )\}%
\text{,}
\end{equation*}%
exist, are independent of $K$ and $\varsigma $, and the convergence (in the
Hausdorff and Monge Kantorovitch metrics, respectively,) is uniform in $%
\sigma $, $K$, and $\varsigma $. The resulting functions 
\begin{equation*}
\mathcal{F}:\Omega \mathbb{\rightarrow H(X)}\text{ and }\mathcal{\tilde{F}}%
:\Omega \mathbb{\rightarrow P(X)}
\end{equation*}%
are continuous.
\end{theorem}

\begin{proof}
Make repeated use of the fact that, for fixed $\sigma \in \Omega $, both
mappings are compositions of contraction mappings of contractivity $l$, by
Theorem 1.
\end{proof}

Let 
\begin{equation}
\mathfrak{H}=\{\mathcal{F}(\sigma )\in \mathbb{H(X)}|\sigma \in \Omega \}%
\text{, and }\mathfrak{\tilde{H}}=\{\mathcal{\tilde{F}}(\sigma )\in \mathbb{%
P(X)}|\sigma \in \Omega \}\text{.}  \label{equationH}
\end{equation}%
Similarly let 
\begin{equation}
\mathfrak{P}=\mathcal{F}(\rho )=\rho \circ \mathcal{F}^{-1}\in \mathbb{P(%
\mathfrak{H})}\text{, and }\mathfrak{\tilde{P}}=\mathcal{\tilde{F}}(\rho
)=\rho \circ \mathcal{\tilde{F}}^{-1}\in \mathbb{P(\mathfrak{\tilde{H}})}%
\text{.}  \label{EquationI}
\end{equation}

\begin{definition}
The sets $\mathfrak{H}$ and $\mathfrak{\tilde{H}}$ are called the sets of
fractal sets and fractal measures, respectively, associated with $\mathcal{F}
$. These random fractal sets and measures are said to be distributed
according to $\mathfrak{P}$ and $\mathfrak{\tilde{P}}$, respectively.
\end{definition}

Random fractal sets and measures are hard to compute. There does not appear
to be a general simple forwards (random iteration) algorithm for practical
computation of approximations to them in two-dimensions with affine maps,
for example. The reason for this difficulty lies with the inconvenient
manner in which the shift operator acts on trees $\sigma \in \Omega $
relative to the expressions (\ref{frasets1}) and (\ref{frameas1}).

\begin{definition}
Both the set of IFSs $\{F^{n}:n=1,2,..,N\}$ and the superIFS $\mathcal{F}$
are said to obey the (uniform) open set condition if there exists a
non-empty open set $O$ such that for each $n\in \{1,2,..,N\}$ 
\begin{equation*}
F^{n}(O)=\bigcup\limits_{m=1}^{M}f_{m}^{n}(O)\subset O,
\end{equation*}%
and%
\begin{equation*}
f_{m}^{n}(O)\cap f_{l}^{n}(O)=\varnothing \text{ }\forall \text{ }m\text{, }l%
\text{ }\in \text{ }\{1,2,..,M\}\text{ with }m\neq l.
\end{equation*}
\end{definition}

For the rest of this section we restrict attention to $(\mathbb{X},d)$ where 
$\mathbb{X\subseteq R}^{K}$ and $d$ is the Euclidean metric. The following
theorem gives a specific value for the Hausdorff dimension for almost all of
the random fractal sets in the case of "non-overlapping" similitudes, see\ 
\cite{Fa86}, \cite{Graf} and \cite{MaWi}.

\begin{theorem}
\label{theorem_above} Let the set of IFSs $\{F^{n}:n=1,2,..,N\}$ consist of
similitudes, i.e. $f_{m}^{n}(x)=s_{m}^{n}O_{m}^{n}x+t_{m}^{n}$ where $%
O_{m}^{n}$ is an orthonormal transformation on $\mathbb{R}^{K}$, $%
s_{m}^{n}\in (0,1)$, and $t_{m}^{n}\in \mathbb{R}^{K},$ for all $n\in
\{1,2,...,N\}$ and $m\in \{1,2,...M\}$. Also let $\{F^{n}:n=1,2,..,N\}$ obey
the uniform open set condition. Then for $\mathfrak{P}$ -almost all $A\in 
\mathfrak{H}$%
\begin{equation*}
\dim _{H}A=\dim _{B}A=D
\end{equation*}%
where $D$ is the unique solution of 
\begin{equation*}
\sum\limits_{n=1}^{N}P_{n}\sum_{m=1}^{M}(s_{m}^{n})^{D}=1\text{.}
\end{equation*}
\end{theorem}

\begin{proof}
This is an application of \cite{Falconer} Theorem 15.2, p.230.
\end{proof}

\section{\label{four}Contraction Mappings on Code Trees and the Space $%
\Omega _{V}$}

\subsection{Construction and Properties of $\Omega _{V}$}

Let $V\in \mathbb{N}$. This parameter will describe the \textit{variability }%
of the trees and fractals that we are going to introduce. Let 
\begin{equation*}
\Omega ^{V}=\underset{V\text{ TIMES }}{\underbrace{\Omega \times \Omega
\times ...\times \Omega }}\text{.}
\end{equation*}%
We refer to an element of $\Omega ^{V}$ as a \textit{grove. }In this section
we describe a certain IFS on $\Omega ^{V}$, and discuss its set attractor $%
\Omega _{V}$: its points are ($V$-tuples of) code trees that we will call $V$%
\textit{-groves}. We will find it convenient to label the trunk of each tree
in a grove by its component index, from the set $\{1,2,...,V\}$.

One reason that we are interested in $\Omega _{V}$ is that, as we shall see
later, the set of trees that occur in its components, called $V$-\textit{%
trees,} provides the appropriate code space for a V-variable superfractal.

Next we describe mappings from $\Omega ^{V}$ to $\Omega ^{V}$ that comprise
the IFS. The mappings are denoted by $\eta ^{a}:\Omega ^{V}\rightarrow
\Omega ^{V}$ for $a\in \mathcal{A}$ where%
\begin{equation}
\mathcal{A}:=\{\{1,2,...,N\}\times \{1,2,...,V\}^{M}\}^{V}\text{.}
\label{index}
\end{equation}%
A typical index $a\in \mathcal{A}$ will be denoted%
\begin{equation}
a=(a_{1},a_{2},..,a_{V})  \label{aaaa}
\end{equation}%
where 
\begin{equation}
a_{v}=(n_{v};v_{v,1},v_{v,2},...,v_{v,M})  \label{aavv}
\end{equation}%
where $n_{v}\in \{1,2,...,N\}$ and $v_{v,m}\in \{1,2,...,V\}$ for $m\in
\{1,2,...,M\}$.

Specifically, algebraically, the mapping $\eta ^{a}$ is defined in Equation (%
\ref{Etas}) below. But it is very useful to represent the indices and the
mappings with trees. See Figure \ref{vtrees1}. Each map $\eta ^{a}$ in
Equation (\ref{Etas}) and each index $a\in \mathcal{A}$ may be represented
by a $V$-tuple of labelled level-1 trees that we call \textit{(level-1)
function trees}. Each function tree has a trunk, a node, and $M$ limbs.
There is one function tree for each component of the mapping. Its trunk is
labelled by the index of the component $v\in \{1,2,...,V\}$ to which it
corresponds. The node of each function tree is labelled by the IFS number $%
n_{v}$(shown circled) of the corresponding component of the mapping. The $%
m^{th}$ limb of the $v^{th}$ tree is labelled by the number $v_{v,m}\in
\{1,2,...,V\}$, for $m\in \{1,2,...,M\}$. We will use the same notation $%
\eta ^{a}$ to denote both a $V$-tuple of function trees and the unique
mapping $\eta ^{a}:\Omega ^{V}\rightarrow \Omega ^{V}$ to which it
bijectively corresponds. We will use the notation $a$ to denote both a $V$%
-tuple of function trees and the unique index $a\in \mathcal{A}$ to which it
bijectively corresponds.\FRAME{ftbpFU}{4.7556in}{1.6881in}{0pt}{\Qcb{Each
map $\protect\eta ^{a}$ in Equation (\protect\ref{Etas}) and each index $%
a\in \mathcal{A}$ may be represented by a $V$-tuple of labelled level-1
trees that we call \textit{(level-1) function trees}. Each function tree has
a trunk, a node, and $M$ limbs. The function trees correspond to the
components of the mapping. The trunk of each function tree is labelled by
the index of the component $v\in \{1,2,...,V\}$ to which it corresponds. The
node of each function tree is labelled by the IFS number $n_{v}$(shown
circled) of the corresponding component of the mapping. The $m^{th}$ limb of
the $v^{th}$ tree is labelled by the domain number $v_{v,m}\in \{1,2,...,V\}$%
, for $m\in \{1,2,...M\}$.}}{\Qlb{vtrees1}}{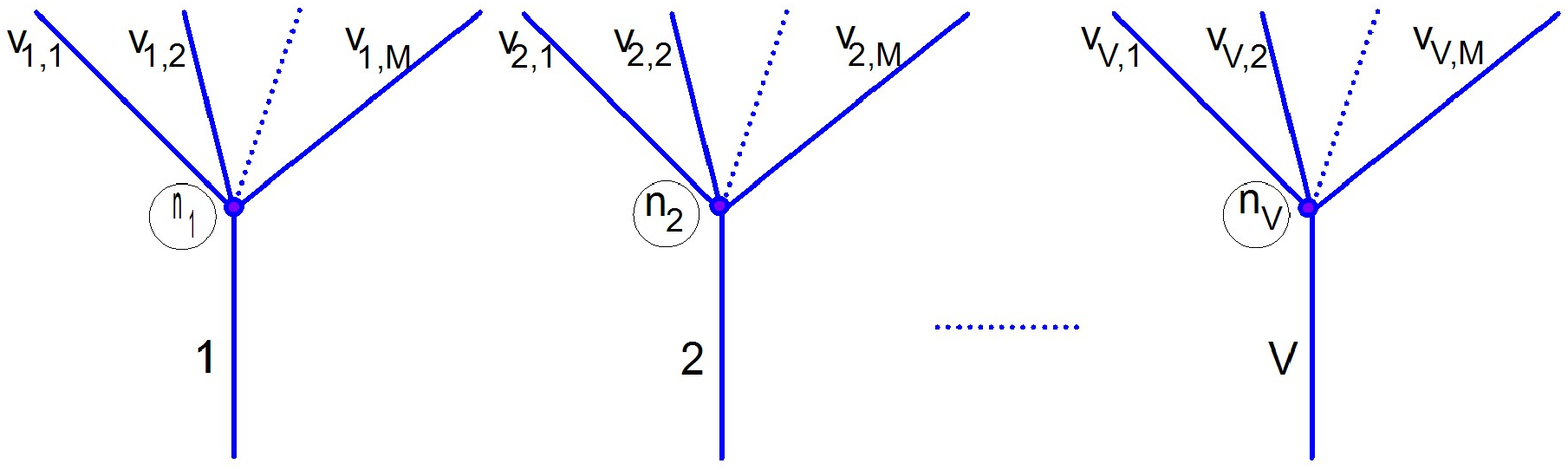}{\raisebox{-1.6881in}{\includegraphics[height=1.6881in]{vtrees1.ps}}}

Now we can describe the action of $\eta ^{a}$ on a grove $\omega \in \Omega
^{V}$. We illustrate with an example where $V=3$, $N=5$, and $M=2$. In
Figure \ref{col_tree1} an arbitrary grove $\omega =(\omega _{1},\omega
_{2},\omega _{3})\in \Omega ^{3}$ is represented by a triple of coloured
tree pictures, one blue, one orange, and one magenta, with trunks labelled
one, two, and \ three respectively. The top left of Figure \ref{col_tree1}
shows the map $\eta ^{a}$, and the index $a\in \mathcal{A}$, where%
\begin{equation}
\eta ^{a}(\omega _{1},\omega _{2},\omega _{3})=(\xi _{1}(\omega _{1},\omega
_{2}),\xi _{5}(\omega _{3},\omega _{2}),\xi _{4}(\omega _{3},\omega _{1}))%
\text{,}  \label{exmap}
\end{equation}%
and 
\begin{equation*}
a=((1;1,2),(5;3,2),(4;3,1))\text{,}
\end{equation*}%
represented by function trees. The functions $\{\xi _{n}:n=1,2,...,5\}$ are
defined in Equation (\ref{definition_jamie}) below. The result of the action
of $\eta ^{a}$ on $\omega $ is represented, in the bottom part of Figure \ref%
{col_tree1}, by a grove whose lowest nodes are labelled by the IFS numbers
1, 5, and 4, respectively, and whose subtrees at level zero consist of trees
from $\omega $ located according to the limb labels on the function trees.
(Limb labels of the top left expression, the function trees of $\eta ^{a}$,
are matched to trunk labels in the top right expression, the components of $%
\omega $.) In general, the result of the action of $\eta ^{a}$ in Figure \ref%
{vtrees1} on a grove $\omega \in \Omega ^{V}$ (represented by $V$ trees with
trunks labelled from $1$ to $N$) is obtained by matching the limbs of the
function trees to the trunks of the $V$-trees, in a similar manner.\FRAME{%
ftbpFU}{4.8058in}{2.8323in}{0pt}{\Qcb{The top right portion of the image
represents a grove $\protect\omega =(\protect\omega _{1},\protect\omega _{2},%
\protect\omega _{3})\in \Omega ^{3}$ by a triple of coloured tree pictures.
The top left portion represents the map $\protect\eta ^{a}$ in Equation (%
\protect\ref{exmap}) using level-1 function trees. The bottom portion
represents the image $\protect\eta ^{a}(\protect\omega )$.}}{\Qlb{col_tree1}%
}{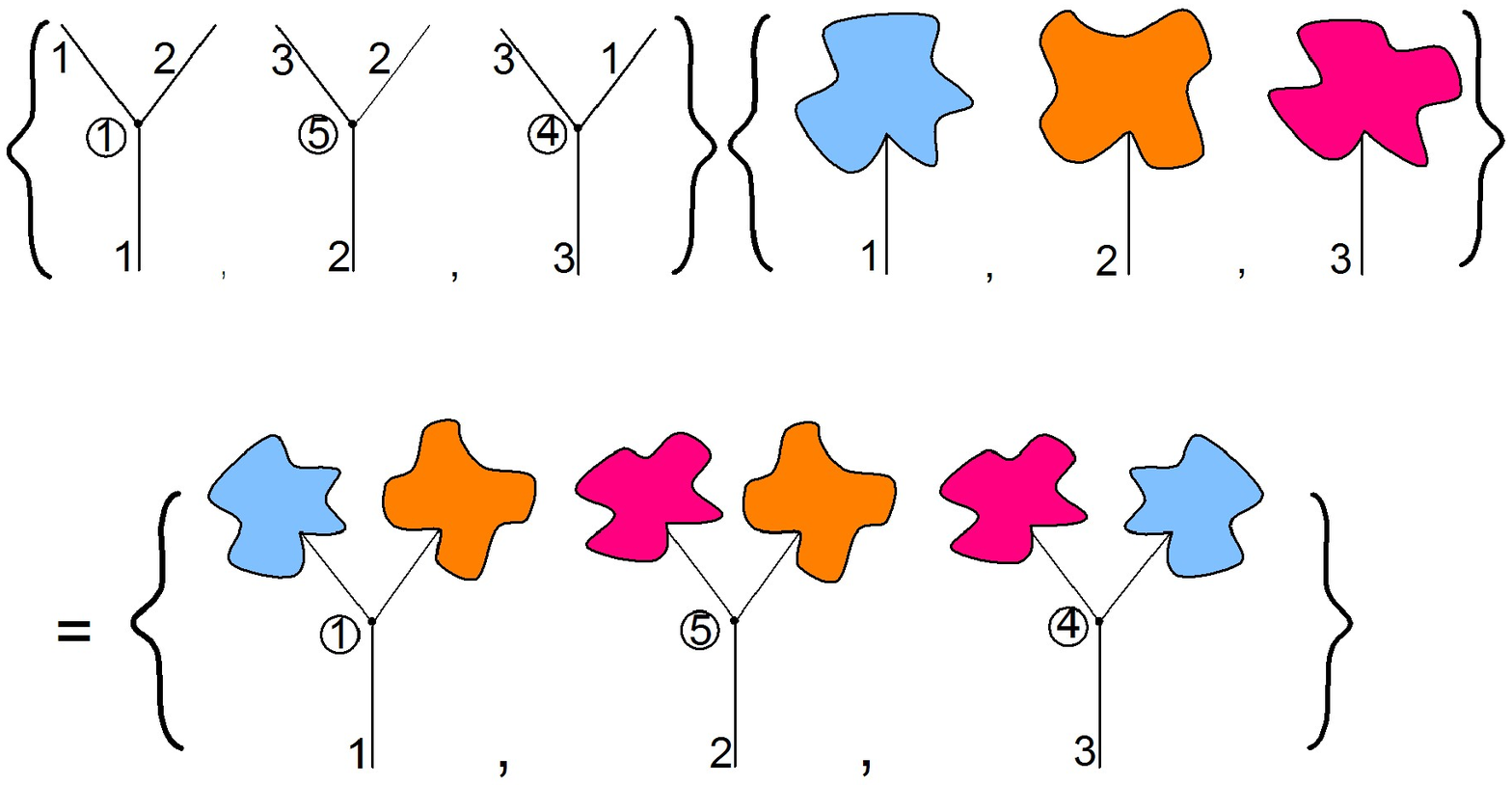}{\raisebox{-2.8323in}{\includegraphics[height=2.8323in]{col_trees1.ps}}}

We are also going to need a set of probabilities $\{\mathcal{P}^{a}|a\in 
\mathcal{A}\}$, with 
\begin{equation}
\sum\limits_{a\in \mathcal{A}}\mathcal{P}^{a}=1\text{, }\mathcal{P}^{a}\geq 0%
\text{ }\forall a\in \mathcal{A}\text{.}  \label{probs}
\end{equation}%
These probabilities may be more or less complicated. Some of our results are
specifically restricted to the case 
\begin{equation}
\mathcal{P}^{a}=\frac{P_{n_{1}}P_{n_{2}}...P_{n_{V}}}{V^{MV}}\text{,}
\label{SuperProbs}
\end{equation}%
which uses only the set of probabilities $\{P_{1},P_{2},...,P_{N}\}$
belonging to the superIFS (\ref{superIFS}). This case corresponds to
labelling all nodes and limbs in a function tree independently with
probabilities such that limbs are labeled according to the uniform
distribution on $\{1,2,...,V\}$, and nodes are labelled $j$ with probability 
$P_{j}$.

or each $n\in \{1,2,...,N\}$ define the $n^{th}$ shift mapping $\xi
_{n}:\Omega ^{M}\rightarrow \Omega $ by 
\begin{equation}
(\xi _{n}(\omega ))(\varnothing )=n\text{ and (}\xi _{n}(\omega
))(mi)=\omega _{m}(i)\text{ }\forall \text{ }i\in T,m\in \{1,2,...,M\}\text{,%
}  \label{definition_jamie}
\end{equation}%
for $\omega =(\omega _{1},\omega _{2},...,\omega _{M})\in \Omega ^{M}$. That
is, the mapping $\xi _{n}$ creates a code tree with its bottom node labelled 
$n$ attached directly to the M trees $\omega _{1},\omega _{2},...,\omega
_{M} $.

\begin{theorem}
\label{jamie} For each $a\in \mathcal{A}$ define $\eta ^{a}:\Omega
^{V}\rightarrow \Omega ^{V}$by%
\begin{equation}
\eta ^{a}(\omega _{1},\omega _{2},...,\omega _{V}):=(\xi _{n_{1}}(\omega
_{v_{1,1}},\omega _{v_{1,2}},...,\omega _{v_{1,M}}),  \label{Etas}
\end{equation}%
\begin{equation*}
\xi _{n_{2}}(\omega _{v_{2,1}},\omega _{v_{2,2}},...,\omega
_{v_{2,M}}),...,\xi _{n_{V}}(\omega _{v_{V,1}},\omega _{v_{V,2}},...,\omega
_{v_{V,M}}))
\end{equation*}%
Then $\eta ^{a}:\Omega ^{V}\rightarrow \Omega ^{V}$ is a contraction mapping
with Lipshitz constant $\frac{1}{M}$.
\end{theorem}

\begin{proof}
Let $a\in \mathcal{A}$. Let $(\omega _{1},\omega _{2},...,\omega _{V})$ and $%
(\widetilde{\omega }_{1,}\widetilde{\omega }_{2,}...,\widetilde{\omega }%
_{V}) $ be any pair of points in $\Omega ^{V}$. Then%
\begin{equation*}
d_{\Omega ^{V}}(\eta ^{a}(\omega _{1},\omega _{2},...,\omega _{V}),\eta ^{a}(%
\widetilde{\omega }_{1,}\widetilde{\omega }_{2,}...,\widetilde{\omega }%
_{V}))=
\end{equation*}%
\begin{equation*}
\max_{v\in \{1,2,...,V\}}d_{\Omega }(\xi _{n_{v}}(\omega _{v_{v,1}},\omega
_{v_{v,2}},...,\omega _{v_{v,M}}),\xi _{n_{v}}(\widetilde{\omega }_{v_{v,1}},%
\widetilde{\omega }_{v_{v,2}},...,\widetilde{\omega }_{v_{v,M}}))
\end{equation*}%
\begin{equation*}
\leq \frac{1}{M}\cdot \max_{v\in \{1,2,...,V\}}d_{\Omega ^{M}}((\omega
_{v_{v,1}},\omega _{v_{v,2}},...,\omega _{v_{v,M}}),(\widetilde{\omega }%
_{v_{v,1},}\widetilde{\omega }_{v_{v,2}}...,\widetilde{\omega }_{v_{v,M}}))
\end{equation*}%
\begin{equation*}
\leq \frac{1}{M}\cdot \max_{v\in \{1,2,...,V\}}d_{\Omega }(\omega _{v},%
\widetilde{\omega }_{v})=\frac{1}{M}\cdot d_{\Omega ^{V}}((\omega
_{1},\omega _{2},...,\omega _{V}),(\widetilde{\omega }_{1,}\widetilde{\omega 
}_{2,}...,\widetilde{\omega }_{V})).
\end{equation*}
\end{proof}

It follows that we can define an IFS $\Phi $ of strictly contractive maps by 
\begin{equation}
\Phi =\{\Omega ^{V};\eta ^{a},\mathcal{P}^{a},a\in \mathcal{A}\}.
\label{specialIFS}
\end{equation}%
Let the set attractor and the measure attractor of the IFS $\Phi $ be
denoted by $\Omega _{V}$, and $\mu _{V}$ respectively. Clearly, $\Omega
_{V}\in \mathbb{H(}\Omega ^{V})$ while $\mu _{V}$ $\in \mathbb{P(}\Omega ^{V}%
\mathbb{)}$. We call $\Omega _{V}$ the \textit{space of }$V$\textit{-groves}%
. The elements of $\Omega _{V}$ are certain $V$-tuples of $M$-fold code
trees on an alphabet of $N$ symbols, which we characterize in Theorem \ref%
{Groves}. But we think of them as special groves of special trees, namely $V$%
-trees.

For all $v\in \{1,2,...,V\}$, let us define $\Omega _{V,v}\subset \Omega
_{V} $ to be the set of $v^{\text{th}}$ components of groves in $\Omega _{V}$%
. Also let $\rho _{V}\in \mathbb{P}(\Omega )$ denote the marginal
probability measure defined by 
\begin{equation}
\rho _{V}(B):=\mu _{V}(B,\Omega ,\Omega ,...,\Omega )\forall B\in \mathbb{B}%
(\Omega \mathbb{)}\text{.}  \label{marginal}
\end{equation}

\begin{theorem}
\label{vtrees}For all $v\in \{1,2,...,V\}$ we have%
\begin{equation}
\Omega _{V,v}=\Omega _{V,1}:=\{\text{set of all V-trees}\}.  \label{setofall}
\end{equation}%
When the probabilities $\{\mathcal{P}^{a}|a\in \mathcal{A}\}$ obey Equation (%
\ref{SuperProbs}), then starting at any initial grove, the random
distribution of trees $\omega \in \Omega $ that occur in the $v^{\text{th}}$
components of groves\ produced by the random iteration algorithm
corresponding to the IFS $\Phi $, after n iteration steps, converge weakly
to $\rho _{V}$ independently of $v$, almost always, as $n\rightarrow \infty $%
.
\end{theorem}

\begin{proof}
Let $\Xi :\Omega $ $^{V}\rightarrow \Omega $ $^{V}$ denote any map that
permutes the coordinates. Then the IFS $\Phi =\{\Omega ^{V};\eta ^{a},%
\mathcal{P}^{a},a\in \mathcal{A}\}$ is invariant under $\Xi $, that is $\Phi
=\{\Omega ^{V};\Xi \eta ^{a}\Xi ^{-1},\mathcal{P}^{a},a\in \mathcal{A}\}$.
It follows that $\Xi \Omega _{V}=\Omega _{V}$ and $\Xi \mu _{V}=\mu _{V}$.
It follows that Equation (\ref{setofall}) holds, and also that, in the
obvious notation, for any $(B_{1},B_{2},...,B_{V})\in (\mathbb{B}(\Omega 
\mathbb{))}^{V}$ we have 
\begin{equation}
\mu _{V}(B_{1},B_{2},...,B_{V})=\mu _{V}(B_{\sigma _{1}},B_{\sigma
_{2}},...,B_{\sigma _{V}}).  \label{assertion}
\end{equation}%
In particular 
\begin{equation*}
\rho _{V}(B)=\mu _{V}(B,\Omega ,\Omega ,...,\Omega )=\mu _{V}(\Omega ,\Omega
,...,\Omega ,B,\Omega ,...,\Omega )\text{ }\forall B\in \mathbb{B}(\Omega 
\mathbb{)}\text{,}
\end{equation*}%
where the $``B$\textquotedblright\ on the right-hand-side is in the $v^{th}$
position. Theorem \ref{jamie} tells us that we can apply the random
iteration algorithm (Theorem \ref{RIA}) to the IFS $\Phi $. This yields
sequences of measures, denoted by $\{\mu _{V}^{(l)}:l=1,2,3...\}$, that
converge weakly to $\mu _{V}$ almost always. In particular $\mu
_{V}^{(l)}(B,\Omega ,\Omega ,...,\Omega )$ converges to $\rho _{V}$ almost
always.
\end{proof}

Let $L$ denote a set of fixed length vectors of labelled trees. We will say
that $L$ and its elements have the property of $V$-\textit{variability}, or
that $L$ and its elements are\textit{\ }$V$-\textit{variable,} if and only
if, for all $\omega \in L$, the number of distinct subtrees of all
components of $\omega $ at level $k$ is at most $V$, for each level $k$ of
the trees.

\begin{theorem}
\label{Groves} Let $\omega \in \Omega ^{V}$. Then $\omega \in \Omega _{V}$
if and only $\omega $ contains at most $V$ distinct subtrees at any level
(i.e. $\omega $ is V-variable). Also, if $\sigma \in \Omega ,$ then $\sigma $
is a V-tree if and only if it is V-variable.
\end{theorem}

\begin{proof}
Let 
\begin{equation*}
S=\{\omega \in \Omega ^{V}\text{ }|\text{ \ }|\{\text{subtrees of components
of }\omega \text{ at level }k\}|\leq V,\text{ }\forall \text{ }k\in \mathbb{%
N\}}\text{.}
\end{equation*}%
Then $S$ is closed: Let $\{s_{n}\in S\}$ converge to $s^{\ast }$. Suppose
that $s^{\ast }\notin S$. Then, at some level $k\in \mathbb{N}$, $s^{\ast }$
has more than $V$ subtrees. There exists $l\in \mathbb{N}$ so that each
distinct pair of these subtrees of $s^{\ast }$ first disagrees at some level
less than $l$. Now choose $n$ so large that $s_{n}$ agrees with $s^{\ast }$
through level $(k+l)$ (i.e. $d_{\Omega ^{V}}(s_{n},s^{\ast })<\frac{1}{%
M^{(k+l)}}$ ). Then $s_{n}$ has more than $V$ distinct subtrees at level $k$%
, a contradiction. So $s^{\ast }\in S$.

Also $S$ is non-empty: Let $\sigma \in \Omega $ be defined by $\sigma (i)=1$
for all $i\in T$. Then $(\sigma ,\sigma ,...,\sigma )\in S$.

Also $S$ is invariant under the IFS $\Phi $: Clearly any $s\in S$ can be
written $s=\eta ^{a}(\widetilde{s})$ for some $a\in \mathcal{A}$ and $%
\widetilde{s}\in S$. Also, if $s\in S$ then $\eta ^{a}(s)\in S$. So $S=\cup
\{\eta ^{a}(S)$ $|$ $a\in \mathcal{A}\}$.

Hence, by uniqueness, $S$ must be \textit{the} set attractor of $\Phi $.
That is, $S=\Omega _{V}$. This proves the first claim in the Theorem.

It now follows that if $\sigma \in \Omega $ is a $V$-tree then it contains
at most $V$ distinct subtrees at level $k$, for each $k\in \mathbb{N}$.
Conversely, it also follows that if $\sigma \in \Omega $ has the latter
property, then $(\sigma ,\sigma ,...,\sigma )\in \Omega _{V}$, and so $%
\sigma \in \Omega _{V,1}$.
\end{proof}

\begin{theorem}
\label{vinfinity} For all 
\begin{equation*}
d_{\mathbb{H}(\Omega \mathbb{)}}(\Omega _{V,1},\Omega )<\frac{1}{V}\text{,}
\end{equation*}%
which implies that 
\begin{equation*}
\underset{V\rightarrow \infty }{\lim }\Omega _{V,1}=\Omega \text{,}
\end{equation*}%
where the convergence is in the Hausdorff metric. Let the probabilities $\{%
\mathcal{P}^{a}|a\in \mathcal{A}\}$ obey Equation (\ref{SuperProbs}). Then 
\begin{equation}
d_{\mathbb{P}(\Omega )}(\rho _{V},\rho )\leq 1.4\left( \frac{M}{V}\right) ^{%
\frac{1}{4}}  \label{tobeproved}
\end{equation}%
which implies 
\begin{equation*}
\underset{V\rightarrow \infty }{\lim }\rho _{V}=\rho \text{,}
\end{equation*}%
where $\rho $ is the stationary measure on trees introduced in Section \ref%
{treesec}, and convergence is in the Monge Kantorovitch metric.
\end{theorem}

\begin{proof}
To prove the first part, let $M^{k+1}>V\geq M^{k}$. Let $\tau $ be any level-%
$k$ code tree. Then $\tau $ is clearly $V$-variable and it can be extended
to an infinite $V$-variable code tree. It follows that $[\tau ]\cap \Omega
_{V,1}\neq \varnothing $. The collection of such cylinder sets $[\tau ]$
forms a disjoint partition of $\Omega $ by subsets of diameter $\frac{1}{%
M^{k+1}}$ , see Equation (\ref{size})), from which it follows that 
\begin{equation*}
d_{\mathbb{H}(\Omega \mathbb{)}}(\Omega _{V,1},\Omega )\leq \frac{1}{M^{k+1}}%
<\frac{1}{V}\text{.}
\end{equation*}%
The first part of the theorem follows at once.

For the proof of the second part, we refer to Section \ref{proofsec}.
\end{proof}

Let $\Sigma _{V}=\mathcal{A}^{\infty }.$ This is simply the code space
corresponding to the IFS $\Phi $ defined in Equation (\ref{specialIFS}).
From Theorem \ref{Add} there exists a continuous onto mapping $\Phi :\Sigma
_{V}\rightarrow $ $\Omega _{V}$ defined by%
\begin{equation*}
\Phi (a_{1}a_{2}a_{3}...)=\underset{k\rightarrow \infty }{\lim }\eta
^{a_{1}}\circ \eta ^{a_{2}}\circ ...\circ \eta ^{a_{k}}(\omega )
\end{equation*}%
for all $a_{1}a_{2}a_{3}...\in \Sigma _{V}$, for any $\omega $. In the
terminology of section \ref{codespace} the sequence $a_{1}a_{2}a_{3}..\in
\Sigma _{V}$ is an address of the $V$-grove $\Phi (a_{1}a_{2}a_{3}...)\in
\Omega _{V}$ and $\Sigma _{V}$ is the code space for the set of V-groves $%
\Omega _{V}$. In general $\Phi :\Sigma _{V}\rightarrow $ $\Omega _{V}$ is
not one-to-one, as we will see in Section \ref{backwards}.

\subsection{\label{backwards}Compositions of the Mappings $\protect\eta ^{a}$%
}

Compositions of the mappings $\eta ^{a}:\Omega ^{V}\rightarrow \Omega ^{V},$ 
$a\in \mathcal{A}$, represented by $V$-tuples of level-$1$ function trees,
as in Figure \ref{vtrees1}, can be represented by higher level trees that we
call \textit{level-}$k$ \textit{function trees}.

First we illustrate the ideas, then we formalize. In Figure \ref{col_tree4}
we illustate the idea of composing $V$-tuples of function trees. In this
example $V=3$, $N=5$, and $M=2$. The top row shows the level-1 function
trees corrresponding to $a,b,c\in \mathcal{A}$ given by%
\begin{equation*}
a=((1;2,3),(3;1,3),(5;2,3))\text{,}
\end{equation*}%
\begin{equation*}
b=((4;3,1),(2;1,2),(3;3,2))\text{,}
\end{equation*}%
and%
\begin{equation*}
c=((1;1,2),(5;1,3),(4;2,3))\text{.}
\end{equation*}%
The first entry in the second row shows the $3$-tuple of level-$2$ function
trees $a\circ b$. The bottom bracketed expression shows the $3$-tuple of
level-$2$ function trees $a\circ b\circ c$.\FRAME{ftbpFU}{4.8058in}{2.9654in%
}{0pt}{\Qcb{Illustrations of compositions of function trees to produce
higher level function trees. Here $V=3$, $N=5$, and $M=2$. We compose the
level-1 function trees corresponding to $a=((1;2,3),(3;1,3),(5;2,3))$, $%
b=((4;3,1),(2;1,2),(3;3,2))$, and $c=((1;1,2),(5;1,3),(4;2,3))$. The top row
shows the separate level-1 function trees, $a,b,$ and $c.$ The second row
shows the level-$2$ function tree $a\circ b,$ and the function tree $c$. The
last row shows the level-$3$ function tree $a\circ b\circ c$. }}{\Qlb{%
col_tree4}}{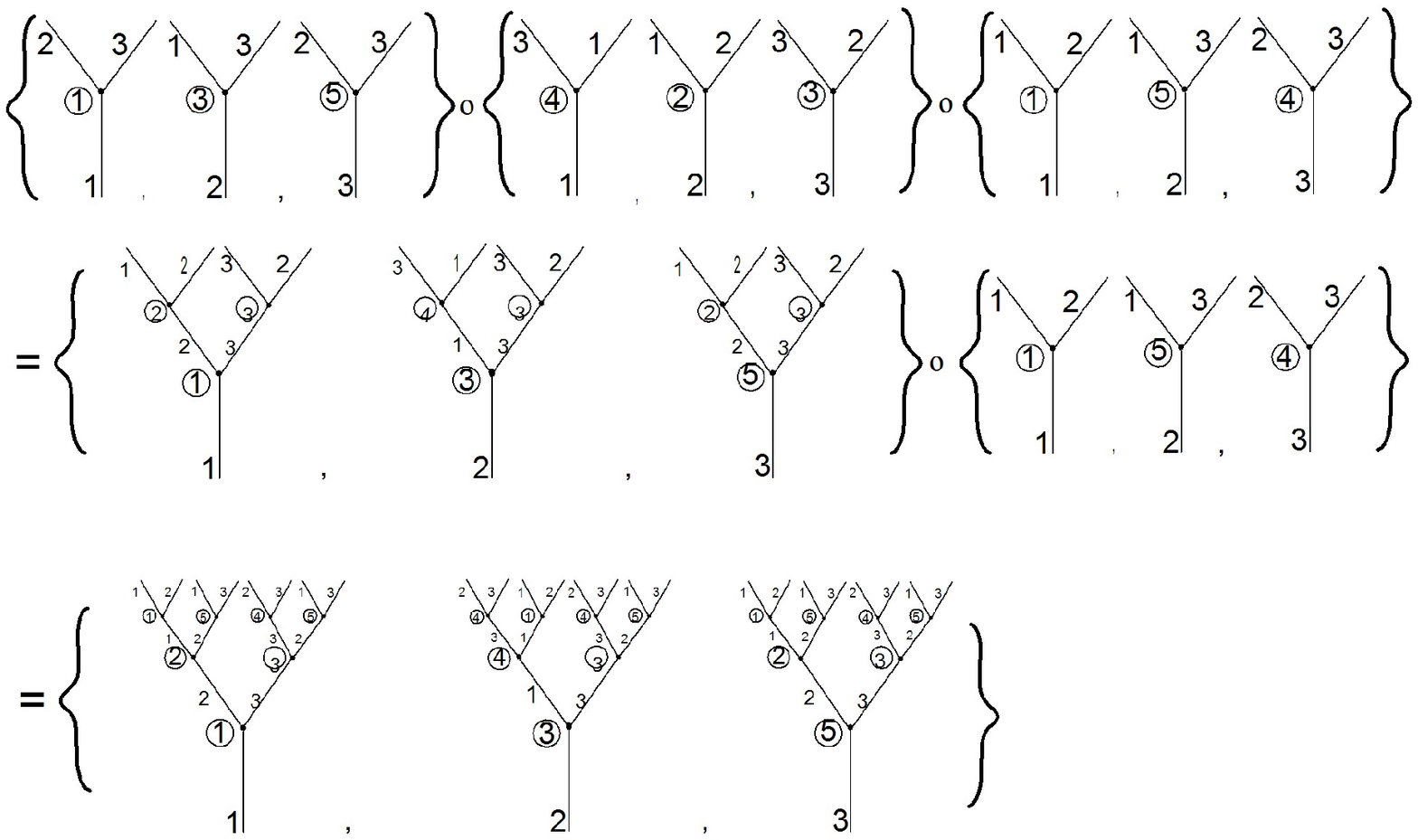}{\raisebox{-2.9654in}{\includegraphics[height=2.9654in]{col_tree4.ps}}}Then in Figure %
\ref{col_tree5}, we have represented $\eta ^{a}\circ \eta ^{b}\circ \eta
^{c}(\omega )$ for $\omega =(\omega _{1},\omega _{2},\omega _{3})\in \Omega
^{3}$. The key ideas are (i) $a\circ b\circ c$ can be converted into a
mapping $\eta ^{a\circ b\circ c}:\Omega ^{3}\rightarrow \Omega ^{3}$and (ii)%
\begin{equation*}
\eta ^{a\circ b\circ c}=\eta ^{a}\circ \eta ^{b}\circ \eta ^{c}\text{.}
\end{equation*}%
The mapping $\eta ^{a\circ b\circ c}(\omega )$ is defined to be the $3$%
-tuple of code trees obtained by attatching the tree $\omega _{v}$ to each
of the top limbs of each level-$3$ function tree in $a\circ b\circ c$ with
label $v$ for all $v\in \{1,2,3\}$ then dropping all of the labels on the
limbs.\FRAME{ftbpFU}{4.7556in}{2.2459in}{0pt}{\Qcb{Illustration of the
composition of the three mappings, $\protect\eta ^{a}\circ $ $\protect\eta %
^{b}\circ \protect\eta ^{c}=$ $\protect\eta ^{a\circ b\circ c}$ applied to $%
\protect\omega =(\protect\omega _{1},\protect\omega _{2},\protect\omega %
_{3})\in \Omega ^{3}$.\ See also Figure \protect\ref{col_tree4} and the text.%
}}{\Qlb{col_tree5}}{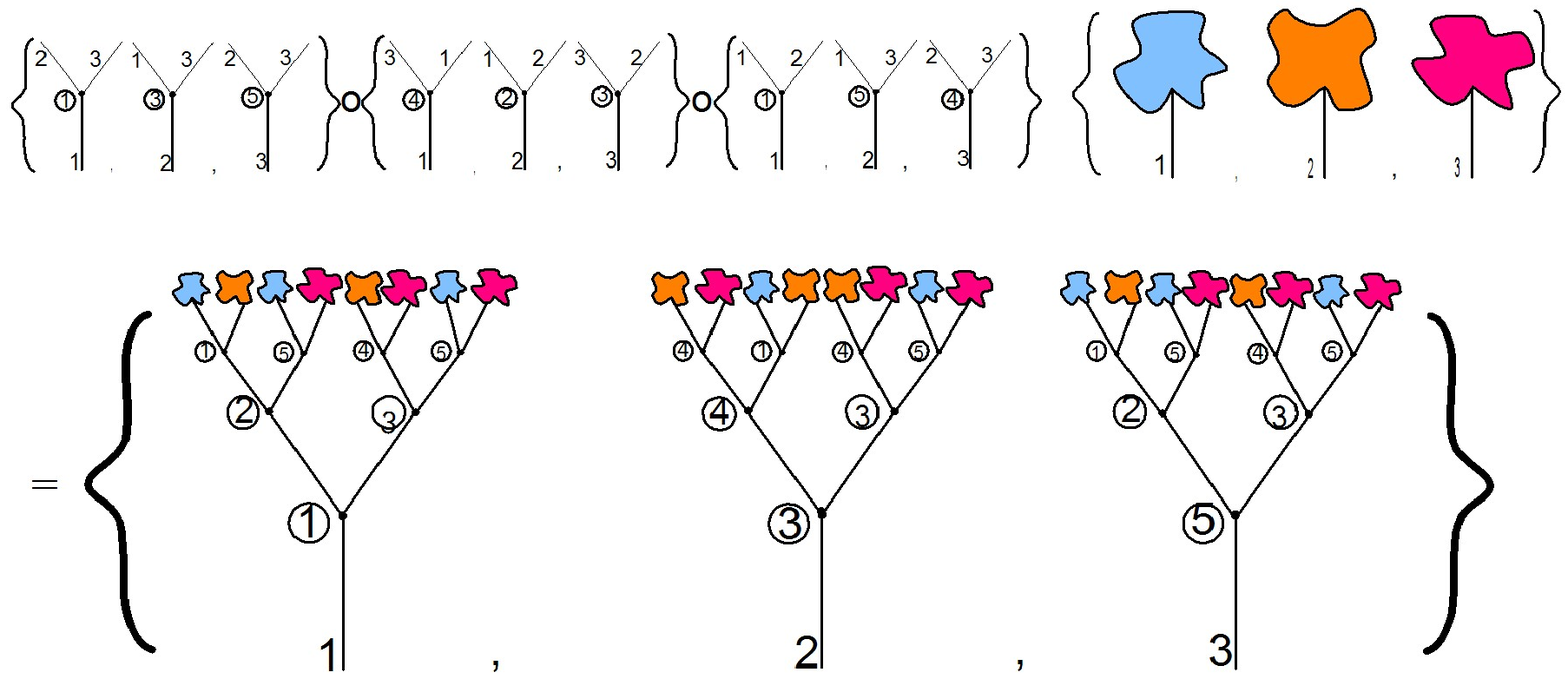}{\raisebox{-2.2459in}{\includegraphics[height=2.2459in]{col-tree5.ps}}}

Next we formalize. Let $k\in \mathbb{N}$. Define a \textit{level-}$k$\textit{%
\ function tree }to be a level-$k$ labelled tree with the nodes of the first 
$k-1$ levels labelled from $\{1,2,...,N\}$ and limbs (i.e. second nodal
values) of all $k$ levels labelled from $\{1,2,...,V\}$. We define a \textit{%
grove} of level-$k$ function trees, say $g$, to be a $V$-tuple of level-$k$
function trees, with trunks labelled according to the component number, and
we define $G_{k}$ to be the set of such $g$. Let $G:=\tbigcup\limits_{k\in 
\mathbb{N}}G_{k}$. We will refer to a component of an element of $G$ simply
as a \textit{function tree. }For $g\in G$ we will write $\left\vert
g\right\vert =k$, where $k\in \mathbb{N}$ is the unique number such that $%
g\in G_{k}$.

Then, for all $g=(g_{1},g_{2},...,g_{V})\in G$, for all $v\in \{1,2,...,V\}$,%
\begin{equation*}
g_{v}(\text{\textit{node} }i)\in \{1,2,...,N\}\text{ }\forall \text{ }i\in T%
\text{ with }|i|\leq \left\vert g\right\vert -1,
\end{equation*}%
and%
\begin{equation*}
g_{v}(\text{\textit{limb} }i)\in \{1,2,...,V\}\text{ }\forall \text{ }i\in T%
\text{ with }|i|\leq \left\vert g\right\vert \text{.}
\end{equation*}

For all $g,h\in G$ we define the composition $g\circ h$ to be the grove of $%
(\left\vert g\right\vert +\left\vert h\right\vert )$-level function trees
given by the following expressions.%
\begin{equation*}
(g\circ h)_{v}(\text{\textit{node} }i)=g_{v}(\text{\textit{node }}i)\text{ }%
\forall \text{ }i\in T\text{ with }|i|\leq \left\vert g\right\vert -1;
\end{equation*}%
\begin{equation*}
(g\circ h)_{v}(\text{\textit{limb} }i)=g_{v}(\text{\textit{limb }}i)\text{ }%
\forall \text{ }i\in T\text{ with }|i|\leq \left\vert g\right\vert ;
\end{equation*}%
\begin{equation*}
(g\circ h)_{v}(\text{\textit{node }}ij)=h_{g_{v}(\text{\textit{limb }}i)}(%
\text{\textit{node}}\mathit{\ }j)\text{ }\forall \text{ }ij\in T\text{ with }%
|i|=\left\vert g\right\vert ,\text{ }|j|\leq \left\vert h\right\vert -1;
\end{equation*}%
\begin{equation*}
(g\circ h)_{v}(\text{\textit{limb} }ij)=h_{g_{v}(\text{\textit{limb }}i)}(%
\text{\textit{limb }}j)\text{ }\forall \text{ }ij\in T\text{ with }%
|i|=\left\vert g\right\vert ,\text{ }|j|\leq \left\vert h\right\vert \text{.}
\end{equation*}

For all $g\in G$ we define $\eta ^{g}:\Omega ^{V}\rightarrow \Omega ^{V}$ by 
\begin{equation*}
(\eta ^{g}(\omega ))_{v}(i)=g_{v}(\text{\textit{node} }i)\text{ }\forall 
\text{ }i\in T\text{ with }|i|\leq \left\vert g\right\vert -1\text{,}
\end{equation*}%
and 
\begin{equation*}
(\eta ^{g}(\omega ))_{v}(ij)=\omega _{g_{v}(\text{\textit{limb} }i)}(\text{%
\textit{node} }j)\text{ }\forall \text{ }ij\in T\text{ with }|i|=\left\vert
g\right\vert ,\text{ }|j|\geq 0\text{.}
\end{equation*}%
This is consistent with the definition of $\eta ^{a}:\Omega ^{V}\rightarrow
\Omega ^{V}$ for $a\in \mathcal{A}$, as the following theorem shows. We will
write $\eta ^{g}$ to denote both the mapping $\eta ^{g}:\Omega
_{V}\rightarrow \Omega _{V}$ and the corresponding unique $V$-tuple of level-%
$k$ code trees for all $g\in G$.

\begin{theorem}
For all $g,h\in G_{k}$ we have%
\begin{equation}
\eta ^{g\circ h}=\eta ^{g}\circ \eta ^{h}\text{.}  \label{compositionone}
\end{equation}%
It follows that the operation $\circ $ between ordered pairs of elements of $%
G$ is associative. In particular, if $(a_{1},a_{2},...,a_{k})\in \mathcal{A}%
^{k}$ then%
\begin{equation}
\eta ^{a_{1}}\circ \eta ^{a_{2}}\circ ...\circ \eta ^{a_{k}}=\eta
^{a_{1}\circ a_{2}\circ ...\circ a_{k}}\text{.}  \label{composition}
\end{equation}
\end{theorem}

\begin{proof}
This is a straightforward exercise in substitutions and is omitted here.
\end{proof}

We remark that Equations (\ref{composition}) and (\ref{compositionone})
allow us to work directly with function trees to construct, count, and track
compositions of mappings $\eta ^{a}$ $\in \mathcal{A}$. The space $G$ also
provides a convenient setting for contrasting the forwards and backwards
algorithms associated with the IFS $\Phi $. For example, by composing
function trees in such as way as to build from the bottom up, which
corresponds to a backwards algorithm, we find that we can construct a
sequence of cylinder set approximations to the first component $\omega _{1}$
of $\omega \in \Omega _{V}$ without having to compute approximations to the
other components.

Let $G_{k}(V)\subseteq G_{k}$ denote the set of elements of $G_{k}$ that can
be written in the form $a_{1}\circ a_{2}\circ ...\circ a_{k}$ for some $%
\mathbf{a}=(a_{1},a_{2},...,a_{k})\in \mathcal{A}^{k}$. (We remark that $%
G_{k}(V)$ is $V$-variable by a similar argument to the proof of Theorem \ref%
{Groves}.) Then we are able to estimate the measures of the cylinder sets $%
[(\eta ^{a_{1}\circ a_{2}\circ ...\circ a_{k}})_{1}]$ by computing the
probabilities of occurence of function trees $g\in G_{k}(V)$ such that ($%
\eta ^{g})_{1}=(\eta ^{a_{1}\circ a_{2}\circ ...\circ a_{k}})_{1}$ built up
starting from level-0 trees, with probabilities given by Equation (\ref%
{SuperProbs}), as we do in the Section \ref{proofsec}.

The labeling of limbs in the approximating grove of function trees of level-$%
k$ of $\Phi (a_{1}a_{2}.....)$ defines the basic $V$-Variable dependence
structure of $\Phi (a_{1}a_{2}.....)$. We call the code tree of limbs of a
function tree the associated dependence tree.

The grove of code trees for $\Phi(a_1 a_2 .....)$ is by definition totally
determined by the labels of the nodes. Nevertheless its grove of dependence
trees contains all information concerning its V-Variable structure. The
dependence tree is the characterizing skeleton of $V$-Variable fractals.

\subsection{A direct characterization of the measure $\protect\rho_V$}

Let$\{\mathbf{F}_{n}\}_{n=0}^{\infty
}=\{F_{1}^{n},...,F_{V}^{n}\}_{n=0}^{\infty }$ be a sequence of random
groves of level-$1$ function trees. Each random function tree can be
expressed as $F_{v}^{n}=(N_{v}^{n},L_{v}^{n}(1),...,L_{v}^{n}(M))$ where the 
$N_{v}^{n}$'s and $L_{v}^{n}$'s corresponds to random labellings of nodes
and limbs respectively.

We assume that the function trees $\{ F_v^n \}$, are independent with $Pr(
N_v^n = j)= P_j$ and $Pr( L_v^n(m) = k)=1/V$ for any $v,k \in \{1,..,V\}$, $%
n \in \mathbb{N}$, $j \in \{1,...,N \}$, and $m \in \{1,...,M \}$.

First the family $\{L_{1}^{n},...,L_{V}^{n}\}$ generates a random dependence
tree, $K:T\rightarrow \{1,..,V\}$, in the following way. Let $K(\emptyset
)=1 $. If $K({i_{1},...,i_{n}})=j$, for some $j=1,...,V$, then we define $K({%
i_{1},...,i_{n},i_{n+1}})={L_{j}^{n}(i_{n+1})}$.

Given the dependence tree, let ${I}_{i}=N_{K(i)}^{n}$ if $|i|=n$.

The following theorem gives an alternative definition of $\rho_{V}$:

\begin{theorem}
$\rho _{V}([\tau ])=Pr( {I}_i = \tau(i),\ \forall \ |i| \leq |\tau| )$,
where $\{I_i\}_{i \in T} $ is defined as above, and $\tau $ is a finite
level code tree.
\end{theorem}

\begin{proof}
It is a straightforward exercise to check that $(\mathbf{F}_0 \circ \mathbf{F%
}_{1} \circ \cdots \circ \mathbf{F}_{k-1} )_1$ is a level-$k$ function tree
with nodes given by $\{I_i\}_{|i| \leq k-1}$, and limbs given by $\{ K(i)
\}_{|i| \leq k}$.

Thus 
\begin{equation*}
Pr(I_i=\tau(i), \ \forall \ i \ {\mathrm{with}}\ |i| \leq k-1 )
\end{equation*}
\begin{equation}  \label{107}
= Pr( (\mathbf{F}_0 \circ \mathbf{F}_{1} \circ \cdots \circ \mathbf{F}_{k-1}
)_1(node\ i) = \tau(i),\ \forall \ i \ {\mathrm{with}}\ |i| \leq k-1 ).
\end{equation}


For 
\begin{equation*}
\mathbf{a}=(a_{1},a_{2},...,a_{k})\in \mathcal{A}^{k}\text{.}
\end{equation*}%
let 
\begin{equation*}
\mathcal{P}^{\mathbf{a}}=\mathcal{P}^{a_{1}}\mathcal{P}^{a_{2}}...\mathcal{P}%
^{a_{k}}
\end{equation*}%
denote the probability of selection of 
\begin{equation*}
\eta ^{\mathbf{a}}=\eta ^{a_{1}}\circ \eta ^{a_{2}}\circ ...\circ \eta
^{a_{k}}=\eta ^{_{a_{1}\circ a_{2}\circ ...\circ a_{k}}}\text{.}
\end{equation*}

By the invariance of $\mu _{V}$%
\begin{equation*}
\mu _{V}=\sum\limits_{\mathbf{a}\in \mathcal{A}^{k}}\mathcal{P}^{\mathbf{a}%
}\eta ^{\mathbf{a}}(\mu _{V})\text{.}
\end{equation*}%
Now let $\tau $ be a level-$(k-1)$ code tree. Then 
\begin{equation*}
\rho _{V}([\tau ])=\mu _{V}(([\tau ],\Omega ,...,\Omega ))=\sum\limits_{%
\mathbf{a}\in \mathcal{A}^{k}}\mathcal{P}^{\mathbf{a}}\mu _{V}((\eta ^{%
\mathbf{a}})^{-1}(([\tau ],\Omega ,...,\Omega ))
\end{equation*}%
\begin{equation*}
=\sum\limits_{\{\mathbf{a}\in \mathcal{A}^{k}|(\eta ^{\mathbf{a}})_{1}(node\
i)=\tau (i),\ \forall \ i\}}\mathbf{\ }\text{ }\mathcal{P}^{\mathbf{a}}.
\end{equation*}

From this and (\ref{107}) it follows that $\rho _{V}([\tau ])=
Pr(I_i=\tau(i), \ \forall \ i \ {\mathrm{with}}\ |i| \leq k-1 )$.
\end{proof}

\subsection{\label{proofsec}Proof of Theorem \protect\ref{vinfinity}
Equation (\protect\ref{tobeproved}).}

\begin{proof}
We say that a dependence tree is \textit{free} up to level $k$, if at each
level $j$, for $1\leq j\leq k$, the $M^{j}$ limbs have distinct labels. If $%
V $ is much bigger than $M$ and $k$ then it is clear that the probability of
being free up to level $k$ is close to unity. More precisely, if $F$ is the
event that dependence tree of $(\mathbf{F}_{0}\circ \mathbf{F}_{1}\circ
\cdots \circ \mathbf{F}_{k-1})_{1}$ is free and $V\geq M^{k}$, then%
\begin{equation*}
\rho _{V}([\tau ])=\mu _{V}(([\tau ],\Omega ,...,\Omega ))=\sum\limits_{%
\mathbf{a}\in \mathcal{A}^{k}}\mathcal{P}^{\mathbf{a}}\mu _{V}((\eta ^{%
\mathbf{a}})^{-1}(([\tau ],\Omega ,...,\Omega ))
\end{equation*}%
\begin{equation*}
\Pr (F)=\prod\limits_{i=1}^{M-1}\left( 1-\frac{i}{V}\right)
\prod\limits_{i=1}^{M^{2}-1}\left( 1-\frac{i}{V}\right)
...\prod\limits_{i=1}^{M^{k}-1}\left( 1-\frac{i}{V}\right)
\end{equation*}%
\begin{equation*}
\geq 1-\frac{1}{V}\left(
\sum\limits_{i=1}^{M-1}i+\sum\limits_{i=1}^{M^{2}-1}i+...+\sum%
\limits_{i=1}^{M^{k}-1}i\right)
\end{equation*}%
\begin{equation*}
\geq 1-\frac{1}{2V}\left( M^{2}+M^{4}+...+M^{2k}\right)
\end{equation*}%
\begin{equation*}
\geq 1-\frac{M^{2(k+1)}}{2V(M^{2}-1)}\geq 1-\frac{2M^{2k}}{3V}\text{.}
\end{equation*}%
In the last steps we have assumed $M\geq 2$.

Let $S$ be the event that $(\eta ^{a_{1}\circ a_{2}\circ ...\circ
a_{k}})_{1}=\tau $. Then, \textit{using the independence of the random
variables labelling the nodes of a free function tree}, and using Equation (%
\ref{treeprobs}), we see that 
\begin{equation*}
\Pr (S|F)=\prod\limits_{\{i\in T|\text{ }1\leq \left\vert i\right\vert \leq
k\}}\mathcal{P}^{\tau (i)}=\rho ([\tau ])\text{.}
\end{equation*}%
Hence 
\begin{equation*}
\rho _{V}(\left[ \tau \right] )=\Pr (S)=\Pr (F)\Pr (S|F)+\Pr (F^{C})\Pr
(S|F^{C})
\end{equation*}%
\begin{equation*}
\leq \Pr (S|F)+\Pr (F^{C})\leq \rho \left( \left[ \tau \right] \right) +%
\frac{2M^{2k}}{3V}.
\end{equation*}%
Similarly,%
\begin{equation*}
\rho _{V}(\left[ \tau \right] )\geq \Pr (F)\Pr (S|F)
\end{equation*}%
\begin{equation*}
\geq \rho \left( \left[ \tau \right] \right) -\frac{2M^{2k}}{3V}.
\end{equation*}%
Hence%
\begin{equation}
\left\vert \rho _{V}(\left[ \tau \right] )-\rho (\left[ \tau \right]
)\right\vert \leqslant \frac{2M^{2k}}{3V}.  \label{winkarlin}
\end{equation}%
In order to compute the Monge Kantorovitch distance $d_{\mathbb{\mathcal{P}}%
(\Omega )}(\rho _{V},\rho )$, suppose $f:\Omega \rightarrow \mathbb{R}$ is
Lipshitz with $Lip$ $f\leq 1,$ i.e. $|f(\omega )-f(\varpi )|\leq $ $%
d_{\Omega }(\omega ,\varpi )$ $\forall $ $\omega ,\varpi \in \Omega $. Since 
$diam(\Omega )=1$ we subtract a constant from $f$ and so can assume $%
\left\vert f\right\vert \leq \frac{1}{2}$ without changing the value of $%
\tint fd\rho -\tint fd\rho _{V}$.

For each level-$k$ code tree $\tau \in T_{k}$ choose some $\omega _{\tau
}\in \left[ \tau \right] \subseteq \Omega $. It then follows that%
\begin{equation*}
\left\vert \int\limits_{\Omega }fd\rho -\int\limits_{\Omega }fd\rho
_{V}\right\vert =\left\vert \sum\limits_{\tau \in T_{k}}\left( \int\limits_{ 
\left[ \tau \right] }fd\rho -\int\limits_{\left[ \tau \right] }fd\rho
_{V}\right) \right\vert =
\end{equation*}%
\begin{equation*}
\left\vert \sum\limits_{\tau \in T_{k}}\int\limits_{\left[ \tau \right]
}\left( f-f(\omega _{\tau }\right) )d\rho -\sum\limits_{\tau \in
T_{k}}\int\limits_{\left[ \tau \right] }(\left( f-f(\omega _{\tau }\right)
)d\rho _{V}+\sum\limits_{\tau \in T_{k}}f(\omega _{\tau })(\rho (\left[ \tau %
\right] )-\rho _{V}(\left[ \tau \right] ))\right\vert
\end{equation*}%
\begin{equation*}
\leq \frac{1}{M^{k+1}}\sum\limits_{\tau \in T_{k}}\rho (\left[ \tau \right]
)+\frac{1}{M^{k+1}}\sum\limits_{\tau \in T_{k}}\rho _{V}(\left[ \tau \right]
)+\sum\limits_{\tau \in T_{k}}\frac{M^{2k}}{3V}
\end{equation*}%
\begin{equation*}
=\varphi (k):=\frac{2}{M^{k+1}}+\frac{M^{3k}}{3V}\text{,}
\end{equation*}%
since $diam$ $\left[ \tau \right] \leq \frac{1}{M^{k+1}}$ from Equation (\ref%
{size}), $\left\vert f(\omega _{\tau })\right\vert \leq \frac{1}{2},$ $Lip$ $%
f\leq 1,$ $\omega _{\tau }\in \left[ \tau \right] ,$ and using Equation \ref%
{winkarlin}. Choose $x$ so that $\frac{2V}{M}=M^{4x}$. This is the value of $%
x$ which minimizes $\left( \frac{2}{M^{x+1}}+\frac{M^{3x}}{3V}\right) $.
Choose $k$ so that $k\leq x\leq k+1.$ Then 
\begin{equation*}
\varphi (k)\leq 2\left( \frac{M}{2V}\right) ^{\frac{1}{4}}+\frac{1}{3V}%
\left( \frac{2V}{M}\right) ^{\frac{3}{4}}=2^{\frac{3}{4}}\left( \frac{M}{V}%
\right) ^{\frac{1}{4}}\left( 1+\frac{1}{3M}\right)
\end{equation*}%
\begin{equation*}
\leq \frac{7}{2^{\frac{1}{4}}3}\left( \frac{M}{V}\right) ^{\frac{1}{4}},%
\text{\ }(M\geq 2)\text{.}
\end{equation*}%
Hence Equation (\ref{tobeproved}) is true.
\end{proof}

\section{\label{five}Superfractals}

\subsection{\label{setsection}Contraction Mappings on $\mathbb{H}^{V}$ and
the Superfractal Set $\mathfrak{H}_{V,1}$.}

\begin{definition}
Let $V\in \mathbb{N}$, let $\mathcal{A}$ be the index set introduced in
Equation (\ref{index}), let $\mathcal{F}$ be given as in Equation (\ref%
{superIFS}), and let probabilities $\{\mathcal{P}^{a}|a\in \mathcal{A}\}$ be
given as in Equation (\ref{probs}). Define%
\begin{equation*}
f^{a}:\mathbb{H}^{V}\rightarrow \mathbb{H}^{V}
\end{equation*}%
by%
\begin{equation}
f^{a}(K)=(\bigcup\limits_{m=1}^{M}f_{m}^{n_{1}}(K_{v_{1,m}}),\bigcup%
\limits_{m=1}^{M}f_{m}^{n_{2}}(K_{v_{2,m}}),...\bigcup%
\limits_{m=1}^{M}f_{m}^{n_{V}}(K_{v_{V,m}}))  \label{faa}
\end{equation}%
$\forall K=(K_{1},K_{2},...,K_{V})\in \mathbb{H}^{V},\forall a\in \mathcal{A}
$. Let 
\begin{equation}
\mathcal{F}_{V}:=\{\mathbb{H}^{V};f^{a},\mathcal{P}^{a},a\in \mathcal{A}\}.
\label{ifsfq}
\end{equation}
\end{definition}

\begin{theorem}
\label{SuperIFS} $\mathcal{F}_{V}$ is an IFS with contractivity factor $l$.
\end{theorem}

\begin{proof}
We only need to prove that the mapping $f^{a}:\mathbb{H}^{V}\rightarrow 
\mathbb{H}^{V}$ is contractive with contractivity factor $l$, $\forall $ $%
a\in \mathcal{A}$. Note that, $\forall K=(K_{1},K_{2},...,K_{M}),$ $%
L=(L_{1},L_{2},...,L_{M})\in \mathbb{H}^{M}$, 
\begin{eqnarray*}
&&d_{\mathbb{H}}(\bigcup\limits_{m=1}^{M}f_{m}^{n}(K_{m}),\bigcup%
\limits_{m=1}^{M}f_{m}^{n}(L_{m})) \\
&\leq &\underset{m}{\max }\{d_{\mathbb{H}%
}(f_{m}^{n}(K_{m}),f_{m}^{n}(L_{m}))\} \\
&\leq &\underset{m}{\max \{}l\cdot d_{\mathbb{H}}(K_{m},L_{m})\} \\
&=&l\cdot d_{\mathbb{H}^{M}}(K,L)\text{.}
\end{eqnarray*}%
Hence, $\forall (K_{1},K_{2},...,K_{V}),(L_{1},L_{2},...,L_{V})\in \mathbb{H}%
^{V},$%
\begin{eqnarray*}
&&d_{\mathbb{H}%
^{V}}(f^{a}(K_{1},K_{2},...,K_{V}),f^{a}(L_{1},L_{2},...,L_{V})) \\
&=&\underset{v}{\max }\{d_{\mathbb{H}}(\bigcup%
\limits_{m=1}^{M}f_{m}^{n_{v}}(K_{v_{v,m}}),\bigcup%
\limits_{m=1}^{M}f_{m}^{n_{v}}(L_{v_{v,m}}))\} \\
&\leq &\underset{v}{\max }\{l\cdot d_{_{\mathbb{H}%
^{M}}}((K_{v_{v,1}},K_{v_{v,2}},...,K_{v_{v,M}}), \\
&&(L_{v_{v,1}},L_{v_{v,2}},...,L_{v_{v,M}}))\} \\
&\leq &l\cdot d_{_{\mathbb{H}%
^{V}}}((K_{1},K_{2},...,K_{V}),(L_{1},L_{2},...,L_{V}))\text{.}
\end{eqnarray*}
\end{proof}

The theory of IFS in Section 2.1 applies to the IFS $\mathcal{F}_{V}$. It
possesses a unique set attractor $\mathfrak{H}_{V}\in $\ $\mathbb{H(H}^{V})$%
, and a unique measure attractor $\mathfrak{P}_{V}$ $\in \mathbb{P(H}^{V}%
\mathbb{)}$. The random iteration algorithm corresponding to the IFS $%
\mathcal{F}_{V}$ may be used to approximate sequences of points ($V$- tuples
of compact sets) in $\mathfrak{H}_{V}$ distributed according to the
probability measure $\mathfrak{P}_{V}$.

However, the individual components of these vectors in $\mathfrak{H}_{V}$,
certain special subsets of $\mathbb{X}$, are the objects we are interested
in. Accordingly, for all $v\in \{1,2,...,V\}$, let us define $\mathfrak{H}%
_{V,v}$ $\subset \mathbb{H}$ to be the set of $v^{\text{th}}$ components of
points in $\mathfrak{H}_{V}$.

\begin{theorem}
\label{componentwise}For all $v\in \{1,2,...,V\}$ we have%
\begin{equation*}
\mathfrak{H}_{V,v}=\mathfrak{H}_{V,1}.
\end{equation*}%
When the probabilities in the superIFS $\mathcal{F}_{V}$ are given by (\ref%
{SuperProbs}), then starting from any initial $V$-tuple of non-empty compact
subsets of $\mathbb{X}$, the random distribution of the sets $K\in \mathbb{H}
$ that occur in the $v^{\text{th}}$ component of vectors produced by the
random iteration algorithm after $n$ initial steps converge weakly to the
marginal probability measure%
\begin{equation*}
\mathfrak{P}_{V,1}(B):=\mathfrak{P}_{V}(B,\mathbb{H},\mathbb{H},...,\mathbb{H%
})\forall B\in \mathbb{B}(\mathbb{H)}\text{,}
\end{equation*}%
independently of $v$, almost always, as $n\rightarrow \infty $.
\end{theorem}

\begin{proof}
The direct way to prove this theorem is to parallel the proof of Theorem \ref%
{vtrees}, using the maps $\{f^{a}:\mathbb{H}^{V}\rightarrow \mathbb{H}^{V}$ $%
|$ $a\in \mathcal{A}\}$ in place of the maps $\{\eta ^{a}:\Omega
^{V}\rightarrow \Omega ^{V}$ $|$ $a\in \mathcal{A}\}$.

However, an alternate proof follows with the aid of the map $\mathcal{F}%
:\Omega \mathbb{\rightarrow H(X)}$ introduced in Theorem \ref{codetree}$.$\
We have put this alternate proof at the end of the proof of Theorem \ref%
{code_space}.
\end{proof}

\begin{definition}
We call $\mathfrak{H}_{V,1}$ a \textit{superfractal set}. Points in $%
\mathfrak{H}_{V,1}$ are called \textit{V-variable fractal sets.}
\end{definition}

\begin{example}
\label{Two}See Figure \ref{redandgreen}. This example is similar to the one
in Section \ref{AnExample}. It shows some of the images produced in a
realization of random iteration of a superIFS with M=N=V=2. Projective
transformations are used in both IFSs, specifically%
\begin{equation*}
f_{1}^{1}(x,y)=(\frac{1.629x+0.135y-1.99}{-0.780x+0.864y-2.569},\frac{%
0.505x+1.935y-0.216}{0.780x-0.864y+2.569})\text{,}
\end{equation*}%
\begin{equation*}
f_{2}^{1}(x,y)=(\frac{1.616x-2.758y+3.678}{1.664x-0.944y+3.883},\frac{%
2.151x+0.567y+2.020}{1.664x-0.944y+3.883})\text{,}
\end{equation*}%
\begin{equation*}
f_{1}^{2}(x,y)=(\frac{1.667x+.098y-2.005}{-0.773x+0.790y-2.575},\frac{%
0.563x+2.064y-0.278}{0.773x-0.790y+2.575})\text{,}
\end{equation*}%
\begin{equation*}
f_{2}^{2}(x,y)=(\frac{1.470x-2.193y+3.035}{2.432x-0.581y+2.872},\frac{%
1.212x+0.686y+2.059}{2.432x-0.581y+2.872}).
\end{equation*}%
One of the goals of this example is to illustrate how closely similar images
can be produced, with \textquotedblleft random\textquotedblright\
variations, so the two IFSs are quite similar. Let us refer to images (or,
more precisely, the sets of points that they represent)\ such as the ones at
the bottom middle and at the bottom right, as \textquotedblleft
ti-trees\textquotedblright . Then each transformation maps approximately the
unit square $\square :=\{(x,y)\,|\,0\leq x\leq 1,\,0\leq y\leq 1\}$, in
which each ti-tree lies, into itself. Both $f_{2}^{1}(x,y)$ and $%
f_{2}^{2}(x,y)$ map ti-trees to lower right branches of ti-trees. Both $%
f_{1}^{1}(x,y)$ and $f_{1}^{2}(x,y)$ map ti-trees to a ti-tree minus the
lower right branch. The initial image for each component, or
\textquotedblleft screen\textquotedblright , is illustrated at the top left.
It corresponds to an array of pixels of dimensions $400\times 400$, some of
which are red, some green, and the rest white. Upon iteration, images of the
red pixels and green pixels are combined as in Example \ref{determ_example}.
The number of iterations increases from left to right, and from top to
bottom. The top middle image corresponds to the fifth iteration. Both the
images at the bottom middle and bottom left correspond to more than thirty
iterations, and are representive of typical images produced after more than
thirty iterations. (We carried out more than fifty iterations.) They
represent images selected from the superfractal $\mathfrak{H}_{2,1}$
according to the invariant measure $\mathfrak{P}_{2,1}$. Note that it is the
support of the red and green pixels that\ corresponds to an element of $%
\mathfrak{H}_{2,1}$. Note too the \textquotedblleft texture
effect\textquotedblright , similar to the one discussed in Example \ref%
{texture_example}.\FRAME{ftbpFU}{4.6942in}{3.1393in}{0pt}{\Qcb{Sequence of
images converging to 2-variable fractals, see Example \protect\ref{Two}.
Convergence to within the numerical resolution has occurred in the bottom
left two images. Note the subtle but real differences between the
silhouettes of \ these two sets. A variant of the \textquotedblleft texture
effect\textquotedblright\ can also be seen. The red points appear to dance
forever on the green ti-trees, while the ti-trees dance forever on the
superfractal.}}{\Qlb{redandgreen}}{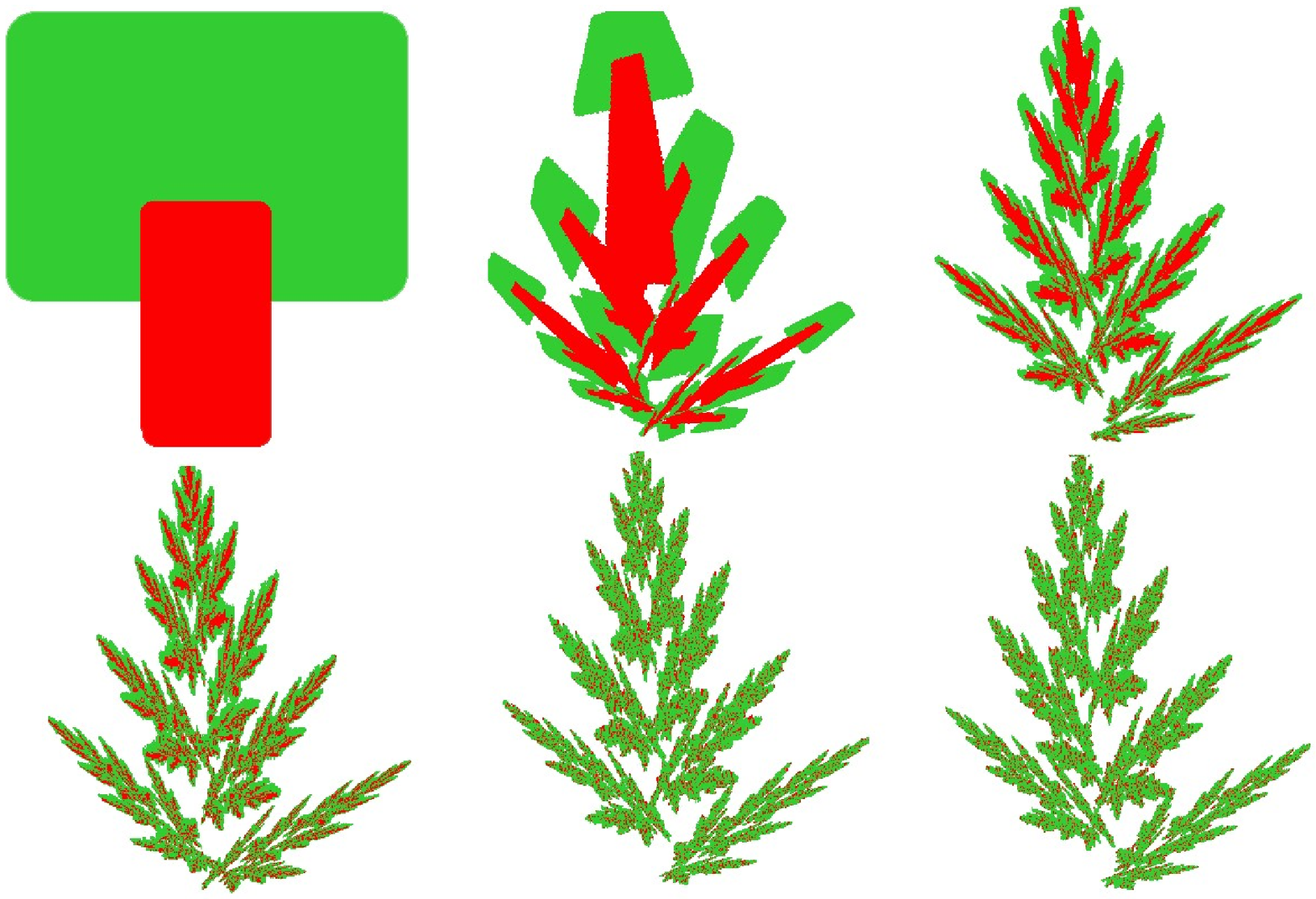}{\raisebox{-3.1393in}{\includegraphics[height=3.1393in]{supertitree.ps}}}
\end{example}

By Theorem \ref{Add}\ there is a continuous mapping $\mathcal{F}_{V}:\Sigma
_{V}\rightarrow \mathfrak{H}_{V}$ that assigns to each address in the code
space $\Sigma _{V}$ $=\mathcal{A}^{\infty }$a V-tuple of compact sets in $%
\mathfrak{H}_{V}$. But this mapping is not helpful for characterizing $%
\mathfrak{H}_{V}$ because $\mathcal{F}_{V}:\Sigma _{V}\rightarrow \mathfrak{H%
}_{V}$ is not in general one-to-one, for the same reason that $\Phi :\Sigma
_{V}\rightarrow \Omega _{V}$ is not one-to-one, as explained in Section \ref%
{backwards}.

The following result is closer to the point. It tells us in particular that
the set of $V$-trees is a useful code space for $V$-variable fractals,
because the action of the IFS $\Phi $ on the space of $V$-tuples of code
trees is conjugate to the action of the IFS $\mathcal{F}_{V}$ acting on
V-tuples of compact sets. (We are concerned here with the mappings that
provide the correspondences between $V$-groves, and $V$-trees, on the one
hand, and points and probability distributions on $\mathfrak{H}_{V}$ and $%
\mathfrak{H}_{V,1}$, on the other.)

\begin{theorem}
\label{code_space} Let $a\in \mathcal{A}$ and $\eta ^{a}:\Omega
^{V}\rightarrow \Omega ^{V}$ be defined as in Theorem \ref{jamie}. Let $%
f^{a}:\mathbb{H}^{V}\rightarrow \mathbb{H}^{V}$ be defined as in Theorem \ref%
{SuperIFS}. Let $\mathcal{F}:\Omega \mathbb{\rightarrow H(X)}$ be the
mapping introduced in Theorem \ref{codetree}. Define $\mathcal{F}:\Omega ^{V}%
\mathbb{\rightarrow (H(X))}^{V}$ by 
\begin{equation*}
\mathcal{F}(\omega _{1},\omega _{2},...,\omega _{V})=(\mathcal{F}(\omega
_{1}),\mathcal{F}(\omega _{2}),...,\mathcal{F}(\omega _{V}))\text{,}
\end{equation*}%
for all $(\omega _{1},\omega _{2},...,\omega _{V})\in \Omega ^{V}$.\ Then 
\begin{equation}
\mathcal{F}(\eta ^{a}(\omega ))=f^{a}(\mathcal{F}(\omega ))\text{ }\forall 
\text{ }a\in \mathcal{A}\text{, }\omega \in \Omega ^{V}\text{.}  \label{key}
\end{equation}%
Also 
\begin{equation}
\mathcal{F}(\Omega _{V})=\mathfrak{H}_{V}\text{ and }\mathcal{F}(\Omega
_{V,1})=\mathfrak{H}_{V,1}\text{,}  \label{custom1}
\end{equation}%
where $\Omega _{V,v}$ denotes the set of $v^{th}$ components of members of $%
\Omega _{V}$. Similarly, when the probabilities in the IFS $\mathcal{F}_{V}$
of Equation (\ref{ifsfq}), are given by Equation (\ref{SuperProbs}) we have%
\begin{equation}
\mathcal{F}(\mu _{V})=\mathfrak{P}_{V}\text{, and }\mathcal{F}(\rho _{V})=%
\mathfrak{P}_{V,1}\text{,}  \label{custom2}
\end{equation}%
where $\rho _{V}$ is the marginal probability distribution given by Equation
(\ref{marginal}).
\end{theorem}

\begin{proof}
We begin by establishing the key Equation (\ref{key}). Note that from
Theorem \ref{codetree}, for any $K\in \mathbb{H(X)}$,%
\begin{equation}
\mathcal{F}(\omega )=(\mathcal{F}(\omega _{1}),\mathcal{F}(\omega _{2}),...,%
\mathcal{F}(\omega _{V}))  \label{comp1}
\end{equation}%
\begin{equation*}
=(\underset{k\rightarrow \infty }{\lim }\{\mathcal{F}_{k}(\omega _{1})(K)\},%
\underset{k\rightarrow \infty }{\lim }\{\mathcal{F}_{k}(\omega
_{2})(K)\},...,\underset{k\rightarrow \infty }{\lim }\{\mathcal{F}%
_{k}(\omega _{V})(K)\})\text{.}
\end{equation*}%
The first component here exemplifies the others; and using Equation (\ref%
{frasets1}) it can be written%
\begin{equation}
\underset{k\rightarrow \infty }{\lim }\{\mathcal{F}_{k}(\omega _{1})(K)\}
\label{comp2}
\end{equation}%
\begin{equation*}
\mathbb{=}\lim_{k\rightarrow \infty }\mathbb{\{}\bigcup\limits_{\{i\in T|%
\text{ }|i|=k\}}f_{i_{1}}^{\omega _{1}(\emptyset )}\circ f_{i_{2}}^{\omega
_{1}(i_{1})}\circ ...\circ f_{i_{k}}^{\omega _{1}(i_{1}i_{2}...i_{k-1})}(K)\}%
\text{.}
\end{equation*}%
Since the convergence is uniform and all of the functions involved are
continuous, we can interchange the $\lim $ with function operation at will.
Look at 
\begin{equation*}
f^{a}(\mathcal{F}(\omega ))=f^{a}(\underset{k\rightarrow \infty }{\lim }\{%
\mathcal{F}_{k}(\omega _{1})(K)\},\underset{k\rightarrow \infty }{\lim }\{%
\mathcal{F}_{k}(\omega _{2})(K)\},...,\underset{k\rightarrow \infty }{\lim }%
\{\mathcal{F}_{k}(\omega _{V})(K)\})
\end{equation*}%
\begin{equation*}
=\underset{k\rightarrow \infty }{\lim }\{f^{a}(\mathcal{F}_{k}(\omega
_{1})(K),\mathcal{F}_{k}(\omega _{2})(K),...,\mathcal{F}_{k}(\omega
_{V})(K))\}\text{.}
\end{equation*}%
By the definition in Theorem \ref{SuperIFS}, Equation \ref{faa}, we have%
\begin{equation*}
f^{a}(\mathcal{F}_{k}(\omega _{1})(K),\mathcal{F}_{k}(\omega _{2})(K),...,%
\mathcal{F}_{k}(\omega _{V})(K))
\end{equation*}%
\begin{equation*}
=(\bigcup\limits_{m=1}^{M}f_{m}^{n_{1}}(\mathcal{F}_{k}(\omega
_{v_{1,m}})(K)),\bigcup\limits_{m=1}^{M}f_{m}^{n_{2}}(\mathcal{F}_{k}(\omega
_{v_{2,m}})(K)),...,\bigcup\limits_{m=1}^{M}f_{m}^{n_{V}}(\mathcal{F}%
_{k}(\omega _{v_{V,m}})(K)))\text{.}
\end{equation*}%
By equation (\ref{frasets1}) the first component here is%
\begin{equation*}
\bigcup\limits_{m=1}^{M}f_{m}^{n_{1}}(\mathcal{F}_{k}(\omega _{v_{1,m}})(K))
\end{equation*}%
\begin{equation*}
=\bigcup\limits_{m=1}^{M}f_{m}^{n_{1}}(\bigcup\limits_{\{i\in T|\text{ }%
|i|=k\}}f_{i_{1}}^{\omega _{v_{1,m}}(\emptyset )}\circ f_{i_{2}}^{\omega
_{v_{1,m}}(i_{1})}\circ ...\circ f_{i_{k}}^{\omega
_{v_{1,m}}(i_{1}i_{2}...i_{k-1})}(K))
\end{equation*}%
\begin{equation*}
=\mathcal{F}_{k+1}(\xi _{n_{1}}(\omega _{v_{1,1}},\omega
_{v_{1,2}},...,\omega _{v_{1,M}}))(K)\text{,}
\end{equation*}%
where we have used the definition in Equation (\ref{definition_jamie}). Hence%
\begin{equation*}
f^{a}(\mathcal{F}(\omega ))=\lim_{k\rightarrow \infty }\{(f^{a}(\mathcal{F}%
_{k}(\omega _{1})(K),\mathcal{F}_{k}(\omega _{2})(K),...,\mathcal{F}%
_{k}(\omega _{V})(K))\}=
\end{equation*}%
\begin{equation*}
\lim_{k\rightarrow \infty }\{(\mathcal{F}_{k+1}(\xi _{n_{1}}(\omega
_{v_{1,1}},\omega _{v_{1,2}},...,\omega _{v_{1,M}}))(K),\mathcal{F}%
_{k+1}(\xi _{n_{2}}(\omega _{v_{2,1}},\omega _{v_{2,2}},...,\omega
_{v_{2,M}}))(K),
\end{equation*}%
\begin{equation*}
...,\mathcal{F}_{k+1}(\xi _{n_{V}}(\omega _{v_{V,1}},\omega
_{v_{V,2}},...,\omega _{v_{V,M}}))(K))\}\text{.}
\end{equation*}%
Comparing with Equations (\ref{comp1}) and (\ref{comp2}), we find that the
right hand side here converges to $\mathcal{F}(\eta ^{a}(\omega ))$ as $%
k\rightarrow \infty $. So Equation (\ref{key}) is true.

Now consider the set $\mathcal{F}(\Omega _{V})$. We have%
\begin{equation*}
\mathcal{F}(\Omega _{V})=\mathcal{F}(\bigcup\limits_{a\in \mathcal{A}}\eta
^{a}(\Omega _{V}))=\bigcup\limits_{a\in \mathcal{A}}\mathcal{F}(\eta
^{a}(\Omega _{V}))=\bigcup\limits_{a\in \mathcal{A}}f^{a}(\mathcal{F}(\Omega
_{V}))\text{.}
\end{equation*}%
It follows by uniqueness that $\mathcal{F}(\Omega _{V})$ must be \textit{the}
set attractor of the IFS $\mathcal{F}_{V}$. Hence $\mathcal{F}(\Omega _{V})=%
\mathfrak{H}_{V}$ which is the first statement in Equation (\ref{custom1}).
Now%
\begin{equation*}
\mathcal{F}(\Omega _{V,1})=\{\mathcal{F}(\omega _{1})|(\omega _{1},\omega
_{2},...,\omega _{V})\in \Omega _{V}\}
\end{equation*}%
\begin{equation*}
=\text{\textit{\ first component of }}\{(\mathcal{F}(\omega _{1}),\mathcal{F}%
(\omega _{2}),...,\mathcal{F}(\omega _{V}))\text{ }|\text{ }(\omega
_{1},\omega _{2},...,\omega _{V})\in \Omega _{V}\}
\end{equation*}%
\begin{equation*}
=\text{\textit{first component of }}\mathcal{F}(\Omega _{V})=\text{\textit{%
first component of }}\mathfrak{H}_{V}=\mathfrak{H}_{V,1}\text{,}
\end{equation*}%
which contains the second statement in Equation (\ref{custom1}).

In a similar manner we consider the push-forward under $\mathcal{F}:\Omega
^{V}\rightarrow \mathbb{H}^{V}$ of the invariant measure $\mu _{V}$ of the
IFS $\Phi =\{\Omega ^{V};\eta ^{a},\mathcal{P}^{a},a\in \mathcal{A}\}$. $%
\mathcal{F}(\mu _{V})$ is normalized, i.e. $\mathcal{F}(\mu _{V})\in \mathbb{%
P}(\mathbb{H}^{V})$, because $\mathcal{F}(\mu _{V})(\mathbb{H}^{V})=\mu _{V}(%
\mathcal{F}^{-1}(\mathbb{H}^{V}))=$ $\mu _{V}(\Omega ^{V})$. We now show
that $\mathcal{F}(\mu _{V})$ is the measure attractor the IFS $\mathcal{F}%
_{V}$. The measure attractor of the IFS $\Phi $ obeys%
\begin{equation*}
\mu _{V}=\sum\limits_{a\in \mathcal{A}}\mathcal{P}^{a}\eta ^{a}(\mu _{V})%
\text{.}
\end{equation*}%
Applying $\mathcal{F}$ to both sides (i.e. constructing the push-fowards) we
obtain%
\begin{equation*}
\mathcal{F}(\mu _{V})=\mathcal{F}(\sum\limits_{a\in \mathcal{A}}\mathcal{P}%
^{a}\eta ^{a}(\mu _{V}))=\sum\limits_{a\in \mathcal{A}}\mathcal{P}^{a}%
\mathcal{F}(\eta ^{a}(\mu _{V}))=\sum\limits_{a\in \mathcal{A}}\mathcal{P}%
^{a}f^{a}(\mathcal{F}(\mu _{V}))\text{,}
\end{equation*}%
where in the last step we have used the key Equation (\ref{key}). So $%
\mathcal{F}(\mu _{V})$ is the measure attractor of the IFS $\mathcal{F}_{V}$%
. Using uniqueness, we conclude $\mathcal{F}(\mu _{V})=\mathfrak{P}_{V}$
which is the first equation in Equation (\ref{custom2}). Finally, observe
that, for all $B\in \mathbb{B}(\mathbb{H)}$, 
\begin{equation*}
\mathfrak{P}_{V,1}(B)=\mathcal{F}(\mu _{V})(B,\mathbb{H},\mathbb{H},...,%
\mathbb{H})=\mu _{V}(\mathcal{F}^{-1}(B,\mathbb{H},\mathbb{H},...,\mathbb{H}%
))
\end{equation*}%
\begin{equation*}
=\mu _{V}((\mathcal{F}^{-1}(B),\mathcal{F}^{-1}(\mathbb{H)},\mathcal{F}^{-1}(%
\mathbb{H)},...,\mathcal{F}^{-1}(\mathbb{H)}))\text{ (using Equation (\ref%
{comp1}))}
\end{equation*}%
\begin{equation*}
=\mu _{V}((\mathcal{F}^{-1}(B),\Omega ,\Omega ,...,\Omega ))\text{ (since }%
\mathcal{F}:\Omega \rightarrow \mathbb{H}\text{)}
\end{equation*}%
\begin{equation*}
=\rho _{V}(\mathcal{F}^{-1}(B))\text{ (by definition (\ref{marginal}))}=%
\mathcal{F}(\rho _{V})(B)\text{.}
\end{equation*}%
This contains the second equation in Equation (\ref{custom2}).

In a similar way, we obtain the alternate proof of Theorem \ref%
{componentwise}. Simply lift Theorem \ref{vtrees} to the domain of the IFS $%
\mathcal{F}_{V}$ using $\mathcal{F}:\Omega ^{V}\rightarrow \mathbb{H}^{V}$.
\end{proof}

The code tree $\Phi (a_{1}a_{2}a_{3}...)$ is called a tree address of the
V-variable fractal $\mathcal{F}(a_{1}a_{2}a_{3}...)$.

The mapping $\mathcal{F}:\Omega _{V,1}\rightarrow \mathfrak{H}_{V,1}$
together with Theorem \ref{Groves} provides a characterization of V-variable
fractals as follows. At any \textquotedblleft
magnification\textquotedblright , any V-variable fractal set is made of $V$
\textquotedblleft forms\textquotedblright\ or \textquotedblleft
shapes\textquotedblright :

\begin{theorem}
\label{forms}Let $M\in \mathfrak{H}_{V,1}$ be any $V$-variable fractal set.
Let $\epsilon >0$ be given. Then $M$ is a finite union of continuous
transformations of at most $V$ distinct compact subsets of $\mathbb{X}$, and
the diameter of each of these transformed sets is at most $\epsilon $.
\end{theorem}

\begin{proof}
Choose $n$ so that $l^{n}<\epsilon $. Note that 
\begin{equation*}
\mathfrak{H}_{V}=\bigcup\limits_{a\in \mathcal{A}}f^{a}(\mathfrak{H}%
_{V})=\bigcup\limits_{(a_{1},a_{2},...,a_{n})\in \mathcal{A}}f^{a_{1}}\circ
f^{a_{2}}\circ ...\circ f^{a_{n}}(\mathfrak{H}_{V})\text{.}
\end{equation*}%
Hence, since $(M,M,...,M)\in \mathfrak{H}_{V}$ it follows that there exists $%
(K_{1},K_{2},...,K_{V})\in \mathfrak{H}_{V}$ such that 
\begin{equation*}
M\in \text{\textit{\ first component of }}\bigcup%
\limits_{(a_{1},a_{2},...,a_{n})\in \mathcal{A}}f^{a_{1}}\circ
f^{a_{2}}\circ ...\circ f^{a_{n}}(K_{1},K_{2},...,K_{V})\text{.}
\end{equation*}%
Each set in the union has diameter at most $\epsilon $.
\end{proof}

\begin{example}
\label{exfive}This example is similar to Example \ref{Two}, with M=N=V=2.
The goal is to illustrate Theorem \ref{forms}. Figure \ref{kipling}\ shows,
from left to right, from top to bottom, a sequence of six successive images,
illustrating successive 2-variable fractals, corresponding to a superIFS of
two IFSs. Each IFS consists of two projective transformations, each mapping
the unit square $\square $ into itself. The images were obtained by running
the random iteration algorithm, as described in Section \ref{AnExample}. The
initial image on each screen was a blue convex region contained in a $%
400\times 400$ array representing $\square $, and the images shown
correspond to one of the discretized screens after forty, forty-one,
forty-two, forty-three, forty-four, and forty-five iterations. The key
features of the transformations can be deduced from the images.\ (For
example, one of the transformations of one of the IFSs, interpreted as a
mapping from one screen to the next, maps the top middle image to the top
two objects in the top left image.) Each of these images, at several scales,
looks as though it is the union of projective transformations of at most two
distinct sets.\FRAME{ftbpFU}{4.6942in}{3.141in}{0pt}{\Qcb{A sequence of
2-variable fractal sets (roughly accurate to viewing resolution),
corresponding to M=N=2. Can you spot the projective transformations? See
Example \protect\ref{Two}. Each of these six images images exhibits
\textquotedblleft 2-variability\textquotedblright : at several scales, each
looks as though it is the union of projective transformations of at most two
distinct sets.}}{\Qlb{kipling}}{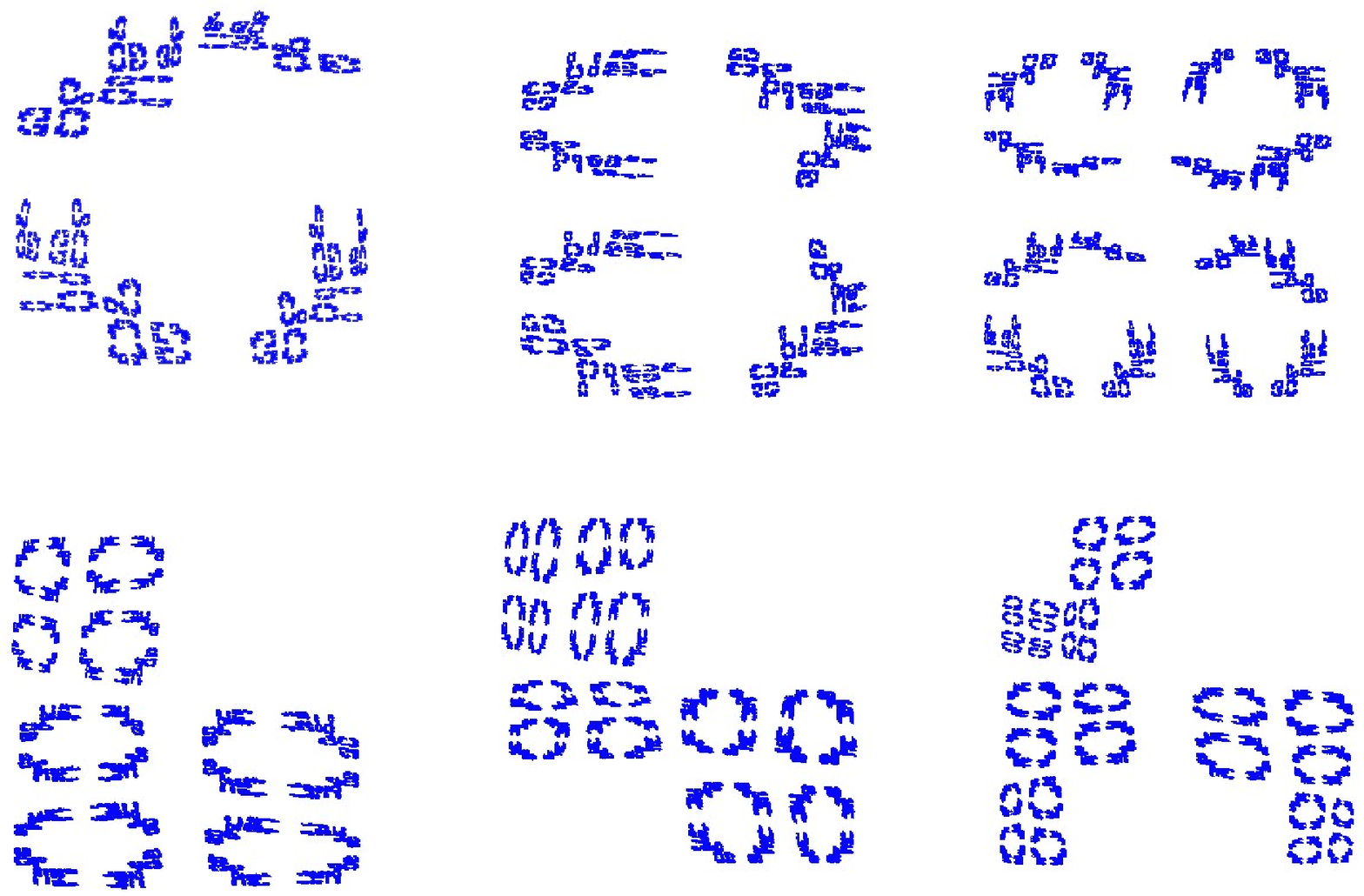}{\raisebox{-3.141in}{\includegraphics[height=3.141in]{superrects.ps}}}
\end{example}

\begin{theorem}
\label{conatinfinity}The set of V-variable fractal sets associated with the
superIFS $\mathcal{F}_{V}$ converges to the set of fractal sets associated
with the superIFS $\mathcal{F}$ introduced in Section \ref{randomsec}; that
is, in the metric of $\mathbb{H(H(X))}$, 
\begin{equation}
\underset{V\rightarrow \infty }{\lim }\mathfrak{H}_{V,1}=\mathfrak{H}\text{.}
\end{equation}%
Moreover, if the probabilities $\{\mathcal{P}^{a}|a\in \mathcal{A}\}$ obey
Equation (\ref{SuperProbs}), then in the metric of $\mathbb{P(H(X))}$ 
\begin{equation*}
\underset{V\rightarrow \infty }{\lim }\mathfrak{P}_{V,1}=\mathfrak{P}\text{,}
\end{equation*}%
where $\mathfrak{P}$ is the stationary measure on random fractal sets
associated with the superIFS $\mathcal{F}$.
\end{theorem}

\begin{proof}
We have, using the mapping $\mathcal{F}:\Omega \rightarrow \mathbb{H(X)}$,%
\begin{equation*}
\underset{V\rightarrow \infty }{\lim }\mathfrak{H}_{V,1}=\lim_{V\rightarrow
\infty }\mathcal{F}(\Omega _{V,1})\text{ (by Theorem \ref{code_space})}
\end{equation*}%
\begin{equation*}
=\mathcal{F}(\lim_{V\rightarrow \infty }\Omega _{V,1})\text{ (since }%
\mathcal{F}:\Omega \rightarrow \mathbb{H(X)}\text{ is continuous by Theorem %
\ref{codetree})}
\end{equation*}%
\begin{equation*}
=\mathcal{F}(\Omega )\text{ (by Theorem \ref{vinfinity})}
\end{equation*}%
\begin{equation*}
=\mathfrak{H}\text{ (by Equation (\ref{equationH})).}
\end{equation*}%
Similarly, using the mapping $\mathcal{F}:\Omega \rightarrow \mathbb{P(X)}$,
have%
\begin{equation*}
\underset{V\rightarrow \infty }{\lim }\mathfrak{P}_{V,1}=\lim_{V\rightarrow
\infty }\mathcal{F}(\rho _{V})\text{ (by Theorem \ref{code_space})}
\end{equation*}%
\begin{equation*}
=\mathcal{F}(\lim_{V\rightarrow \infty }\rho _{V})\text{ (since }\mathcal{F}%
:\Omega \rightarrow \mathbb{P(X)}\text{ is continuous by Theorem \ref%
{codetree})}
\end{equation*}%
\begin{equation*}
=\mathcal{F}(\rho )\text{ (by Theorem \ref{vinfinity})}
\end{equation*}%
\begin{equation*}
=\mathfrak{P}\text{ (by Equation (\ref{EquationI})).}
\end{equation*}
\end{proof}

\subsection{Contraction Mappings on $\mathbb{P}^{V}$ and the Superfractal
Measures $\mathfrak{\tilde{H}}_{V,1}$.}

Recall that $\mathbb{P}=\mathbb{P(X)}$. Let $\mathbb{P}^{V}=(\mathbb{P(X))}%
^{V}$ \ In this section we follow the same lines as in Section \ref%
{setsection}, constructing an IFS using the individual IFSs of the superIFS $%
\mathcal{F}$, except that here the underlying space is $\mathbb{P}^{V}$
instead of $\mathbb{H}^{V}.$

\begin{definition}
Let $V\in \mathbb{N}$, let $\mathcal{A}$ be the index set introduced in
Equation (\ref{index}), let $\mathcal{F}$ be given as in Equation (\ref%
{superIFS}), and let probabilities $\{\mathcal{P}^{a}|a\in \mathcal{A}\}$ be
given as in Equation (\ref{probs}). Define%
\begin{equation*}
f^{a}:\mathbb{P}^{V}\rightarrow \mathbb{P}^{V}
\end{equation*}%
by%
\begin{equation}
f^{a}(\mu )=(\sum\limits_{m=1}^{M}p_{m}^{n_{1}}f_{m}^{n_{1}}(\mu
_{v_{1,m}}),\sum\limits_{m=1}^{M}p_{m}^{n_{2}}f_{m}^{n_{2}}(\mu
_{v_{2,m}}),...,\sum\limits_{m=1}^{M}p_{m}^{n_{V}}f_{m}^{n_{V}}(\mu
_{v_{V,m}}))  \label{measop}
\end{equation}%
$\mu =(\mu _{1},\mu _{2},...,\mu _{V})\in \mathbb{P}^{V}.$ Let%
\begin{equation}
\mathcal{\tilde{F}}_{V}:=\{\mathbb{P}^{V};f^{a},\mathcal{P}^{a},a\in 
\mathcal{A\}}\text{.}  \label{ifsfvm}
\end{equation}
\end{definition}

\begin{theorem}
\label{SuperIFS_meas}$\mathcal{\tilde{F}}_{V}$ is an IFS with contractivity
factor $l$.
\end{theorem}

\begin{proof}
We only need to prove that the mapping $f^{a}:\mathbb{P}^{V}\rightarrow 
\mathbb{P}^{V}$ is contractive with contractivity factor $l$, $\forall $ $%
a\in \mathcal{A}$. Note that, $\forall \mu =(\mu _{1},\mu _{2},...,\mu
_{M}), $ $\varphi =(\varphi _{1},\varphi _{2},...,\varphi _{M})\in \mathbb{P}%
^{M}$,%
\begin{equation*}
d_{\mathbb{P}}(\sum\limits_{m=1}^{M}p_{m}^{n}f_{m}^{n}(\mu
_{m}),\sum\limits_{m=1}^{M}p_{m}^{n}f_{m}^{n}(\varphi _{m}))
\end{equation*}%
\begin{equation*}
\leq \sum\limits_{m=1}^{M}d_{\mathbb{P}}(p_{m}^{n}f_{m}^{n}(\mu
_{m}),p_{m}^{n}f_{m}^{n}(\varphi _{m}))
\end{equation*}%
\begin{equation*}
\leq l\cdot \max_{m}\{d_{\mathbb{P}}(\mu _{m},\varphi _{m})\}
\end{equation*}%
\begin{equation}
=l\cdot d_{\mathbb{P}^{M}}(\mu ,\varphi )\text{.}  \notag
\end{equation}%
Hence, $\forall (\mu _{1},\mu _{2},...,\mu _{V}),(\varphi _{1},\varphi
_{2},...,\varphi _{V})\in \mathbb{P}^{V},$%
\begin{eqnarray*}
&&d_{\mathbb{P}^{V}}(f^{a}(\mu _{1},\mu _{2},...,\mu _{V}),f^{a}(\varphi
_{1},\varphi _{2},...,\varphi _{V})) \\
&=&\underset{v}{\max }\{d_{\mathbb{P}}(\sum%
\limits_{m=1}^{M}p_{m}^{n_{v}}f_{m}^{n_{v}}(\mu
_{v_{v,m}}),\sum\limits_{m=1}^{M}p_{m}^{n_{v}}f_{m}^{n_{v}}(\varphi
_{v_{v,m}}))\} \\
&\leq &\underset{v}{\max }\{l\cdot d_{_{\mathbb{P}^{M}}}((\mu _{v_{v,1}},\mu
_{v_{v,2}},...,\mu _{v_{v,M}}), \\
&&(\varphi _{v_{v,1}},\varphi _{v_{v,2}},...,\varphi _{v_{v,M}}))\} \\
&\leq &l\cdot d_{_{\mathbb{P}^{V}}}((\mu _{1},\mu _{2},...,\mu
_{M}),(\varphi _{1},\varphi _{2},...,\varphi _{M}))\text{.}
\end{eqnarray*}
\end{proof}

The set attractor of the IFS\ $\mathcal{\tilde{F}}_{V}$ is $\mathfrak{\tilde{%
H}}_{V}\in \mathbb{H(P}^{V}),$ a subset of $\mathbb{P}^{V},$ a set of $V$%
-tuples of probability measures on $\mathbb{X}$. As we will see, each of
these measures is supported on a $V$-variable fractal set belonging to the
superfactal $\mathfrak{\tilde{H}}_{V,1}$. The measure attractor of the IFS $%
\mathcal{\tilde{F}}_{V}$ is a probability measure $\mathfrak{\tilde{P}}_{V}$ 
$\in \mathbb{P(P}^{V}\mathbb{)}$, namely a probability measure on a set of $%
V $-tuples of normalized measures, each one a random fractal measure. The
random iteration algorithm corresponding to the IFS $\mathcal{\tilde{F}}_{V}$
may be used to approximate sequences of points in $\mathfrak{\tilde{H}}_{V}$%
, namely vectors of measures on $\mathbb{X}$, distributed according to the
probability measure $\mathfrak{\tilde{P}}_{V}.$

As in Section \ref{setsection}, we define $\mathfrak{\tilde{H}}_{V,v}$ to be
the set of $v^{th}$ components of sets in $\mathfrak{\tilde{H}}_{V}$, for $%
v\in \{1,2,...,V\}.$

\begin{theorem}
\label{newrand}For all $v\in \{1,2,...,V\}$ we have%
\begin{equation*}
\mathfrak{\tilde{H}}_{V,v}=\mathfrak{\tilde{H}}_{V,1}.
\end{equation*}%
When the probabilities in the IFS $\mathcal{\tilde{F}}_{V}$ are given by
Equation (\ref{SuperProbs}), then starting at any initial $V$-tuple of
probability measures on $\mathbb{X}$, the probability measures $\mu \in 
\mathbb{P(X)}$ that occur in the $v^{\text{th}}$ component of points
produced by the random iteration algorithm after $n$ steps converge weakly
to the marginal probability measure%
\begin{equation*}
\mathfrak{\tilde{P}}_{V,1}(B):=\mathfrak{\tilde{P}}_{V}(B,\mathbb{P},\mathbb{%
P},...,\mathbb{P})\forall B\in \mathbb{B}(\mathbb{P)}\text{,}
\end{equation*}%
independently of $v$, almost always, as $n\rightarrow \infty $.
\end{theorem}

\begin{proof}
The direct way to prove this theorem is to parallel the proof of Theorem \ref%
{vtrees}, using the maps $\{f^{a}:\mathbb{P}^{V}\rightarrow \mathbb{P}^{V}$ $%
|$ $a\in \mathcal{A}\}$ in place of the maps $\{\eta ^{a}:\Omega
^{V}\rightarrow \Omega ^{V}$ $|$ $a\in \mathcal{A}\}$.

It is simpler however to lift Theorem \ref{vtrees} to the domain of the IFS $%
\widetilde{\mathcal{F}}_{V}$ using $\widetilde{\mathcal{F}}:\Omega
^{V}\rightarrow \mathbb{P}^{V}$ which is defined in Theorem \ref%
{codespace-meas}\ with the aid of the mapping $\widetilde{\mathcal{F}}%
:\Omega \mathbb{\rightarrow P=P(X)}$ introduced in Theorem \ref{codetree}.
We omit the details as they are straightforward.
\end{proof}

We call $\mathfrak{\tilde{H}}_{V,1}$ a superfractal set of measures (of
variability V). Points in $\mathfrak{\tilde{H}}_{V,1}$ are called V-variable
fractal measures.

\begin{example}
\label{one}See Figure \ref{greenandblack}. This example corresponds to the
same superIFS as in Example \ref{Two}. The probabilities of the functions in
the IFSs are $p_{1}^{1}=p_{1}^{2}=0.74$, and $p_{2}^{2}=p_{2}^{2}=0.26$. The
IFSs are assigned probabilities $P_{1}=P_{2}=0.5$. \FRAME{ftbpFU}{4.6942in}{%
1.5835in}{0pt}{\Qcb{Three successive 2-variable fractal measures computed
using the random iteration algorithm in Theorem \protect\ref{newrand}
applied to the superIFS in Example \protect\ref{one}. The pixels in the
support of each measure are coloured either black or a shade of green, using
a similar technique to the one used in Example \protect\ref{greeny}. The
intensity of the green of a pixel is a monotonic increasing function of the
measure of the pixel.}}{\Qlb{greenandblack}}{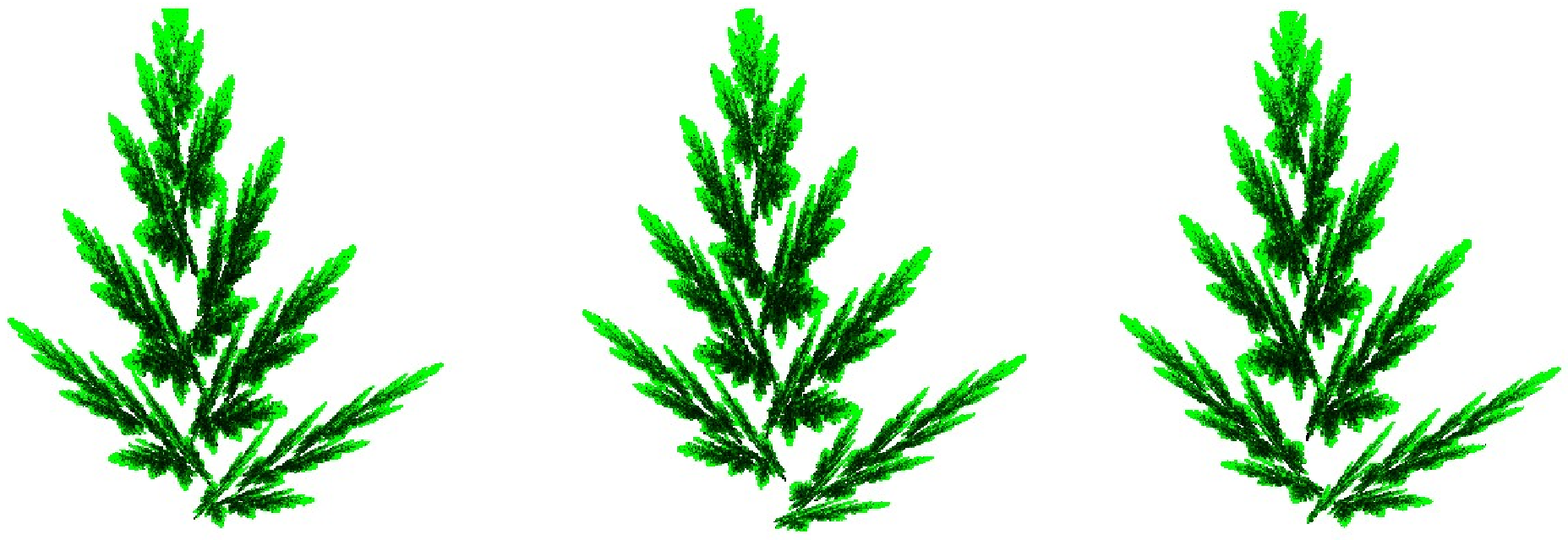}{\raisebox{-1.5835in}{\includegraphics[height=1.5835in]{measuretrees.ps}}}
\end{example}

The following Theorem tells us in particular that the set of $V$-trees is a
useful code space for $V$-variable fractal measures, because the action of
the IFS $\Phi $ on the space of $V$-tuples of code trees is conjugate to the
action of the IFS $\widetilde{\mathcal{F}}_{V}$ acting on V-tuples of
normalized measures.

\begin{theorem}
\label{codespace-meas} Let $a\in \mathcal{A}$ and $\eta ^{a}=\Omega
^{V}\rightarrow \Omega ^{V}$ be defined as in Theorem \ref{jamie}. Let $%
f^{a}:\mathbb{P}^{V}\rightarrow \mathbb{P}^{V}$ be defined as in Theorem \ref%
{SuperIFS_meas}. Let $\widetilde{\mathcal{F}}:\Omega \mathbb{\rightarrow
P=P(X)}$ be the mapping introduced in Theorem \ref{codetree}. Define $%
\widetilde{\mathcal{F}}:\Omega ^{V}\mathbb{\rightarrow P}^{V}=\mathbb{(P(X))}%
^{V}$ by 
\begin{equation*}
\widetilde{\mathcal{F}}(\omega _{1},\omega _{2},...,\omega _{V})=(\widetilde{%
\mathcal{F}}(\omega _{1}),\widetilde{\mathcal{F}}(\omega _{2}),...,%
\widetilde{\mathcal{F}}(\omega _{V}))\text{,}
\end{equation*}%
for all $(\omega _{1},\omega _{2},...,\omega _{V})\in \Omega ^{V}$.\ Then 
\begin{equation*}
\widetilde{\mathcal{F}}(\eta ^{a}(\omega ))=f^{a}(\widetilde{\mathcal{F}}%
(\omega ))\text{ }\forall \text{ }a\in \mathcal{A}\text{, }\omega \in \Omega
^{V}\text{.}
\end{equation*}%
Also 
\begin{equation*}
\widetilde{\mathcal{F}}(\Omega _{V})=\widetilde{\mathfrak{H}}_{V}\text{ and }%
\widetilde{\mathcal{F}}(\Omega _{V,1})=\widetilde{\mathfrak{H}}_{V,1}\text{,}
\end{equation*}%
where $\Omega _{V,v}$ denotes the set of $v^{th}$ components of members of $%
\Omega _{V}$. Similarly, when the probabilities in the IFS $\widetilde{%
\mathcal{F}}_{V}$ of Equation (\ref{ifsfvm}), are given by Equation (\ref%
{SuperProbs}) we have%
\begin{equation*}
\widetilde{\mathcal{F}}(\mu _{V})=\widetilde{\mathfrak{P}}_{V}\text{, and }%
\widetilde{\mathcal{F}}(\rho _{V})=\widetilde{\mathfrak{P}}_{V,1}\text{,}
\end{equation*}%
where $\rho _{V}$ is the marginal probability distribution given by Equation
(\ref{marginal}).
\end{theorem}

\begin{proof}
The proof is entirely parallel to the proof of Theorem \ref{code_space},
using $\widetilde{\mathcal{F}}$ in place of $\mathcal{F}$ and is omitted.
\end{proof}

\begin{definition}
The code tree $\Phi (a_{1}a_{2}a_{3}...)$ is called a tree address of the
V-variable fractal measure $\mathcal{\tilde{F}}_{V}(a_{1}a_{2}a_{3}...)$.
\end{definition}

The mapping $\mathcal{\tilde{F}}:\Omega _{V,1}\rightarrow \widetilde{%
\mathfrak{H}}_{V,1}$ together with Theorem \ref{Groves} allows us to
characterize V-variable fractals as follows:

\begin{theorem}
\label{forms_meas}Let $\widetilde{\mu }\in \widetilde{\mathfrak{H}}_{V,1}$
be any $V$-variable fractal measure. Let $\epsilon >0$ be given. Then $%
\widetilde{\mu }$ is a finite weighted superposition of continuous
transformations of at most $V$ distinct normalized measures supported on
compact subsets of $\mathbb{X}$, and the diameter of the support of each of
these transformed measures is at most $\epsilon $.
\end{theorem}

\begin{proof}
Choose $n$ so that $l^{n}<\epsilon $. Note that 
\begin{equation*}
\widetilde{\mathfrak{H}}_{V}=\bigcup\limits_{a\in \mathcal{A}}f^{a}(%
\widetilde{\mathfrak{H}}_{V})=\bigcup\limits_{(a_{1},a_{2},...,a_{n})\in 
\mathcal{A}}f^{a_{1}}\circ f^{a_{2}}\circ ...\circ f^{a_{n}}(\widetilde{%
\mathfrak{H}}_{V})\text{.}
\end{equation*}%
Hence, since $(\widetilde{\mu },\widetilde{\mu },...,\widetilde{\mu })\in 
\mathfrak{H}_{V}$ it follows that there exists $(\varpi _{1},\varpi
_{2},...,\varpi _{V})\in \mathfrak{H}_{V}$ such that 
\begin{equation*}
\widetilde{\mu }\in \text{\textit{\ first component of }}\bigcup%
\limits_{(a_{1},a_{2},...,a_{n})\in \mathcal{A}}f^{a_{1}}\circ
f^{a_{2}}\circ ...\circ f^{a_{n}}(\varpi _{1},\varpi _{2},...,\varpi _{V})%
\text{.}
\end{equation*}%
\ Inspection of Equation (\ref{measop}) shows that each of the measures in
the set of measures on the right-hand-side here is as stated in the theorem.
\end{proof}

\begin{example}
\label{exten}See Figure \ref{tworecsmes}. This corresponds to the same
superIFS as used in Example \ref{exfive} but here the measure is rendered in
shades of blue to provide a pictorial illustration of Theorem \ref%
{forms_meas}. The three successive images were computed with the aid of the
random iteration algorithm in Theorem \ref{newrand}, a new rendered measure
theoretic image being produced at each iteration. At each discernable scale,
approximately, each picture appears to have the property that it a
superposition of a number of \textquotedblleft little
pictures\textquotedblright\ belonging to one of two equivalence classes.
Pictures belonging to an equivalence class in this case are related by a
projective transformation together with a scaling of brightness. \FRAME{%
ftbpFU}{4.6942in}{1.5835in}{0pt}{\Qcb{Three successive 2-variable fractal
measures, in shades of blue. Illustrates the \textquotedblleft
shapes\textquotedblright\ and \textquotedblleft forms\textquotedblright\
theorem. See Example \protect\ref{exten}.}}{\Qlb{tworecsmes}}{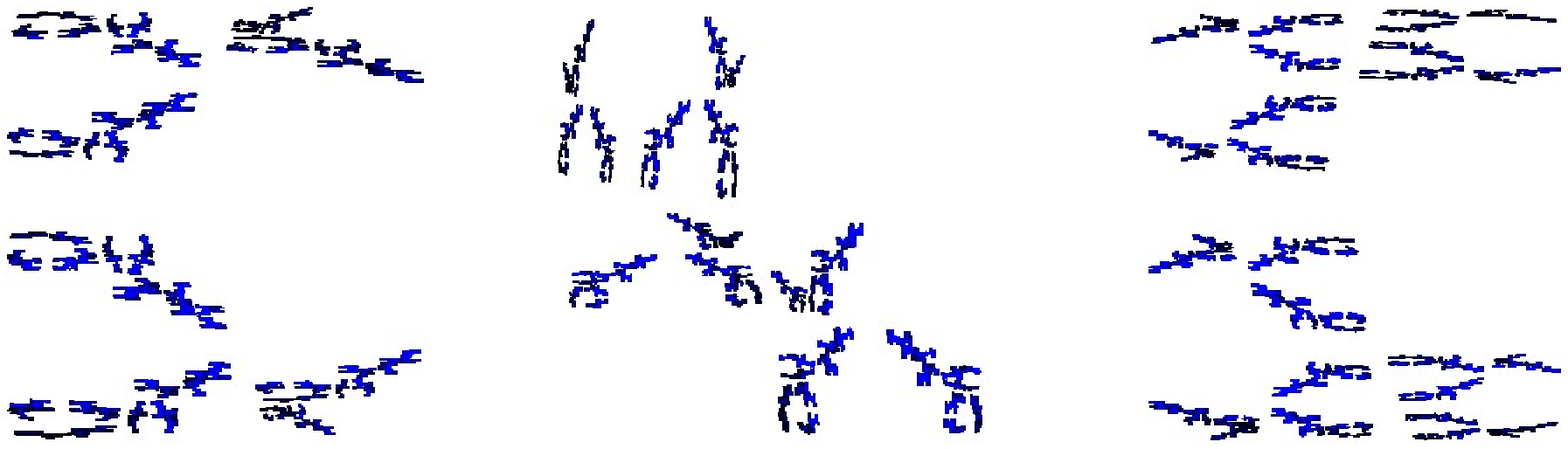%
}{\raisebox{-1.5835in}{\includegraphics[height=1.5835in]{tworecsmeas.ps}}}
\end{example}

\begin{theorem}
\label{theorem21}The set of V-variable fractal measures associated with the
superIFS $\mathcal{\tilde{F}}_{V}$ converges to the set of fractal measures
introduced in Section \ref{randomsec}; that is, in the metric of $\mathbb{%
H(P(X))}$%
\begin{equation*}
\underset{V\rightarrow \infty }{\lim }\mathfrak{\tilde{H}}_{V,1}=\mathfrak{%
\tilde{H}}.
\end{equation*}%
If the probabilities $\{\mathcal{P}^{a}|a\in \mathcal{A}\}$ obey Equation (%
\ref{SuperProbs}), then in the metric of $\mathbb{P(P(X))}$ 
\begin{equation*}
\underset{V\rightarrow \infty }{\lim }\widetilde{\mathfrak{P}}_{V,1}=%
\widetilde{\mathfrak{P}}\text{,}
\end{equation*}%
where $\widetilde{\mathfrak{P}}$ is the stationary measure on fractal sets
associated with the superIFS $\mathcal{F}$.
\end{theorem}

\begin{proof}
We have, using the mapping $\widetilde{\mathcal{F}}:\Omega \rightarrow 
\mathbb{P}=\mathbb{P(X)}$, 
\begin{equation*}
\underset{V\rightarrow \infty }{\lim }\widetilde{\mathfrak{H}}%
_{V,1}=\lim_{V\rightarrow \infty }\widetilde{\mathcal{F}}(\Omega _{V,1})%
\text{ (by Theorem \ref{codespace-meas})}
\end{equation*}%
\begin{equation*}
=\widetilde{\mathcal{F}}(\lim_{V\rightarrow \infty }\Omega _{V,1})\text{
(since }\widetilde{\mathcal{F}}:\Omega \rightarrow \mathbb{P}=\mathbb{P(X)}%
\text{ is continuous by Theorem \ref{codetree})}
\end{equation*}%
\begin{equation*}
=\widetilde{\mathcal{F}}(\Omega )\text{ (by Theorem \ref{vinfinity})}
\end{equation*}%
\begin{equation*}
=\widetilde{\mathfrak{H}}\text{ (by Equation (\ref{equationH})).}
\end{equation*}%
We have, using the mapping $\widetilde{\mathcal{F}}:\Omega \rightarrow 
\mathbb{P(P)}$, 
\begin{equation*}
\underset{V\rightarrow \infty }{\lim }\widetilde{\mathfrak{P}}%
_{V,1}=\lim_{V\rightarrow \infty }\widetilde{\mathfrak{P}}(\rho _{V})\text{
(by Theorem \ref{codespace-meas})}
\end{equation*}%
\begin{equation*}
=\widetilde{\mathcal{F}}(\lim_{V\rightarrow \infty }\rho _{V})\text{ (since }%
\widetilde{\mathcal{F}}:\Omega \rightarrow \mathbb{P(P)}\text{ is continuous
by Theorem \ref{codetree})}
\end{equation*}%
\begin{equation*}
=\widetilde{\mathcal{F}}(\rho )\text{ (by Theorem \ref{vinfinity})}
\end{equation*}%
\begin{equation*}
=\widetilde{\mathfrak{P}}\text{ (by Equation (\ref{EquationI})).}
\end{equation*}
\end{proof}

\subsection{Fractal Dimensions}

Here we quantify and compare the Hausdorff dimensions of fractals
corresponding to a (super) IFS of similitudes on $\mathbb{R}^{K}$ for some $%
K\in \mathbb{N}$ that obeys the open set condition in the following four
cases: deterministic fractals, standard random fractals, homogeneous\ random
fractals ($V=1$), and $V$-variable fractals ($V>1$). The functions of the
IFS $F^{n}$ are of the form $f_{m}^{n}(x)=s_{m}^{n}O_{m}^{n}x+t_{m}^{n}$
where $O_{m}^{n}$ is an orthonormal transformation, $s_{m}^{n}\in (0,1)$,
and $t_{m}^{n}\in \mathbb{R}^{K},$ for all $n\in \{1,2,...,N\}$ and $m\in
\{1,2,...M\}.$

\subsubsection{Deterministic Fractals}

In this case there is only one IFS, say $F^{1}$. By Theorem \ref%
{dimension_theorem} the Hausdorff dimension of the corresponding fractal set 
$A$ is $D$, the unique solution of 
\begin{equation*}
\sum_{m=1}^{M}(s_{m}^{1})^{D}=1\text{.}
\end{equation*}

\begin{example}
\label{dimexone}Suppose $K\geq 2$. Let the IFS $F^{1}$ consists of three
similitudes with $s_{1}^{1}=s_{2}^{1}=s_{3}^{1}=\frac{1}{2}$ and that the
fixed points are not collinear. Then the set attractor of $F^{1}$ is the
Sierpinski triangle with vertices at the three fixed points. Its fractal
dimension $D_{1}$ is given by $\ \,3\,\frac{1}{2^{D_{1}}}=1$ which implies $%
D_{1}=\frac{\ln 3}{\ln 2}=1.585$.

Let the IFS $F^{2}$ consist of three similitudes with $%
s_{1}^{2}=s_{2}^{2}=s_{3}^{2}=\frac{1}{3}$ and the same fixed points as $%
F^{1}$. Then the fractal dimension $D_{2}$ of the set attractor of $F^{2}$
is the given by $\ \,3\,\frac{1}{3^{D_{2}}}=1$ which implies $D_{2}=1$.
\end{example}

\subsubsection{Random Fractals}

By Theorem \ref{theorem_above} the Hausdorff dimension $D_{R}$ of $\mathfrak{%
P}$-almost all of the random fractals sets for the superIFS $\mathcal{F}$ is
given by 
\begin{equation*}
\sum\limits_{n=1}^{N}P_{n}\sum_{m=1}^{M}(s_{m}^{n})^{D_{R}}=1\text{.}
\end{equation*}

\begin{example}
\label{dimextwo}Let the superIFS be $\{\square
;F^{1},F^{2};P_{1}=P_{2}=0.5\} $ where the IFS's are defined in Example \ref%
{dimexone}. Then the fractal dimension $D_{R}$ of $\mathfrak{P}$-almost all
of the random fractals in the set is given by $\frac{1}{2}\,3\,\frac{1}{%
2^{D_{R}}}+\frac{1}{2}\,3\,\frac{1}{3^{D_{R}}}=1\Longrightarrow D_{R}=1.262.$
\end{example}

\subsubsection{Homogeneous Random Fractals ($V=1$)}

The case of homogeneous\ random fractals corresponds to $V=1.$Each run of
the experiment gives a different random Sierpinski triangle.

\begin{theorem}
\cite{Ham}. Let the superIFS $\mathcal{F}$ be as specified as in Theorem \ref%
{theorem_above}. Let $V=1$. Then for $\mathfrak{P}_{1,1}$ almost all $A\in 
\mathfrak{H}_{1,1}$%
\begin{equation*}
\dim _{H}A=D
\end{equation*}%
where $D$ is the unique solution of 
\begin{equation*}
\sum\limits_{n=1}^{N}P_{n}\ln \sum_{m=1}^{M}(s_{m}^{n})^{D}=1\text{.}
\end{equation*}
\end{theorem}

\begin{example}
For the case of the superIFS in Example \ref{dimextwo}, whose $1$-variable
fractal sets we refer to as homogeneous random Sierpinski triangles, the
Hausdorff dimension $D$ of almost all of them is given by $\frac{1}{2}\log
\left( 3\frac{1}{2^{D}}\right) +\frac{1}{2}\log \left( 3\frac{1}{3^{D}}%
\right) =0,$ $\Longrightarrow d=2\log 3/(\log 2+\log 3)=1.226.$
\end{example}

\subsubsection{$V$-Variable Fractals ($V\geq 1)$}

Let $(a_{1},a_{2},...)\in $ $\mathcal{A}^{\infty }$ denote an i.i.d.
sequence of indices, with probabilities $\{\mathcal{P}^{a}|a\in \mathcal{A}%
\} $ given in terms of the probabilities $\{P_{1},P_{2,}...,P_{V}\}$
according to Equation (\ref{SuperProbs}). Define, for $\alpha \in \lbrack
0,\infty )$ and $a\in \mathcal{A}$, the $V\times V$ \textit{flow matrix}%
\begin{equation*}
M_{v,w}^{a}(\alpha )=\sum\limits_{\{m|v_{v,m}=w\}}(s_{m}^{n_{v}\text{ }%
})^{\alpha }\text{,}
\end{equation*}%
and let us write%
\begin{equation*}
M_{v,w}^{k}=M_{v,w}^{k}(\alpha )=M_{v,w}^{a_{k}}(\alpha ).
\end{equation*}%
We think of $s_{m}^{n_{v}\text{ }}$ as being the \textquotedblleft
flow\textquotedblright\ through the $m^{th}$ channel from screen $v$ to
screen $w$.where $v_{v,m}=w$. The sequence of random matrices $M_{v,w}^{1}$, 
$M_{v,w}^{2},...$ is i.i.d., again with probabilities induced from $%
\{P_{1},P_{2,}...,P_{V}\}$. For any real square matrix $M$ we define the
norm $\Vert M\Vert $ to be the sum of the absolute values of its entries. By
the Furstenberg Kesten Theorem \cite{furstenberg},%
\begin{equation*}
\gamma (\alpha ):=lim_{k\rightarrow \infty }k^{-1}\log \Vert M^{1}(\alpha
)\circ \dots \circ M^{k}(\alpha )\Vert
\end{equation*}%
exists and has the same value with probability one. Provided that the
superIFS\ obeys the open set condition, we have shown in \cite{BHS2} that
the unique value of $D\in \lbrack 0,\infty )$ such that 
\begin{equation*}
\gamma (D)=0
\end{equation*}%
is the Hausdorff dimension of $\mathfrak{P}_{V,1}$ almost all $A\in 
\mathfrak{H}_{V,1}$.

Kingman remarks that in general the calculation of $\gamma $
\textquotedblleft has pride of place among the unsolved problems of
subadditive ergodic theory\textquotedblright , \cite{kingman}, p.897.
However it is possible to estimate numerically. Namely, generate random
copies of $M^{k}$ and iteratively compute $M^{1},M^{2},...,M^{k}$ and hence $%
k^{-1}\log \Vert M^{1}(\alpha )\circ \dots \circ M^{k}(\alpha )\Vert $ for $%
k=1,2,...$ The limit will give $\gamma (\alpha )$. (Even for large $V$ this
will be quick since the $M^{k}$ are sparse.) One could now use the bisection
method to estimate $D$.

\section{\label{six}Applications}

Fractal geometry plays some role in many application areas, including the
following. In biology: breast tissue patterns, structure and development of
plants, blood vessel patterns, and morphology of fern fronds. In chemistry:
pattern-forming alloy solidification, and diffusion processes. In physics:
transport in porous media, patterns formed during phase transitions in
statistical mechanics, dynamical systems, turbulence and wave propagation.
In geology: particle size distribution in soil, and landscape habitat
diversity. In computer science: digital image compression and watermarking,
compression of digital audio signals, image segmentation, and computer
graphics. In engineering: wavelets, stochastic processes, rough surfaces,
antennae and frequency selective surfaces, stochastic optimal control,
signal processing, and fragmentation of thin plates.

In many of these areas it is clearly desireable to use random fractals; for
example random fractals can be used in connection with diverse mathematical
modeling application areas including Brownian motion, oil-wells, critical
phenomena in statistical physics, for example associated with lattice gasses
and percolation, stock-market prices in finance, and in computer graphics
they can be used to represent diverse picture types including natural images
and textures. But random fractals are hard to compute, which may have held
up the development of some applications, while deterministic fractals, which
can be computed relatively easily, may not be rich enough to provide
convenient models for the applications to which one would want to apply them.

Thus we believe that $V$-variable fractals could find many applications;
they can be computed easily, with rapid access to many examples, contain a
controllable amount of \textquotedblleft randomness\textquotedblright , and
have many of the advantages of fractals in general: for similitudes, with an
open set condition, their fractal dimension may be computed, they are
resolution independent, and they are in general geometrically complex at all
levels of magnification, while being expressed with relatively small amounts
of information, coefficients of affine transformations and some
probabilities, for example.

\subsection{Space-filling curves}

Space-filling curves can be constructed with the aid of IFS theory, see for
example \cite{sagan}, Chapter 9. These curves have many applications,
including adaptive multigrid methods for numerical computation of solutions
of PDEs and hierarchical watermarking of digital images. \ Here we note that
interesting $V$-variable space-filling curves, and finite resolution
approximants to them, can be produced.

\begin{example}
\label{fillerex} Let $M=3,V=2,N=2$. The IFS $F^{1}=\{\square
;f_{1}^{1},f_{2}^{1},f_{3}^{1}\}$ consists of affine maps whose actions we
explain with the aid of the left-hand diagram in Figure \ref{spacefiller}. $%
\square $ is the unit square in the diagram, while $f_{1}^{1}(\square )$ is
the lower left square, $f_{2}^{1}(\square )$ is the upper left square, and $%
f_{3}^{1}(\square )$ is the rectangle on the right. The transformations are
chosen so that $f_{1}^{1}(\overline{OC})=\overline{OA}$, $f_{2}^{1}(%
\overline{OC})=\overline{AB},$ and $f_{3}^{1}(\overline{OC})=\overline{BC}$.
Specifically $f_{1}^{1}(x,y)=(\frac{1}{2}y,\frac{1}{2}x),f_{2}^{1}(x,y)=(-%
\frac{1}{2}y+\frac{1}{2},-\frac{1}{2}x+1),f_{3}^{1}(x,y)=(\frac{1}{2}x+\frac{%
1}{2},-y+1).$

The IFS$\ F^{2}=\{\square ;f_{1}^{2},f_{2}^{2},f_{3}^{2}\}$ is explained
with the aid of the right-hand diagram in Figure \ref{spacefiller}; $%
f_{1}^{2}(\square )$ is the lower left rectangle, $f_{2}^{2}(\square )$ is
the upper left rectangle, and $f_{3}^{2}(\square )$ is the rectangle on the
right; such that $f_{1}^{2}(\overline{OC})=\overline{OA^{\prime }}$, $%
f_{2}^{2}(\overline{OC})=\overline{A^{\prime }B^{\prime }},$ and $f_{3}^{2}(%
\overline{OC})=\overline{B^{\prime }C}.$ Specifically $f_{1}^{2}(x,y)=(\frac{%
2}{3}y,\frac{1}{2}x),f_{2}^{2}(x,y)=(-\frac{2}{3}y+\frac{2}{3},-\frac{1}{2}%
x+1),f_{3}^{2}(x,y)=(\frac{1}{3}x+\frac{2}{3},-y+1).$

Neither of the IFSs here is strictly contractive, but each is contractive
\textquotedblleft on the average\textquotedblright , for any assignment of
positive probabilities to the constituent functions. We assign probabilities 
$P_{1}=P_{2}=0.5$ to the individual IFSs.\ An initial image consisting of
the line segment $\overline{OC}$ is chosen on both screens, and \ the random
iteration algorithm is applied; typical images produced after five
iterations are illustrated in Figure \ref{filler20}; an image produced after
seven iterations is shown in Figure \ref{filler18}. Each of these images
consists of line segments that have been assigned colours according to the
address of the line segment, in such as way as to provide some consistency
from one image to the next. \FRAME{ftbpFU}{4.644in}{1.97in}{0pt}{\Qcb{%
Transformations used for space-filling curves, see Example \protect\ref%
{fillerex}}}{\Qlb{spacefiller}}{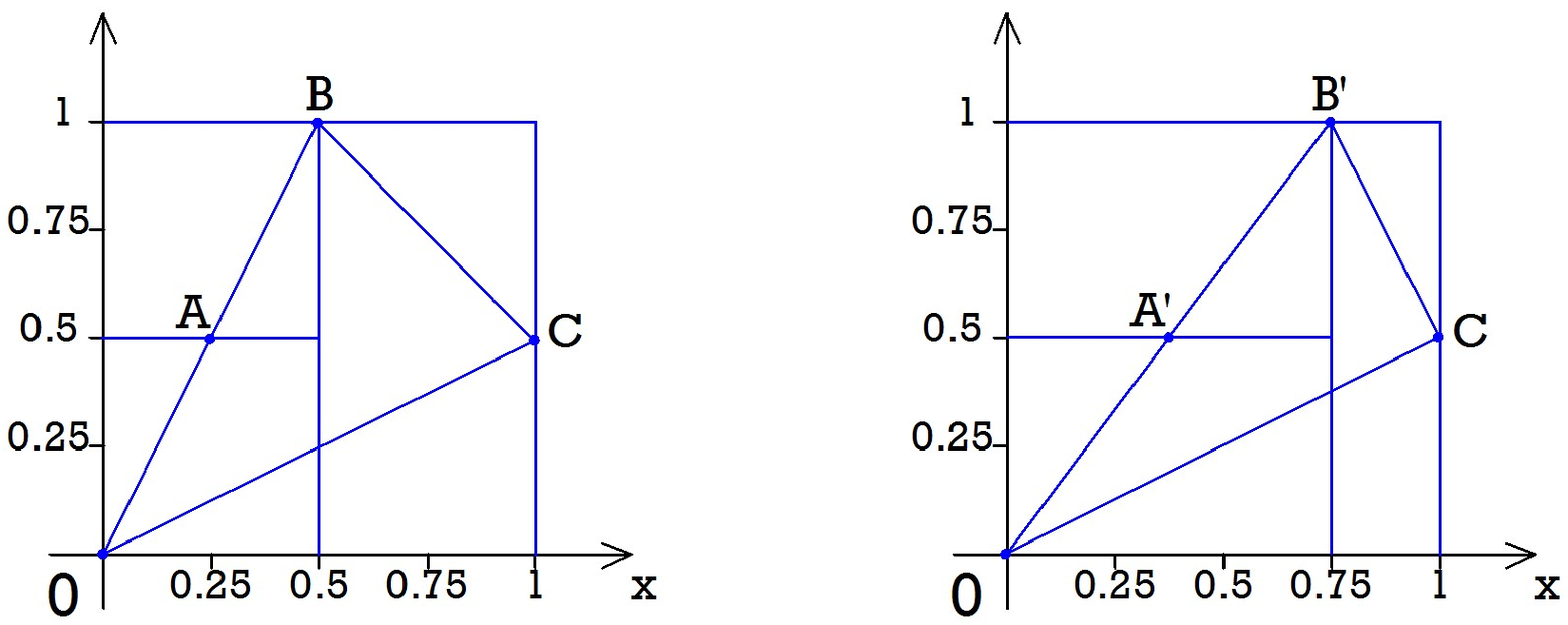}{\raisebox{-1.97in}{\includegraphics[height=1.97in]{trianglefig.ps}}}\FRAME{ftbpFU}{4.5282in}{1.6985in}{0pt}{\Qcb{Low-order
approximants to two 2-variable space filling curves, belonging to the same
superfractal, see Example \protect\ref{fillerex}}}{\Qlb{filler20}}{%
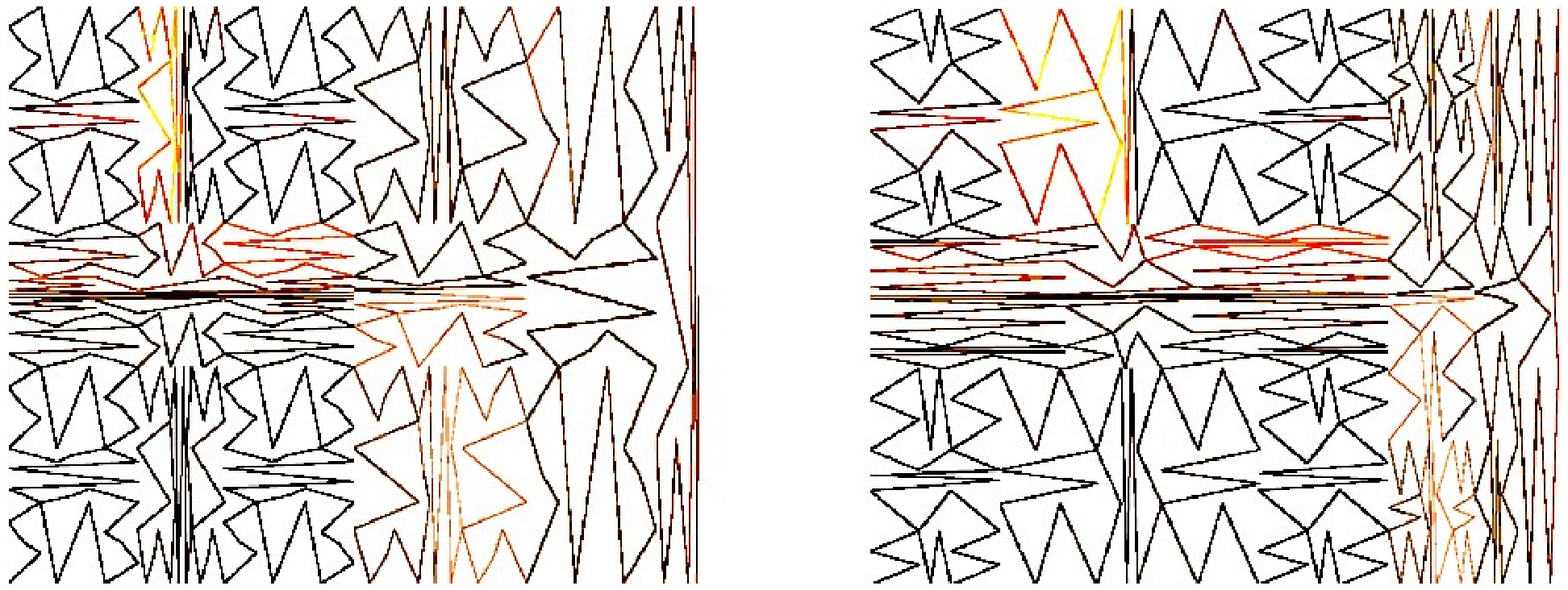}{\raisebox{-1.6985in}{\includegraphics[height=1.6985in]{spacefiller20.ps}}}\FRAME{ftbpFU}{4.5083in}{%
3.7671in}{0pt}{\Qcb{Finite-resolution approximation to a 2-variable
space-filling curve. See Example \protect\ref{fillerex}.}}{\Qlb{filler18}}{%
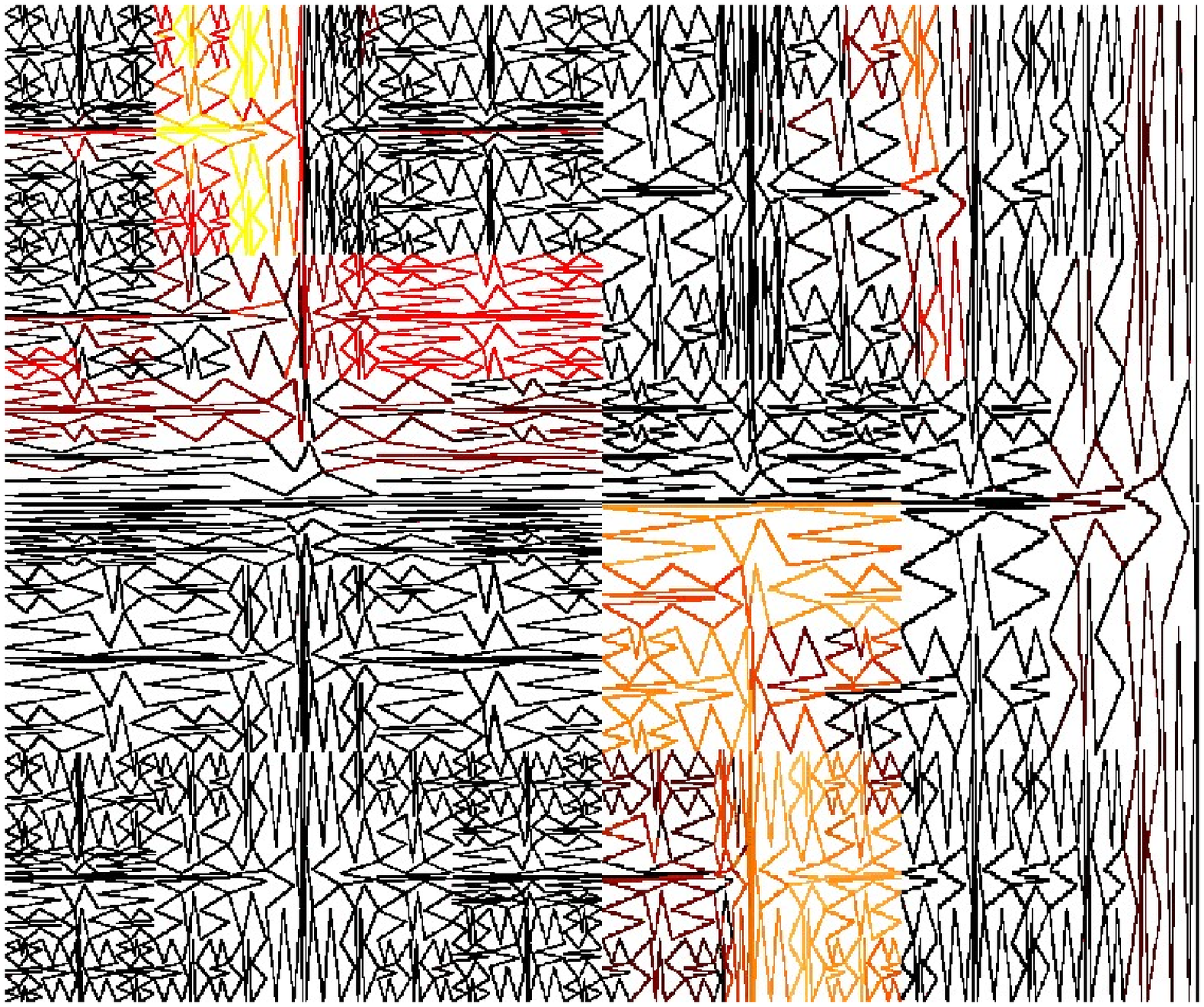}{\raisebox{-3.7671in}{\includegraphics[height=3.7671in]{spacefiller18.ps}}}
\end{example}

\subsection{Computer Graphics}

New techniques in computer graphics are playing an increasingly important
role in the digital content creation industry, as evidenced by the
succession of successes of computer generated films, from \textquotedblleft
Toy Story\textquotedblright\ to \textquotedblleft Finding
Nemo\textquotedblright . Part of the appeal of such films is the artistic
quality of the graphics. Here we point out that $V$-variable fractals are
able to provide new types of rendered graphics, significantly extending
standard IFS\ graphics \cite{Ba88}.

\begin{example}
\label{grafex}Here N=2, V=2, M=4. The two IFSs are given by 
\begin{equation*}
F^{n}=\text{\ }\{\square ;f_{1}^{n},f_{2}^{n},f_{3}^{n},f_{4}^{n};\}\text{ }%
n\in \{1,2\},
\end{equation*}%
where $\square \subset \mathbb{R}^{2},$ and each $f_{m}^{n}:\square
\rightarrow \square $ is a projective transformation. The colours were
obtained as follows. A computer graphics rendering of the set attractor of $%
F^{1}$ is shown in Figure \ref{Fig1}, and of $F^{2}$ in Figure \ref{Fig2}.

The colouring of each of these two figures was obtained with the aid of an
auxiliary IFS acting on the cube $C:=[0,255]^{3}\subset \mathbb{R}^{3}$
given by $\mathcal{G}:=\{C;g_{1}^{n},g_{2}^{n},g_{3}^{n},g_{4}^{n}\}$ where
each $g_{m}$ is a contractive (in the Euclidean metric) affine
transformation, represented by a $3\times 3$ matrix and a $3\times 1$
vector. For $n\in \{1,2\}$ discretized approximations, of the same
resolution,\ to the attractors of both IFSs $F^{n}$ and $\mathcal{G}$ were
calculated via the deterministic algorithm (Corollary \ref{detalgcor}); each
pixel on the attractor of the IFS $F^{n}$ was assigned the colour whose red,
green, and blue components, each an integer from $0$ to $255$, were the
three coordinates of the point on the attractor of $\mathcal{G}$ with the
same code space address. At those points in the attractor of $F^{n}$ with
multiple code space addresses, the lowest address was chosen.\FRAME{ftbpFU}{%
4.0283in}{3.3702in}{0pt}{\Qcb{The rendered set attractor of the IFS $F^{1}$%
in Example \protect\ref{grafex}.}}{\Qlb{Fig1}}{ferncolorx.ps}{\raisebox{-3.3702in}{\includegraphics[height=3.3702in]{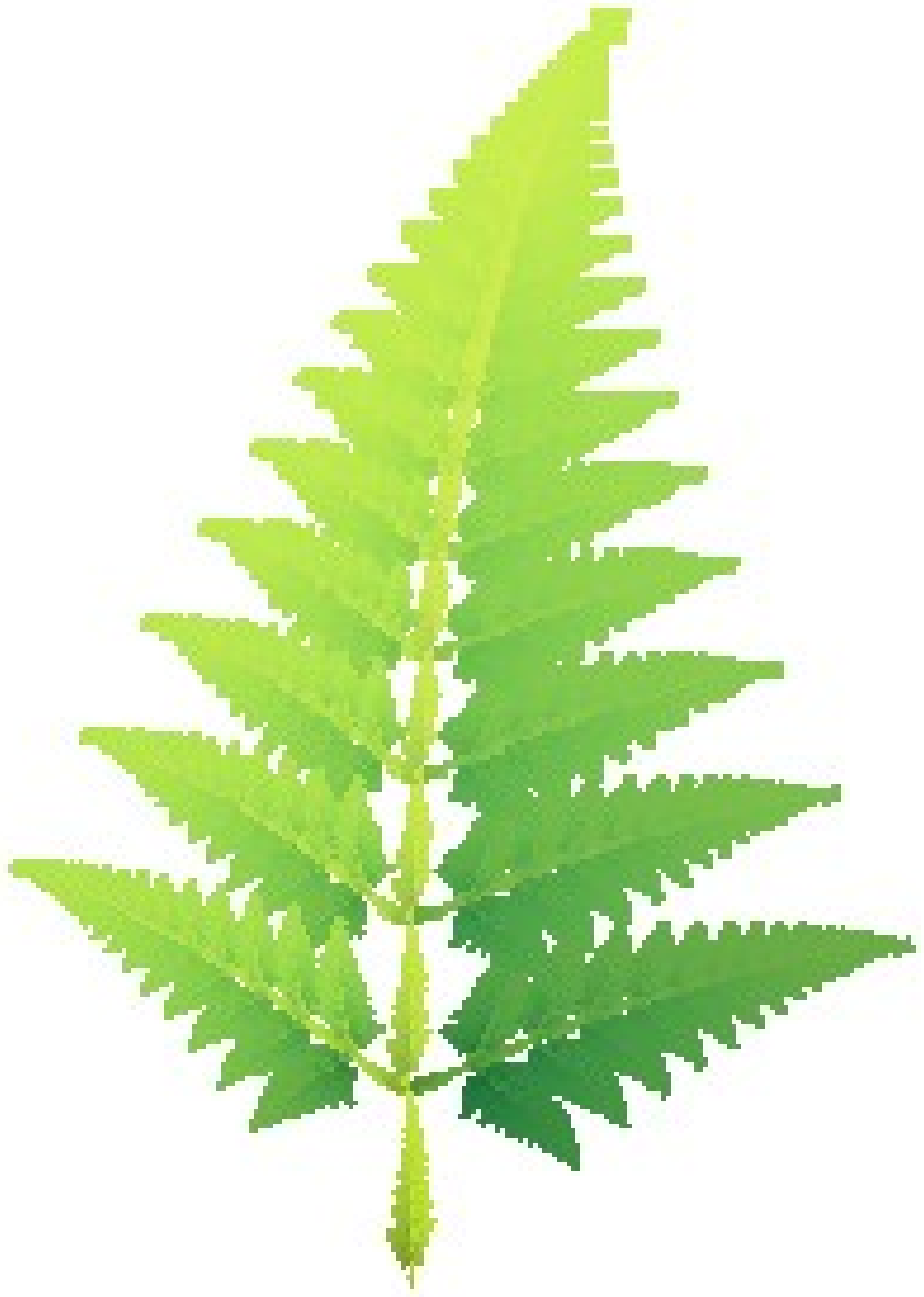}}}\FRAME{ftbpFU}{4.0283in}{3.3702in}{0pt}{\Qcb{The rendered set
attractor of the IFS $F^{2}$in Example \protect\ref{grafex}.}}{\Qlb{Fig2}}{%
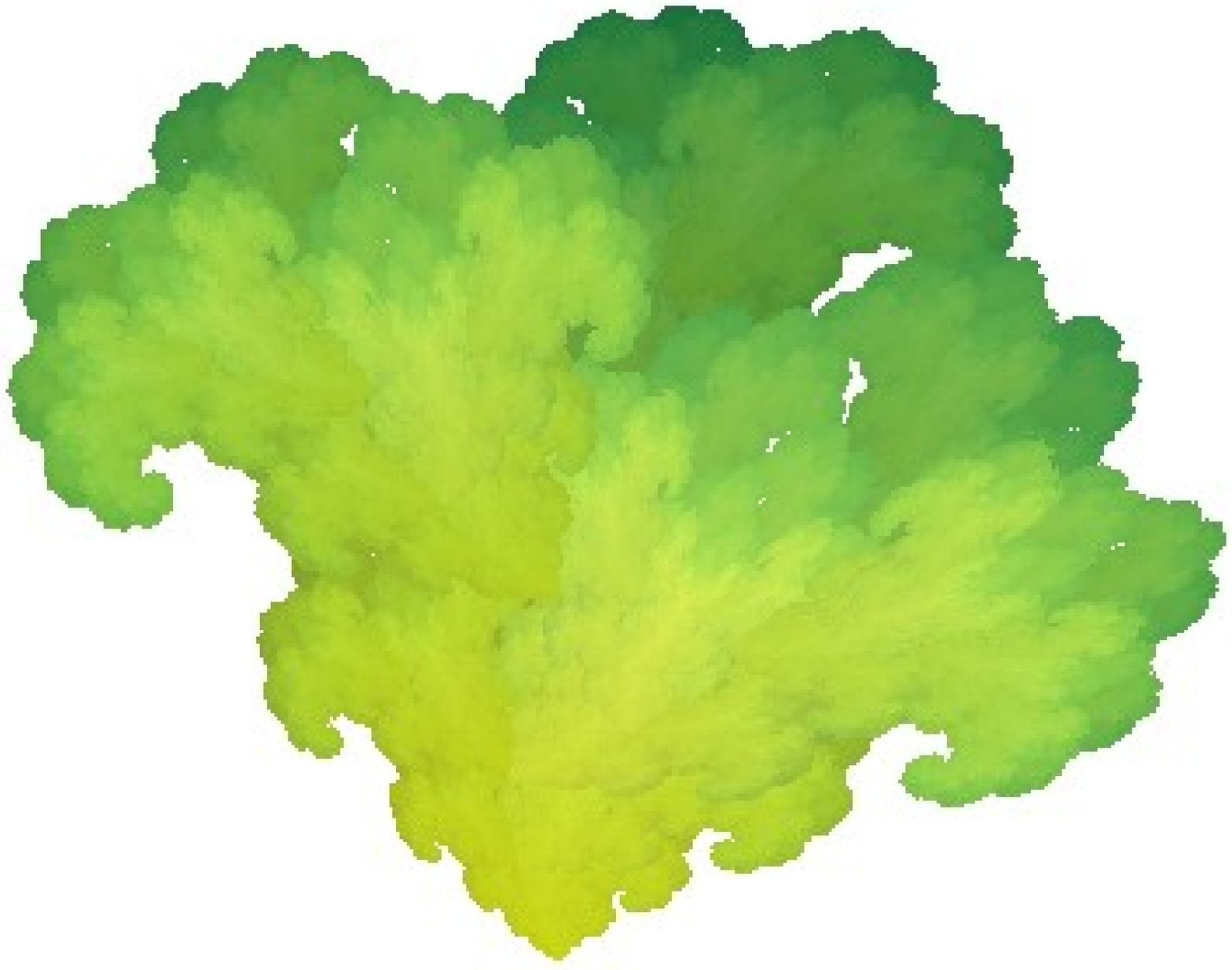}{\raisebox{-3.3702in}{\includegraphics[height=3.3702in]{soartistic.ps}}} \ 

The superIFS\ we use is 
\begin{equation*}
\mathcal{F}=\{\square ;F^{1},F^{2};P^{1}=0.5,P^{2}=0.5\}
\end{equation*}%
with $V=2.$ Then Figures \ref{Fig3} and \ref{Fig5} show two examples, from
among many different but similar ones, all equally visually complex, of
computer graphics of 2-variable fractals for this superIFS, computed using
the new random iteration algorithm. The images were rendered in much the
same way as the images of the attractor sets of $F^{1}$ and $F^{2}$ were
rendered above. The essential difference is the meaning of a
\textquotedblleft code space address\textquotedblright\ of a point on a $V$%
-variable fractal, which we define to be the sequence of lower indices of a
sequence of functions that converges to the point; for example, the point 
\begin{equation*}
\lim_{k\rightarrow \infty }f_{2}^{1}\circ f_{1}^{2}\circ f_{2}^{2}\circ
f_{1}^{1}\circ f_{2}^{2}\circ f_{2}^{1}\circ ...\circ f_{m_{k}}^{n_{k}}(x)
\end{equation*}%
corresponds to the address $212122...m_{k}...,$ in the obvious notation. 
\FRAME{ftbpFU}{4.0283in}{3.3702in}{0pt}{\Qcb{A 2-variable fractal set for
the superIFS $\mathcal{F}$ in Example \protect\ref{grafex}.}}{\Qlb{Fig3}}{%
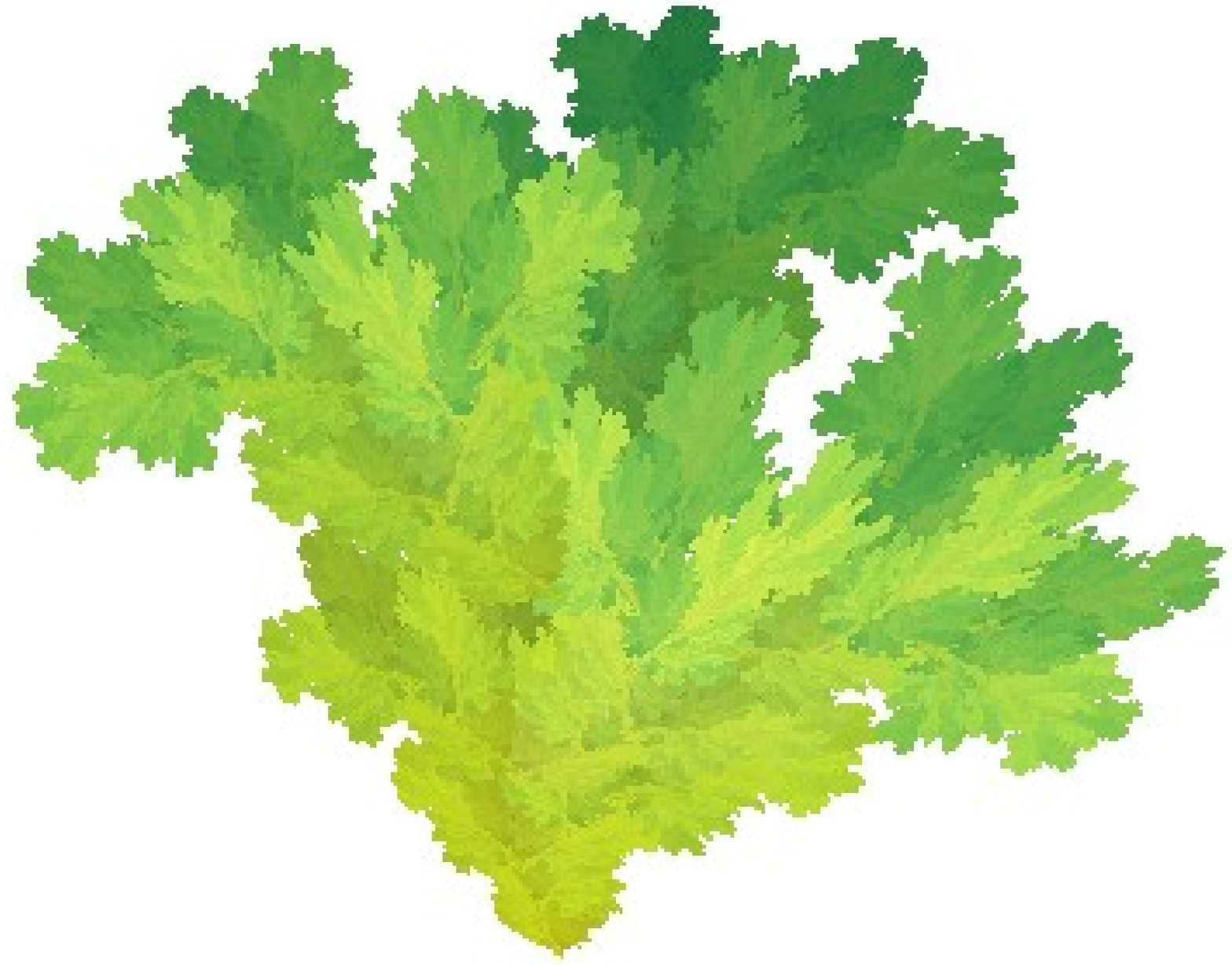}{\raisebox{-3.3702in}{\includegraphics[height=3.3702in]{leaftree11.ps}}}\FRAME{ftbpFU}{4.0283in}{3.3702in%
}{0pt}{\Qcb{Another 2-variable fractal set for the superIFS $\mathcal{F}$ in
Example \protect\ref{grafex}.}}{\Qlb{Fig5}}{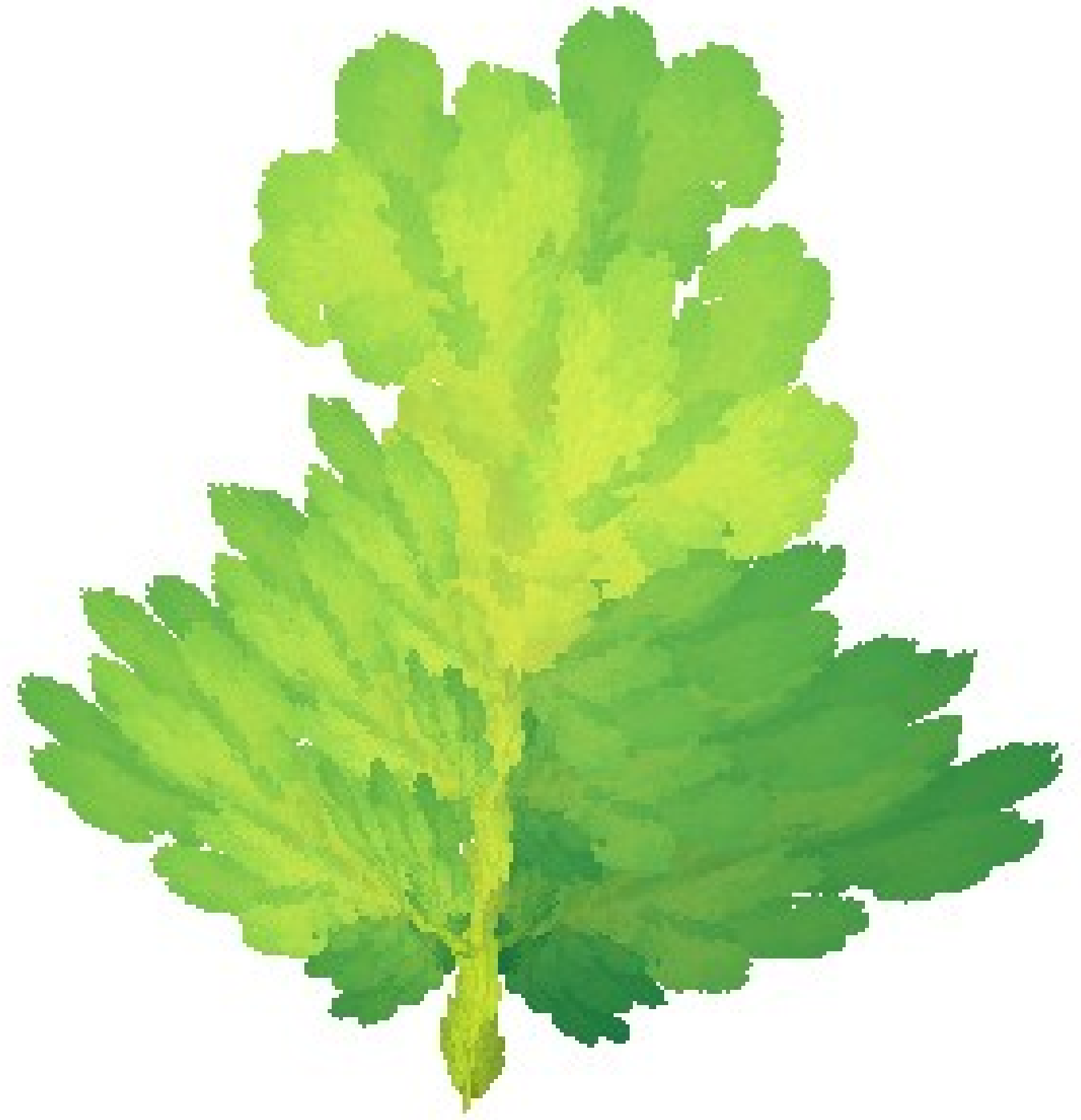}{\raisebox{-3.3702in}{\includegraphics[height=3.3702in]{leaftree15.ps}}}
\end{example}

\subsection{V-variable Fractal Interpolation}

The technique of fractal interpolation has many applications including
modelling of speech signals, altitude maps in geophysics, and stock-market
indices. A simple version of this technique is as follows. Let a set of real
interpolation points in $\{(x_{i},y_{i})\in $ $\mathbb{R}^{2}|i=0,1,...,I\}$
be given. It is desired to find a continuous function $f:[x_{0},x_{I}]%
\rightarrow $ $\mathbb{R}$ such that $f(x_{i})=y_{i}$ $\forall i\in
\{0,1,...,M\},$ such that its graph $G=\{(x,y)\in $ $\mathbb{R}^{2}:y=f(x)\}$
is a fractal, possibly with specified fractal dimension. Introduce the IFS\ 
\begin{equation*}
F=\{\mathbb{R}^{2};f_{1},f_{2},...,f_{M}\}
\end{equation*}%
with 
\begin{equation*}
f_{m}(x,y)=(a_{m}x+e_{m},c_{m}x+d_{m}y+g_{m}),
\end{equation*}%
where the real coefficients $a_{m},e_{m},c_{m},d_{m}$ and $e_{m}$ are chosen
so that 
\begin{equation*}
f_{m}(x_{0},y_{0})=y_{m-1},f_{m}(x_{0},y_{0})=y_{m},
\end{equation*}%
and $d_{m}\in \lbrack 0,1),$ for $m\in \{1,2,...,M\}.$ Then the attractor of
the IFS is the graph of a function $f$ with the desired properties, its
dimension being a function of the free parameters $\{d_{m}:m=1,2,...,M\}.$

Now let the superIFS $\mathcal{F}=\{\square
;F^{1},F^{2};P^{1}=0.5,P^{2}=0.5\}$ for some $V$, consist of two IFSs both
of which provide fractal interpolations of the data. Then all of the
elements of the corresponding superfractal will be graphs of continuous
functions that interpolate the data, have the property of $V$-variability,
and may be sampled using the random iteration algorithm.

\section{\label{seven}Generalizations}

It \ is natural to extend the notions of $V$-variable fractals, superIFS and
superfractal to include the case of maps contractive on the average, more
than a finite number of maps, more than a finite number of\ IFSs, IFSs with
a variable number of maps, IFSs operating on sets which are not necessarily
induced by point maps, other methods of constructing the probabilities for a
superIFS, probabilities that are dependent upon position etc. But for
reasons of simplicity and in order to illustrate key features we have not
treated these generalizations at any length.

\end{document}